\documentstyle[11pt]{article}
\newcommand{\df}{\displaystyle\frac}

\newcommand{\hsp}{\hspace*{\parindent}}

\newcommand{\sig}{\sigma}

\newcommand{\ra}{\rightarrow}
\newcommand{\In}{\infty}
\newcommand{\RR}{{\Bbb R}}
\newcommand{\ZZ}{{\Bbb Z}}
\newcommand{\NN}{{\Bbb N}}
\newcommand{\QQ}{{\Bbb Q}}
\newcommand{\CC}{{\Bbb C}}

\newcommand{\sF}{{\cal F}}

\newcommand{\sA}{{\cal A}}

\newcommand{\sM}{{\cal M}}
\newcommand{\sP}{{\cal P}}
\newcommand{\sB}{{\cal B}}
\newcommand{\sS}{{\cal S}}

\newcommand{\bb}{{\bf b}}

\newcommand{\bs}{{\bf s}}

\newcommand{\bv}{{\bf v}}

\input amssym.def
\input amssym.tex

\newcommand{\beql}[1]{\begin{equation}\label{#1}}
\newcommand{\eeq}{\end{equation}}

\catcode`\@=11
\renewcommand{\section}{
        \setcounter{equation}{0}
        \@startsection {section}{1}{\z@}{-3.5ex plus -1ex minus
        -.2ex}{2.3ex plus .2ex}{\large\bf}
        }
\catcode`@=12
\makeatletter
\def\eqalignno#1{\displ@y \tabskip\@centering
  \halign to\displaywidth{\hfil$\@lign\displaystyle{##}$\tabskip\z@skip
    &$\@lign\displaystyle{{}##}$\hfil\tabskip\@centering
    &\llap{$\@lign##$}\tabskip\z@skip\crcr
    #1\crcr}}
\makeatletter
\def\eqalignno#1{\displ@y \ta {\bf s} kip\@centering
  \halign to\displaywidth{\hfil$\@lign\displaystyle{##}$\ta {\bf s} kip\z@skip
  \ldots \@lign\displaystyle{{}##}$\hfil\ta {\bf s} kip\@centering
  \ldots llap{$\@lign##$}\ta {\bf s} kip\z@skip\crcr
    #1\crcr}}
\makeatother

\makeatletter
\def\@sect#1#2#3#4#5#6[#7]#8{\ifnum #2>\c@secnumdepth
     \def\@svsec{}\else 
     \refstepcounter{#1}\edef\@svsec{\csname the#1\endcsname.\hskip .75em }\fi
     \@tempskipa #5\relax
      \ifdim \@tempskipa>\z@ 
        \begingroup #6\relax
          \@hangfrom{\hskip #3\relax\@svsec}{\interlinepenalty \@M #8\par}%
        \endgroup
       \csname #1mark\endcsname{#7}\addcontentsline
         {toc}{#1}{\ifnum #2>\c@secnumdepth \else
                      \protect\numberline{\csname the#1\endcsname}\fi
                    #7}\else
        \def\@svsechd{#6\hskip #3\@svsec #8\csname #1mark\endcsname
                      {#7}\addcontentsline
                           {toc}{#1}{\ifnum #2>\c@secnumdepth \else
                             \protect\numberline{\csname the#1\endcsname}\fi
                       #7}}\fi
     \@xsect{#5}}
\def\@begintheorem#1#2{\it \trivlist \item[\hskip \labelsep{\bf #1\ #2.}]}
\makeatother

\begin{document}
\begin{center}
{\Large {\bf The $3x + 1$ Problem: An Annotated Bibliography (1963--1999)\\
(Sorted by Author) }} \\
\vspace{\baselineskip}
{\large {\em Jeffrey C. Lagarias}} \\
\vspace*{.2\baselineskip}
Department of Mathematics \\
University of Michigan \\
Ann Arbor, MI 48109--1109\\
{\tt lagarias@umich.edu}\\
\vspace*{1\baselineskip}
(January 1, 2011 version) \\
\vspace*{2\baselineskip}
\end{center}

\noindent{\bf ABSTRACT.} 
The $3x+1$ problem concerns iteration of the map $T: \ZZ \rightarrow \ZZ$
given by
$$
T(x) = \left\{
\begin{array}{cl}
\df{3x+1}{2} & \mbox{if} ~~x \equiv 1~~ (\bmod~2) ~. \\
~~~ \\
\df{x}{2} & \mbox{if}  ~~x \equiv 0~~ (\bmod~2)~.
\end{array}
\right.
$$
The $3x+1$ Conjecture asserts that each $m \geq 1$ has some iterate
$T^{(k)} (m) = 1$.
This is an annotated bibliography of work done on the $3x+1$ problem
and related problems from 1963 through 1999.
At present the  $3x+1$ Conjecture remains unsolved.\medskip

This version of the $3x+1$
bibliography  sorts the papers alphabetically by the first author's surname. 
A different version of this bibliography appears in the recent book: 
 {\em The Ultimate Challenge: The $3x+1$ Problem}, Amer. Math. Soc, Providence, RI  2010.
 In it the papers are sorted by decade (1960-1969), (1970-1979), (1980-1989), (1990-1999).\\

%****************************************************
%
% Section 1. Terminology
%
%****************************************************
%
\setlength{\baselineskip}{1.0\baselineskip}
\section{Introduction}
\hsp

The $3x+1$ problem is most simply stated in terms
of the {\em Collatz function} $C(x)$ defined on integers as
``multiply by three and add one'' for odd integers
and ``divide by two'' for even integers. That is,
$$
C(x) = 
\left\{
\begin{array}{cl}
3x+1 & \mbox{if}~ x \equiv 1~~ (\bmod ~2 ) ~, \\
~~~ \\
\df{x}{2} & \mbox{if} ~~x \equiv 0~~ (\bmod~2) ~,
\end{array}
\right.
$$
The 
{\em $3x+1$ problem}, or {\em Collatz problem}, 
is to prove that starting from any positive integer,
some iterate of this function takes the value $1$.

Much work on the problem is stated in terms 
of the  {\em $3x+1$  function} 
$$
T(x) = \left\{
\begin{array}{cl}
\df{3x+1}{2}  & \mbox{if} ~~x \equiv 1~~ (\bmod ~ 2) \\
~~~ \\
\df{x}{2} & \mbox{if}~~ x \equiv 0~~ (\bmod ~2 )~.
\end{array}
\right.
$$
The  $3x  +  1$ {\em conjecture} 
states that every $m \geq  1$ has some iterate
$T^{(k)} (m) = 1$.

The $3x+1$ problem
is generally attributed to Lothar Collatz and 
has reportedly circulated since at least the early 1950's.
It also goes under the names {\em Syracuse Problem}, {\em Hasse's Algorithm},
 {\em Kakutani's Problem}
and {\em Ulam's Problem}.  The first published reference
to a $3x+1$-like function that I am aware of is Klamkin (1963),
which concern's Collatz's original function, which is a permutation
of the integers. Collatz has stated that he studied this function 
in the 1930's, see  Lagarias (1985).  The $3x+1$ problem itself
was reportedly described in an informal lecture of Collatz in 1950 at the
International Math. Congress in Cambridge, Massachusetts
(see Shanks (1965) and  Trigg et al (1976)).
The first journal appearance of the $3x+1$
problem itself seems to be the text of  a  1970 lecture of
H. S. M. Coxeter, which appeared in  Coxeter (1971). This was followed by 
Beeler et al (1972), Conway (1972),  Gardner (1972), Kay (1972)
and Ogilvy (1972). See Bryan Thwaites (1985) for his
assertion  to have formulated the problem in 1952.
See Collatz (1986) for his  assertions on  formulating the $3x+1$ problem
prior to 1952.

The $3x+1$ problem can also be rephrased as a  
problem concerning sets of integers generated using affine maps. 
Let $T$ be the smallest set of integers
including $1$ and closed under iteration of the affine maps
$x \mapsto 2x$ and $3x+2 \mapsto 2x+1$. Here the latter map
is the affine map $y \to \frac{2y-1}{3}$,  with input restricted 
  to integers  $y\equiv 2~(\bmod ~3)$, so that the output
 is an integer.  The $3x+1$ conjecture 
asserts that $T$ is  the set of all positive integers. Therefore this 
bibliography includes work on sets of integers generated by
iteration of affine maps, tracing back to Isard and Zwicky (1970)
and Klarner and Rado (1974), which includes  a problem of Erd\H{o}s,
described in Klarner (1982).

As of 1999 the  $3x+1$ conjecture was verified
 up to $2 \times 10^{16}$ by a  computation 
of  Oliveira e Silva (1999).  It has since been
verified to at least  $ 2 \times 10^{18}$ 
in independent computations by
Oliveira e Silva and Roosendaal.
At present the $3x+1$ conjecture  remains unsolved. \\

This bibliography includes some surveys  on results on the $3x+1$ problem  during 1963--1999:~
Lagarias (1985), M\"{u}ller (1991), and 
the first chapter of Wirsching (1998a).
%A more recent update  appears in Chamberland (2002). 

%****************************************************
%
% Section 2. Terminology
%
%****************************************************
%
\section{Terminology}
\hsp

We  use the following definitions.
The {\em trajectory} or {\em forward orbit} of an integer $m$ is the set
$$
\{ m,~ T(m)~, ~~ T^{(2)} (m) , \ldots \}~.
$$
The
{\em stopping time} $\sigma(m)$
of $m$ is the least $k$ such that $T^{(k)} (m) < m$,
and is $\infty$ if no such $k$ exists.
The
{\em total stopping time }
$\sig_\infty (m)$ is the least $k$ such that
$m$ iterates to $1$ under $k$ applications of the function  $T$ i.e.
$$ \sig_\infty (m) : = \inf~ \{k~:~ T^{(k)} (m) = 1\}. $$ 
The {\em scaled total stopping time} or {\em gamma value} $\gamma (m)$ is
$$ \gamma (m) := \frac {\sig_\infty (m)} {log~ m} $$
The
{\em height} $h(m)$  the least $k$ for which the Collatz function $C(x)$ 
has $C^{(k)} (m) =  1$. It is easy to show that 
$$h(m) = \sig_\infty (m) + d(m),$$ 
where $d(m)$ counts the
number of iterates $T^{(k)} (m)$ in 
$0 \leq k < \sigma_\infty (m)$ that are odd integers. 
%If the orbit of $n$ is a cycle, then $d(n)$ counts  the number of
%odd integers in the cycle, which equals the number of 
% {\em turning points} in the cycle, i. 
Finally, the function $\pi_a (x)$ counts the number of $n$ with
$|n| \leq x$ whose forward orbit under $T$ includes
$a$.

%****************************************************
%
% Section 3. Bibliography
%
%****************************************************
%
\newpage
\section{Bibliography}
\hsp

This annotated bibliography mainly covers research 
articles and survey articles on the $3x+1$ problem and
related problems. It provides additional 
references to earlier history, much of it appearing as problems.
  It also includes  a few influential technical
reports that were never published. Chinese references 
are written with surnames first. 

\begin{enumerate}
\item
Sergio Albeverio, Danilo Merlini and Remiglio Tartini (1989), 
{\em Una breve introduzione a diffusioni su insiemi frattali e ad
alcuni essempi di sistemi dinamici semplici,}
Note di matematica e fisica,
Edizioni Cerfim Locarno {\bf 3} (1989), 1--39. \\
\newline
\hspace*{.25in}
This paper discusses dimensions of some simple fractals,
starting with the Sierpinski gasket and  the Koch snowflake.
These arise as a fixed set from combining several linear iterations. 
It then considers the Collatz function  iteration as analogous, as it is given by  set of
two linear transformations.  It looks at the tree of inverse iterates ("chalice")
and estimates emprically the number of leaves at depth at most
$k$ as growing like $c^k$ for $c \approx 1.265$. It discusses various
cascades of points that arrive at the cycle $\{ 1, 4,2\}$.

It then looks at iteration schemes of type: start with vertices of equilateral triangle,
and a new point $x_0$ not on the triangle. Then form
$f(x_0) = \frac{x_0 + x_k}{\xi}$ where $\xi >0$ is a  parameter. then iterate
with $k$ varying in any order through $1,2,3$. Depending on the value
of $\xi$ for $1 < \xi < 2$ get a fractal set, for $\xi=2$ get an orbit supported
on Sierpinski gasket. It discusses critical exponents, analogy with Julia sets,
and then the Ising model.

\item
Jean-Paul Allouche (1979),
{\em Sur la conjecture de ``Syracuse-Kakutani-Collatz''},
S\'{e}minaire de Th\'{e}orie des Nombres 1978--1979, Expose
No. 9, 15pp., CNRS Talence (France), 1979.
(MR 81g:10014). \\
\newline 
\hspace*{.25in}
This paper studies generalized $3x+1$ functions of the
form proposed by Hasse. These have the form
$$
T(n) = T_{m, d, R}(n) := \left\{
\begin{array}{cl}
\df{mn+r_j}{d}  & \mbox{if} ~~n \equiv j~~ (\bmod ~ d),~1 \le j \le d-1 \\
~~~ \\
\df{x}{d} & \mbox{if}~~ n \equiv 0~~ (\bmod ~d )~.
\end{array}
\right.
$$
in which the parameters $(d, m)$ satisfy
$d \ge 2$, $gcd(m, d) =1$,
and  the set $R= \{ r_j: 1 \le j \le d -1\}$ has each
$r_j \equiv -mj~(\bmod~d)$.  The author  notes that
it is easy to show  that all maps in Hasse's class with $1 \le m < d$
have a finite number of cycles, and for these maps
all orbits eventually enter one of these cycles.
Thus we may assume that $m > d$.

The paper proves two theorems. Theorem 1 
improves on  the results of Heppner (1978). Let $T(\cdot)$
be a function in Hasse's class with parameters $d, m$.
Let $a > 1$ be fixed and set $k = \lfloor \frac{\log x}{a \log m} \rfloor.$
Let $A, B$ be two rationals with $A < B$ that are 
not of the form $\frac{m^i}{d^j}$ for 
integers $i \ge 0, j \ge 1$,
and consider the counting function
$$
F_{a, A, B}(x) := \sum_{n=1}^{x} \chi_{(A,B)}\left( \frac{ T^{(k)}(n)}{n}\right),
$$
in which $\chi_{(A, B)}(u) = 1$ if $A < u < B$ and $0$ otherwise.
Let $C$ be the maximum of the denominators of $A$ and $B$. Then:

(i) If $A < B < 0$ then 
$$
F_{a, A, B}(x) = O\left( C x^{\frac{1}{a}} \right),
$$
where the constant implied by  the $O$-symbol is absolute.

(ii) If $B > 0$ and there exists $\epsilon > 0$ such that
$$
\frac{d-1}{d} \ge  \frac{\log d}{\log m} + \frac{\log B}{k \log m} + \epsilon,
$$
then for $k > k_0(\epsilon)$, there holds 
$$
F_{a, A, B}(x) = O \left( C x^{\frac{1}{a}} + 
x^{1 - \frac{|\log \eta|}{a \log m}} \right)
$$
for some $\eta$ with $0 < \eta < 1$ which depends
on $\epsilon.$  This is true in particular
when $m > d^{\frac{d}{d-1}}$ with $B$ fixed and $x \to \infty$.

(iii) If $B > A > 0$ and if there exists $\epsilon$ such that
$$
 \frac{d-1}{d} \le \frac{\log d}{\log m} +
\frac{\log A}{k \log m} -  \epsilon,
$$
then for $k > k_0(\epsilon)$, there holds
$$
F_{a, A, B}(x) = O \left( Cx^{\frac{1}{a}} + 
x^{1 - \frac{|\log \eta|}{a \log m}} \right)
$$
for some $\eta$ with $0 < \eta < 1$ which depends
on $\epsilon.$ This is true in particular
when  $m < d^{\frac{d}{d-1}}$, with $A$ fixed and $x \to \infty$.

Theorem 2  constructs for given values 
$d, m$ with $gcd(m, d) =1$   with $m > d \ge 2$ 
two functions $F(\cdot)$
and $G(\cdot)$ which fall slightly outside Hasse's
class and have the following properties:

(i) The function $F(\cdot)$ has a finite number of
periodic orbits, and every $n$ when iterated under $F(\cdot)$ 
eventually enters one of these orbits.

(ii) Each orbit of the  function $G(\cdot)$ is divergent, i.e.
$|G^{(k)}(n)| \to \infty$ as $k \to \infty$, for all
but a finite number of initial values $n$.

To define the first function $F(\cdot)$ pick an integer $u$ such that
$d <m < d^u$, and set 
$$
F(x) = \left\{
\begin{array}{cl}
\df{mx + r_j}{d}  & \mbox{if} ~~n \equiv j  (\bmod ~ d^u),~ 
gcd( j, d) = 1, \\
~~~ \\
\df{x}{d} + s_j & \mbox{if}~~ x \equiv j~~ (\bmod ~d^u),~~j 
\equiv 0 (\bmod~d).
\end{array}
\right.
$$
in which $0 < r_j < d^u$ is determined by 
 $mj + r_j \equiv 0 (\bmod~ d^u)$,
and $0 \le s_j < d^{u-1}$ is determined by
$\frac{j}{d} +  r_j \equiv 0 (\bmod~ d^{u-1})$.
The second function $G(x)$ is defined by $G(x) = F(x) + 1$.

These functions fall outside
Hasse's class because each 
 is linear on  residue classes $n (\bmod~d^u)$
for some $u \ge 2$, rather than linear on residue classes $(\bmod~d)$. 
However both  these functions exhibit behavior qualitatively like functions
in Hasse's class: There is a 
constant $C$ such that 
$$
|F(n) - \frac{n}{d}| \le C ~~\mbox{~if~}~~ n \equiv 0 ~(\bmod~d).
$$
$$
|F(n) - \frac{mn}{d}| \le C ~~\mbox{~if~}~~ n \not\equiv 0 ~(\bmod~d).
$$
and similarly for $G(n)$, taking $C= d^{u-1} + 1$.
The important difference is that functions  in Hasse's class are mixing on 
residue classes $(\bmod~d^k)$ for all powers of $k$, while
the functions $F(\cdot)$ and $G(\cdot)$ are not mixing in
this fashion. The nature of the non-mixing behaviors
of these functions underlies the proofs
of properties (i),  resp. (ii) for $F(\cdot)$, resp. $G(\cdot)$.

%\newpage
\item
Amal S.  Amleh, Edward A. Grove, Candace M. Kent, and Gerasimos Ladas (1998),
{\em On some difference equations with eventually periodic solutions,}
J. Math. Anal. Appl. {\bf 223} (1998), 196--215. (MR 99f:39002)
\newline 
\hspace*{.25in}

The authors study  the boundedness and periodicity of solutions
of the set of difference equations
$$
x_{n+1} = \left\{\begin{array}{cl}
\frac{1}{2}(\alpha x_n + \beta x_{n-1}) & \mbox{if}~~x \equiv 0~~(\bmod~2 ) ~.\\
~~~ \\
\gamma x_n + \delta x_{n-1}& \mbox{if}~~x \equiv 1~~(\bmod~2 ) ~, 
\end{array}
\right.
$$
where the parameters $\alpha, \beta, \gamma, \delta \in \{ -1, 1\}$, and the
initial conditions $(x_0, x_1)$ are integers. There are $16$ possible such iterations. 
Earlier Clark and Lewis (1995) considered the case 
$(\alpha, \beta, \gamma, \delta)=(1, 1, 1, -1)$, and showed that all orbits 
with initial conditions integers $(x_0, x_1)$ with $gcd(x_0, x_1)=1$ converge
to one of three periodic orbits .
Here the authors consider  all $16$ cases, showing first a duality between solutions of
$(\alpha, \beta, \gamma, \delta)$ and $(-\alpha, \beta, -\gamma, \delta)$, 
taking a solution $\{ x_n\}$ of one to $\{ (-1)^{n+1} x_n\}$ of the other.
This reduces to considering the eight cases with 
$\alpha=+1$. They resolve six of these cases, as follows, leaving open the cases
$(\alpha, \beta, \gamma, \delta)= (1, -1, 1,1)$
and $(1, -1, -1, -1)$.

For the parameters $(1,1,1,1)$ they show all orbits are eventually constant or unbounded,
and that unbounded orbits occur. 

For the  parameters $(1,1,1,-1)$, the work of  Clark and Lewis (1995) showed all orbits
are eventually periodic.

For the parameters $(1,1,-1, 1)$ all solutions are eventually periodic, and there are
five relatively prime cycles, the fixed points $(1), (-1)$, the $4$-cycles $(2, -1, 3,1),~(-2,1,-3,1)$
and the $6$-cycle $(1,0,1,-1,0,-1)$.

For the parameters $(1,1,-1, -1)$ all solutions are eventually periodic, and there are
four relatively prime cycles,  the fixpoints $(1), (-1)$ and the $3$-cycles $ (-1,0,1),~(1,0,-1)$.

For parameters $(1, -1, 1,-1)$, all orbits are eventually periodic, with two relatively prime cycles,
 the $6$-cycle $(-1, 0,1,1, 0, -1)$.

For parameters $(1, -1,-1, 1)$ all orbits are eventually periodic, with one relatively prime cycle, the
$8$-cycle $(0, -1, 1, 1, 0, 1, -1,-1)$.

The authors conjecture that for the two remaining cases  $(1, -1, 1,1)$
and $(1, -1, -1, -1)$, that
all orbits are eventually periodic.

\item
S. Anderson (1987), 
{\em Struggling with the $3x  +  1$ problem,} 
Math. Gazette {\bf 71} (1987), 271--274. \\
\newline 
\hspace*{.25in}
This paper studies simple analogues of the $3x  +  1$ function such as
$$
g(x) = \left\{
\begin{array}{cl}
x+k & \mbox{if}~~x \equiv 1~~(\bmod~2 ) ~, \\
~~~ \\
\df{x}{2} & \mbox{if}~~x \equiv 0~~(\bmod~2 ) ~.
\end{array}
\right.
$$
For $k =  1$ when iterated this map gives the binary expansion of $x$.
The paper also reformulates the $3x  +  1$ Conjecture using the function:
$$
f(x) =  \left\{
\begin{array}{cl}
\df{x}{3} & \mbox{if}~~x \equiv 0~~(\bmod~3) ~, \\
~~~ \\
\df{x}{2} & \mbox{if}~~x \equiv 2 ~\mbox{or}~4~ (\bmod ~6) ~, \\
~~~ \\
3x+1 & \mbox{if}~~x \equiv 1~ \mbox{or}~ (\bmod ~6) ~.
\end{array}
\right.
$$

\item
Stefan Andrei and Cristian  Masalagiu (1998),
{\em About the Collatz Conjecture,}
Acta Informatica {\bf 35} (1998), 167--179. (MR 99d:68097). \\
\newline 
\hspace*{.25in}
This paper describes two recursive algorithms for computing 
$3x+1$-trees, starting from a given base node. 
A $3x+1$-tree is a tree of
inverse iterates of the function $T(.)$. The second algorithm
shows a speedup of a factor of about three over the ``naive' first 
algorithm.

\item
David Applegate and Jeffrey C. Lagarias (1995a),
{\em Density Bounds for the $3x+1$ Problem I. Tree-Search Method,}
Math. Comp., {\bf 64} (1995), 411--426. 
(MR 95c:11024) \\
\newline 
\hspace*{.25in}
Let $ n_k (a)$ count the number
of integers $n$ having $T^{(k)} (n) =  a$.
Then for any $a \not\equiv 0~~(\bmod~3)$ and
sufficiently large $k$,
$(1.299)^k \leq n_k (a) \leq (1.361)^k$.
Let $\pi_k (a)$ count the number of $| n | \leq x$ which eventually
reach $a$ under iteration by $T$.
If $a \not\equiv 0~~(\bmod~3)$
then $\pi_a (x)  >  x^{.643}$ for all sufficiently large $x$.
The extremal distribution of number of leaves in $3x+1$ trees with
root $a$ and depth $k$ (under iteration of $T^{-1}$)
as $a$ varies
are computed for $k \leq 30$.
The proofs are computer-intensive.

\item
David Applegate and Jeffrey C. Lagarias (1995b), 
{\em Density Bounds for the $3x+1$ Problem II.\ Krasikov Inequalities,}
Math. Comp., {\bf 64} (1995). 427--438.
(MR 95c:11025) \\
\newline 
\hspace*{.25in}
Let $\pi_a (x)$ count the number of $|n| \leq x$ 
which eventually
reach $a$ under iteration by $T$.
If $a \not\equiv 0~~(\bmod~3)$, then
$\pi_a (x)  >  x^{.809}$
for all sufficiently large $x$.
It is shown that the inequalities of Krasikov (1989) can be used to
construct nonlinear programming problems which yield lower bounds for the
exponent
$\gamma$ in $\pi_a (x)  > x^\gamma$.
The exponent above was derived by computer for such a nonlinear program
having about 20000 variables.

\item
David Applegate and Jeffrey C. Lagarias (1995c), 
{\em On the distribution of $3x+1$ trees, }
Experimental Mathematics
{\bf 4} (1995), 101--117.
(MR 97e:11033). \\
\newline 
\hspace*{.25in}
The extremal distribution of the number of leaves in $3x+1$ 
trees with root $a$ and
depth $k$ (under iteration of $T^{-1}$) as $a$ varies were computed for
$k \leq 30$ in Applegate and Lagarias (1995a).
These data appear to have a much narrower spread around the mean value
$( \frac{4}{3} )^k$ of leaves in a $3x+1$ tree of depth $k$ than is 
predicted by
(repeated draws from) the branching process models of 
Lagarias and Weiss (1992).
Rigorous asymptotic results are given for the branching process models.
\newline
\hspace*{.25in}
The paper also derives formulas for the expected number of 
leaves in a
$3x+1$ tree of depth $k$ whose root node is $a$ $( \bmod ~3^\ell )$.
A 3-adic limit is proved to exist almost everywhere
as $k \rightarrow \infty$, the
expected number of leaves being $W_\infty (a) \left( \frac{4}{3} \right)^k$
where the function $W_\infty : \ZZ_3^\times \rightarrow \RR$
 almost everywhere
satisies the 3-adic functional equation
$$
W_\infty( \alpha ) = \df{3}{4} \left( W_\infty ( 2 \alpha ) + 
\psi ( \alpha \bmod ~9) W_\infty \left(
\df{2 \alpha -1}{3} \right) \right) ~, \eqno{\mbox{($\ast$)}}
$$
in which $\psi ( \alpha ) = 1$ if $\alpha \equiv 2$ or $8 ( \bmod ~9 )$
and is 0 otherwise.
(Here $\ZZ_3^\ast = \{ \alpha \in \ZZ_3 : 
\alpha \not\equiv 0~ ( \bmod ~3) \}$).
It is conjectured that $W_\infty$ is continuous and everywhere nonzero.
It is an open problem to characterize solutions of the 
functional equation $( \ast )$.

\item
Jacques Arsac (1986), 
{\em Algorithmes pour v\'erifier la conjecture de Syracuse,}
C. R. Acad. Sci. Paris
{\bf 303},Serie I, no. 4, (1986), 155--159.
[Also: RAIRO, Inf. Th\'{e}or. Appl. {\bf 21} (1987), 3--9.]
(MR 87m:11128). \\
\newline 
\hspace*{.25in}
This paper studies 
 the computational complexity of algorithms 
to compute stopping times
of $3x  +  1$ function on all integers below a given bound $x$.

\item
Charles Ashbacher (1992), 
{\em Further Investigations of the Wondrous Numbers,}
J. Recreational Math. {\bf 24} (1992), 1--15. \\
\newline 
\hspace*{.25in}
This paper numerically studies the ``MU'' functions $F_D(x)$ 
of Wiggin (1988) on
$x \in \ZZ^+ $ for $2 \leq D \leq 12$.
It finds no exceptions for Wiggin's conjecture that all cycles of
$F_D $ on $\ZZ^+ $ contain an integer smaller than $D$,
for $x < 1.4 \times 10^7 $.
It tabulates integers in this range that have a large stopping time, 
and observes
various patterns.
These are easily explained by observing that,
for $n \not\equiv 0~~(\bmod ~D)$,
$F_D^{(2)} (n) = \frac{n(D+1)-R}{D}$
if
$n  \equiv  R~(\bmod ~D)$, $-1 \leq R \leq D -2$,
hence, for most $n$, $F_D^{(2)} (n)  >  n$, although $F_D $
decreases iterates on the average.

\item
Arthur Oliver  Lonsdale Atkin (1966), 
{\em Comment on Problem $63-13^{*}$},
SIAM Review {\bf 8} (1966), 234--236. \\
\newline 
\hspace*{.25in}
This comment gives more information
of the  problem of Klamkin (1963)
concerning iteration of the original
Collatz function, which is a permutation
of the positive integers. 
Adding to the  comment of Shanks (1965),
he notes there is a method which in
principle can determine all the cycles
of a given period $p$ of this map. 
This method determines upper and lower bounds on the
integers that can appear in such a cycle.
By computer calculation he shows that
aside from the known cycles of periods
$1$, $2$, $5$ and $12$ on the nonnegative
integers, there are no other cycles of
period less than $200$. 

He gives an example casting some doubt on  the heuristic
of Shanks (1965) concerning the  possible
lengths of periods. He shows that for the 
related permutation
$f(3n) = 4n+3, f(3n+1) = 2n, f(3n+2) = 4n+1$,
which should obey a similar heuristic, 
that it has  a cycle of period $94$ ( least term $n=140$),
and  $94$
is not a denominator of the continued
fraction convergent to $\log_2 3$.

Atkin  presents a heuristic argument asserting
that the Collatz permutation  should
only have a finite number of cycles,
since the iterates grow ``on average''
at an exponential rate.

\item
Michael R. Avidon (1997), 
{\em On primitive $3$-smooth partitions of $n$},
Electronic J. Combinatorics {\bf 4} (1997), no.1 , 10pp.
(MR 98a:11136). \\
\newline 
\hspace*{.25in}
The author studies the number $r(n)$ of representations of $n$ as
sums of numbers of the form $2^a 3^b$ which are
primitive (no summand divides another). Iterates of 
$3x+1$ function applied to $n$ that get to $1$
produces a representation of $n$ of this kind.
The author proves results 
about the maximal and average order of this function. 
See also Blecksmith, McCallum and Selfredge (1998) for more
information.

\item
Claudio Baiocchi (1998),
{\em 3N+1, UTM e Tag-Systems} (Italian),
Dipartimento di Matematica dell'Universit\`{a} "La Sapienza"
di Roma {\bf 98/38} (1998). \\
\newline
\hspace*{.25in}
This technical report constructs small state Turing machines that 
simulate the $3x+1$ problem. Let $T(k, l)$ denote the class of
one-tape Turing machines with $k$ state, with $l$-symbols,
with one read head, and the tape is infinite in two directions.
The author constructs Turing machines for simulating 
the $3x+1$ iteration in the classes
$T(10,2), T(5,3), T(4,4), T(3,5)$ and $T(2,8)$, working on unary inputs.
It follows that  no method is curently known to decide the reachability  problem for such
machines. The author then produces a universal Turing machine in the class 
$T(22, 2)$.

{\em Note.} This work was motivated by  a conference paper of  M. Margenstern (1998),
whose journal version is: M. Margenstern,  Theor. Comp. Sci. {\bf 231} (2000), 217-251.

\item
Ranan B. Banerji (1996),
{\em Some Properties of the $3n+1$ Function,}
Cybernetics and Systems
{\bf 27} (1996), 473--486. \\
\newline 
\hspace*{.25in}
The paper derives elementary results on forward iterates
of the $3x+1$ function viewed as binary integers, and on
backward iterates of the map $g$ taking odd integers
to odd integers, given by
$$
g(n) : = \df{3n+1}{2^k} ~,~~\mbox{where~~ $2^k || 3n+1$}~.
$$
Integers $n \equiv 0 (\bmod~3)$ have no
preimages under $g$.
If $n \not\equiv 0 ( \bmod~3)$ define $g^{-1} (n)$ to be the
unique integer $t$ such that $g(t) = n$ and
$t \not\equiv 5 ( \bmod~8)$.
Note that
each odd $n$ there are infinitely many $\tilde{t}$
with $g( \tilde{t}) = n$.
If $d(n) = 4n+1$, these preimages are just 
$\{d^{(j)} g^{-1} (n) : j \geq 1 \}$.
The ternary expansion of $g^{-1} (n)$ is asserted to be computable from
the ternary expansion of $n$ by a finite automaton.
The author conjectures that given any odd integer $n$, there is some
finite $k$ such that $(g^{-1} )^{(k)} (n) \equiv 0 ( \bmod~3)$.
Here we are iterating the partially defined map
\begin{eqnarray*}
g^{-1} (6n+1) & = & 8n+1 ~, \\
g^{-1} (6n+5) & = & 4n+3~, 
\end{eqnarray*}
and asking if some
iterate is $0 ( \bmod~3)$.
The problem resembles Mahler's $Z$-number iteration
[J. Australian Math. Soc. {\bf 8} (1968), 313--321].

%\newpage
\item
Enzo Barone (1999),
{\em A heuristic probabilistic argument for the Collatz sequence.}
(Italian), Ital. J. Pure Appl. Math. {\bf 4} (1999), 151--153.
(MR 2000d:11033). \\
\newline 
\hspace*{.25in}
This paper presents a  heuristic probabilistic argument
which argues that iterates of the $3x+1$-function should decrease
on average by a multiplicative
factor $(\frac{3}{4})^{1/2}$ at each step. Similar arguments
appear earlier in Lagarias (1985) and many other places,
and trace back to the original work of Terras, Everett and Crandall.

\item
Michael Beeler, William Gosper and Richard Schroeppel (1972),
{\em HAKMEM,}  Memo 239, Artificial Intelligence Laboratory,
MIT, 1972. (ONR Contract N00014-70-A-0362-0002).\\
\newline 
\hspace*{.25in}
The influential memorandum, never published in a journal,  is a collection of
problems and results. The list 
contains solved and unsolved problems, computer programs to write,
programming hacks,  computer hardware to design. There are 191 items in all.
Example:

"{\bf Problem 95: } Solve {\em chess.} There are about $10^{40}$ possible positions;
in most of them, one side is hopelessly lost."

It contains one of the earliest statements of the $3x+1$ problem, which
 appears as item 133. It was
contributed by A. I. Lab members
Richard Schroeppel, William Gosper, William Henneman and Roger Banks.
It  asks if there are any other cycles on the integers other than the five known ones. It asks if
any orbit diverges. 

[This memo is currently available online at: \\
{\tt http://www.inwap.com/pdp10/hbaker/hakmen/hakmem.html}]

\item
Edward Belaga (1998),
{\em Reflecting on the $3x+1$ Mystery: Outline of a Scenario-
Improbable or Realistic?}
U. Strasbourg report 1998-49, 10 pages.\\
({\tt http://hal.archives-ouvertes.fr/IRMA-ACF}, file hal-0012576)\\
%http://www-irma.u-strasbg.fr/irma/publications/1998/98049.shtml \\
\newline 
\hspace*{.25in}
This is an expository paper, discussing the possibility
that the $3x+1$ conjecture is an undecidable problem.
Various known results, pro and con, are presented.

\item
Edward Belaga and Maurice Mignotte (1999),
{\em Embedding the $3x+1$ Conjecture in a $3x+d$ Context, }
Experimental Math. {\bf7}, No. 2 (1999), 145--151. 
(MR 200d:11034). \\
\newline
\hspace*{.25in}
The paper studies iteration on the positive integers 
of the $3x+d$ function
$$
T(x) = \left\{
\begin{array}{cl}
\df{3x+d}{2} & \mbox{if}~~x \equiv 1 ~~ (\bmod ~2)~, \\
~~\\
\df{x}{2} & \mbox{if} ~~x \equiv 0 ~~ (\bmod ~2) ~,
\end{array}
\right.
$$
where $d \geq -1$ and $d \equiv \pm 1~ ( \bmod~6)$.
It proves that there is an absolute constant $c$ such that
there are at most $dk^c$ periodic orbits which contain at
most $k$ odd integers.
Furthermore $c$ is effectively computable.
This follows using a transcendence result of Baker and
W\"ustholz [J. reine Angew. {\bf 442} (1993), 19--62.]

\item
Stefano Beltraminelli, Danilo Merlini and Luca Rusconi (1994),
{\em Orbite inverse nel problema del $3n+1$,}
Note di matematica e fisica,
Edizioni Cerfim Locarno
{\bf 7} (1994), 325--357. \\
\newline
\hspace*{.25in}
This paper discusses the tree of inverse iterates of the Collatz function,
which it terms the "chalice." It states the Collatz conjecture, and notes
that it fails for the map
$$
C_{3,5}(x) := 
\left\{
\begin{array}{cl}
3x+5 & \mbox{if}~ x \equiv 1~~ (\bmod ~2 ) ~, \\
~~~ \\
\df{x}{2} & \mbox{if} ~~x \equiv 0~~ (\bmod~2) ~,
\end{array}
\right.
$$
where there are at least two periodic orbits $\{ 1, 8, 4, 2\}$ and $\{5, 20, 10, 5\}$;
here orbits with initial term $x \equiv 0 ~(\bmod~5)$ retain this property
throughout the iteration. It studies patterns of inverse iterates in the tree
with numbers written in binary.

\item
Lothar Berg and G\"{u}nter Meinardus (1994), 
{\em Functional equations connected with the Collatz problem,}
Results in Math. {\bf 25} (1994), 1--12.
(MR 95d:11025). \\
\newline
\hspace*{.25in}
The $3x+1$ Conjecture is
stated as Conjecture 1. The paper
proves its equivalence 
to each of Conjectures 2 and 3 below, which involve generating functions
encoding iterations of the $3x+1$  function $T(x)$.
For $m,n \ge0$ define $f_m(z) = \sum_{n=0}^\infty  T^{ (m)} (n) z^n$ and 
$g_n(w) = \sum_{n=0}^\infty  T^{ (m)} (n) w^m.$ The paper shows that
each $f_m(z)$ is a rational function of   form
$$
f_m(z) = \frac{p_m(z)}{(1-z^{2^m})^2},
$$
where $p_m(z)$ is a polynomial of degree $2^{m+1}-1$ with integer coefficients. Conjecture 2
asserts that each $g_n(w)$ is a rational function of  the form
$$
g_n(w) = \frac{q_n(w)}{1-w^2},
$$
where $q_n(w)$ is a polynomial with integer
coefficients, with no bound assumed on its degree.
Concerning functional equations,  the authors show first that the $f_m(z)$ satisfy the
recursions
$$
f_{m+1}(z^3)= f_m(z^6)+ \frac{1}{3z} \sum_{j=0}^2 \omega^j f_m(\omega^j z^2),
$$
in which $\omega:= \exp\left( \frac{2 \pi i}{3} \right)$
is a  nontrivial cube root of unity. They also consider the
bivariate generating function $F(z,w) = \sum_{m=0}^\infty
\sum_{n=0}^\infty T^{(m)} (n) z^n w^m$, which converges for $|z| < 1$
and $|y| < \frac{2}{3}$ to an analytic function of two complex variables.
The authors show that it satisfies the functional equation 
$$
F(z^3 ,w) = \frac{z^3}{(1-z^3 )^2} + wF (z^6 , w)
+ \frac{w}{3z} \sum_{j=0}^2 \omega^j F(\omega^j z^2, w)~.
$$
They prove that this functional equation determines $F(z,w)$ uniquely, i.e. there
is only one analytic function of two variables in a neighborhood of $(z,w)=(0,0)$
satisfying it. Next they consider a one-variable  functional
equation obtained from this one
by formally setting $w=1$ (note this falls outside the known region
of analyticity of the function), and  dropping the non-homogenous term. This
functional equation is 
$$
h(z^3) = h(z^6) + \frac{1}{3z} \sum_{j=0}^2 \omega^j h( \omega^j z^2 ).
$$
Conjecture 3 asserts that the only 
solutions $h(z)$ of this functional equation that are
analytic in the unit disk $|z| <    1$ are $h(z)= c_0 + c_1(\frac{ z}{1-z})$ for complex
constants $c_0 , c_1$.

\item
Lothar  Berg and G\"{u}nter  Meinardus (1995), 
{\em The 3n+1 Collatz Problem and Functional Equations, }
Rostock Math. Kolloq. {\bf 48} (1995), 11-18.
(MR 97e:11034). \\
\newline
\hspace*{.25in}
This paper reviews the results of Berg and Meinardus (1994) and adds
some new results. The first new result  considers the functional equation
$$
h(z^3) = h(z^6) + \frac{1}{3z} \sum_{j=0}^2 \omega^{j }h( \omega^{j }z^2 ).
$$
with $\omega:= \exp\left( \frac{2 \pi i}{3} \right)$, and shows 
that the only solutions
$h(z)$ that are entire functions are  constants. 
 The authors next transform 
this functional equation 
%for functions $h(z)$ analytic on the unit disk  
to an equivalent
 system of two functional equations:
\begin{eqnarray*}
h(z)+ h(-z) &= &  2 h(z^2)\\
h(z^3) - h(-z^3) &=& \frac{2}{3z} \sum_{j=0}^2 \omega^{j }h( \omega^{j}z^2).
\end{eqnarray*}
They observe that  analytic solutions  on the open unit disk to
these  two functional equations can be studied separately. 
The second one is the most interesting. 
Set $\Phi(z):= \int_{0}^z h(z) dz$ for $|z|<1.$ Making the change of
variable $z= e^{\frac{2\pi i}{3} \xi}$, the unit disk $|z|<1$ is mapped to the upper
half plane $Im(\xi)>0$. Letting $\phi(\xi):= \Phi ( e^{\frac{2 \pi i}{3} \xi})$ 
the second functional equation above becomes
$$
\phi(3 \xi) + \phi(3\xi+ \frac{3}{2}) = \phi(2 \xi) + \phi(2\xi+1) + \phi(2 \xi +2).
$$
Here we also require $\phi(\xi)= \phi(\xi+3)$. The authors remark that it might also
be interesting to study solutions to this functional equation for $\xi$ on the
real axis. The paper concludes with new formulas for the rational
functions $f_m(z)$ studied in Berg and Meinardus (1994).

\item
Daniel J. Bernstein (1994), 
{\em A Non-Iterative 2-adic Statement of the $3x +  1$ Conjecture,} 
Proc. Amer. Math. Soc., {\bf 121} (1994), 405--408.
(MR 94h:11108). \\
\newline
\hspace*{.25in}
Let $\ZZ_2$ denote the 2-adic integers, 
and for $x  \in  \ZZ_2$ write
$x =  \sum_{i=0}^\infty 2^{d_i}$ with $0 \leq d_0  < d_1 < d_2  <  \cdots$.
Set $ \Phi (x) =  - \sum_{j=0}^\infty  \frac{1}{3^{j+1}}  2^{d_j} $.
The map $ \Phi $ is shown to be a homeomorphism of the 2-adic integers
to itself, which is the inverse of the map $Q_\infty$ defined in 
Lagarias (1985).
The author proves in Theorem 1
a result equivalent to  $ \Phi^{-1} \circ  T \circ  \Phi =  S$,
where $T$ is the $(3x  +  1)$-function on $\ZZ_2$, and $S$ is the
shift map
$$
S(x) =  \left\{
\begin{array}{cl}
\df{x-1}{2} & \mbox{if}~~x \equiv 1~~(\bmod~2 ) ~, \\
~~~ \\
\df{x}{2} & \mbox{if}~~x \equiv 0~~(\bmod~2 ) ~.
\end{array}
\right.
$$
He shows that the $3x  +  1$ Conjecture is equivalent to the conjecture that
$\ZZ^+   \subseteq   \Phi ( \frac{1}{3} \ZZ )$.
He rederives the known results that 
$\QQ  \cap  \ZZ_2  \subseteq  Q_\infty ( \QQ  \cap  \ZZ_2 )$, 
and that $Q_\infty$ is nowhere differentiable, cf. M\"{u}ller (1991).

\item
Daniel J. Bernstein and Jeffrey C. Lagarias (1996), 
{\em The $3x+1$ Conjugacy Map,}
Canadian J. Math., {\bf 48} (1996), 1154-1169. 
(MR 98a:11027). \\
\newline
\hspace*{.25in}
This paper studies the map 
$ \Phi : \ZZ_2  \rightarrow  \ZZ_2$ of Bernstein (1994)
that conjugates the 2-adic shift map to the $3x+1$ function.
This is the inverse of the map $Q_\infty$ in Lagarias (1985); 
see also Akin (2004).
The map $ \bar{\Phi}_n  \equiv   \Phi ~ (\bmod  ~2^n)$ is a permutation of
$\ZZ / 2^n \ZZ$.
This permutation is shown to have order $2^{n-4}$ for $n  \geq  6$.
Let $ \hat{\Phi}_n$ denote the restriction of this permutation to
$( \ZZ / 2^n \ZZ )^\ast =  \{ x: x   \equiv  1~  (\bmod~2) \}$.
The function $ \Phi $ has two odd fixed points $x =  -1$ and $x =  1/3$ 
and the 2-cycle
$\{ 1,  - 1/3 \}$, hence each $ \hat{\Phi}_n$ inherits 
two 1-cycles and a 2-cycle
coming from these points.
Empirical evidence indicates that $ \hat{\Phi}_n$ has about $2n$
fixed points for $n \leq 1000$.
A heuristic argument based on this data suggests that 
$-1$ and $1/3$ are the only
odd fixed points of $ \Phi $.
The analogous conjugacy map $ \Phi_{25, -3}$ for the `$25x-3$' problem
is shown to have no nonzero fixed points.

\item
Jacek B\l a\.{z}ewitz and Alberto Pettorossi (1983),
{\em Some properties of binary sequences useful for proving Collatz's
conjecture,}
J. Found. Control Engr. {\bf 8} (1983), 53--63.
(MR 85e:11010, Zbl. 547.10000). \\
\newline
\hspace*{.25in}
This paper studies
 the $3x+1$ Problem interpreted
as a strong termination property of a term
rewriting system. They view the problem as transforming binary strings
into new binary strings and look in partcular at its action on the patterns
$1^n$, $0^n$ and $(10)^n$ occurring inside strings. The $3x+1$ map
exhibits regular behavior relating these patterns.

\item
Richard Blecksmith, Michael McCallum, and John L. Selfridge (1998),
{\em 3-Smooth Representations of Integers,}
American Math. Monthly {\bf 105} (1998), 529--543.
(MR 2000a:11019). \\
\newline
\hspace*{.25in} 
A {\em 3-smooth representation} of an integer $n$ is
a representation as a sum of distinct positive integers each of which has the
form $2^a 3^b$, and no term divides any other term.  This paper
proves a conjecture of Erdos and Lewin that for each
integer $t$ all sufficiently large integers have a  3-smooth representation
with all individual terms larger than $t$. They note a
connection of 3-smooth representations to the $3x + 1$-problem,
 which is that a number $m$ iterates
to $1$ under the $3x+1$ function if and only if there are positive
integers $e$ and $f$
such that $n = 2^e - 3^f m$ is a positive integer that 
 has a 3-smooth
representation with $f$ terms in which there is one term exactly 
divisible by each 
power of three from $0$ to $f - 1$. The choice of $e$ and $f$ is not
unique, if it exists.

\item
Corrado B\"{o}hm and Giovanna Sontacchi (1978), 
{\em On the existence of cycles of given length in integer sequences
like $x_{n+1}= x_n/2$ if $x_n$ even, and $x_{n+1}= 3x_n + 1$ otherwise,}
Atti Accad. Naz. Lincei Rend. Cl. Sci. Fis. Mat. Natur. {\bf 64}
(1978), 260--264.
(MR 83h:10030) \\
\newline
\hspace*{.25in}
The authors are primarily concerned
with  cycles of a generalization of the $3x+1$
function. They consider the recursion in which $x_{n+1} = a x_n +b$,
if a given  recursive predicate
$P(x_n)$ is true, and $x_{n+1} = c x_n + d$ if the predicate $P(x_n)$
false, where $a, b, c, d$ and $x_n$
are rational numbers.  They observe that as a consequence of linearity alone
there are at most $2^k$ possible cycles of period $k$, corresponding
to all possible sequences of ``true'' and ``false'' of length $n$.
Furthermore one can effectively determine the $2^n$ rationals that are the
solutions to each of these equations and check if they give cycles.
Thus in principle one can determine all cycles below any given 
finite bound. They observe that a rational number $x$ in
a cycle of the $3x+1$-function $T(\cdot)$ 
of period $n$ necessarily has the form
$$
x = \frac{\sum_{k=0}^{n-1} 3^{m-k-1} 2^{v_k}}{2^n - 3^m}
$$
with $0 \le v_0 < v_1 < ... < v_m =n$. They deduce  that
every integer $x$ in a cycle of length $n$ necessarily has
$|x| < 3^n.$  

Further study of rational cycles of the $3x+1$ function
appears in Lagarias (1990).

\item
David W. Boyd (1985), 
{\em Which rationals are ratios of Pisot sequences?,}
Canad. Math. Bull. {\bf 28 } (1985), 343--349.
(MR 86j:11078). \\
\newline
\hspace*{.25in}The Pisot sequence $E(a_0 ,  a_1 )$ is defined by
$a _{n+2} = \left[ \frac{a_{n+1}^{2}}{a_n }  +  \frac{1}{2} \right]$,
where $a_0 ,  a_1$ are integer starting values.
If $0 <  a_0  <  a_1$ then $\frac{a_n}{a_n+1}$ converges to
a limit $\theta$
as $n \rightarrow  \infty$.
The paper asks: which rationals $\frac{p}{q}$ can occur as a limit?
If $\frac{p}{q}  >  \frac{q}{2}$ then $\frac{p}{q}$ must be an integer.
If $\frac{p}{q}  <  \frac{q}{2}$ the question is related to
a stopping time problem resembling the
$3x  +  1$ problem.

\item
Monique Bringer (1969), 
{\em Sur un probl\`{e}me de R. Queneau},
Math\'{e}matiques et Sciences Humaines
[Mathematics and Social Science], {\bf 27} Autumn 1969, 13--20.\\
\newline
\hspace*{.25in}
This paper considers a  problem proposed by Queneau (1963)
in connection with rhyming patterns in poetry. It concerns,
 for a fixed $n \ge 2$,
iteration of the map on the integers
$$
\delta_n(x) := \left\{ 
 \begin{array}{cl}
\df{x}{2} & ~~\mbox{if} ~~x  ~~\mbox{is ~even}\\
~&~\\
\df{2n+ 1-x}{2} & ~~\mbox{if} ~~x  ~~\mbox{is ~odd}\\
\end{array}
%\.
\right.
$$
This map acts as a permutation on the integers $\{1, 2, ..., n\}$ and
it also has the fixed point $\delta_n(0)=0$. It is called by the author
a  spiral permutation of $\{1, 2, ..., n\}$.  
The paper studies for which $n$ this spiral  is a cyclic permutation, and calls such
numbers {\em admissible}. 

The motivation for this problem was that this permutation for $n=6$ represents
a poetic stanza pattern, the sestina,  used in poems by an 11th century Troubadour, Arnaut Daniel. 
This pattern for general $n$ was studied by Raymond Queneau (1963),
% a member of the French mathematical-literary group Oulipo. 
who determined  small values
of $n$ giving  a cyclic permutation. His colleague Jacques Roubaud (1969) termed
these rhyme schemes $n-$ines or {\em quenines}.
Later he called these numbers  {\em Queneau numbers}, cf. Robaud (1993). 

In this paper the author, a student of Roubaud, shows that a necessary condition for a number $n$ to
be admissible  is that $p=2n+1$ be prime. She shows that
a sufficient condition to be admissible is that
$2$ be a primitive root $(\bmod ~p)$. 
She deduces that if $n$ and $2n+1$ are  both primes then $n$ is admissible, and 
that the numbers $n=2^k$ and $n=2^k -1$ are never admissible.
Finally she shows that all $p \equiv 1 ~(\bmod ~8)$
are not admissible.

{\em Note.} The function $\delta_n(x)$ is defined on the integers, and  is  of
$3x+1$ type (i.e. it is a periodically linear function).  Its long term behavior under iteration is 
 analyzable because the function $\delta_n(x)$ 
 decreases absolute value on each iteration for any $x$ with $|x|> 2n+1.$  Thus all
 orbits eventually enter $-2n-1\le x\le 2n+1$ and become eventually periodic.
 One can further show that all orbits eventually enter the region $\{0, 1,..., n\}$ on which 
 $\delta_n(x)$ is a permutation.

\item
Stefano Brocco (1995),
{\em A Note on Mignosi's Generalization of the $3x+1$ Problem, }
J. Number Theory, {\bf 52} (1995), 173--178. (MR 96d:11025). \\
\newline
\hspace*{.25in}
F. Mignosi (1995) studied the function 
$T_\beta : \NN  \rightarrow  \NN$ defined by
$$
T_\beta (n) =  \left\{
\begin{array}{cl}
\lceil \beta n \rceil & \mbox{if}~~n \equiv 1~~(\bmod~2 ) \\
~~~ \\
\df{n}{2} & \mbox{if}~n \equiv 0~~(\bmod~2 ) ~,
\end{array}
\right.
$$
where $\lceil x \rceil$ denotes the smallest integer $n  \geq  x$.
He also formulated Conjecture $C_\beta$ asserting that $T_\beta$
has finitely many cycles and that every
$n  \in  \NN$ eventually enters a cycle under $T_\beta$.
This paper shows that Conjecture $C_\beta$ is false whenever $\beta$
is a Pisot number or a Salem number.
The result applies further to functions
$$
T_{\beta , \alpha} (n) = \left\{
\begin{array}{cl}
\lceil \beta n + \alpha \rceil & \mbox{if}~~ n \equiv 1 ~( \bmod ~2)~, \\
~~~ \\
\df{n}{2} & \mbox{if}~~ n \equiv 0 ~ (\bmod ~2)~,
\end{array}
\right.
$$
for certain ranges of values of $\beta$ and $\alpha$.

\item
Serge Burckel (1994),
{\em Functional equations associated with congruential functions,}
Theoretical Computer Science {\bf 123} (1994), 397--406.
(MR 94m.11147). \\
\newline
\hspace*{.25in}
The author proves undecidability results for periodically 
linear functions generalizing those of
Conway (1972).
% [Proc. 1972 Number Theory Conference, Boulder, Colorado, 49--52].
A periodically linear function $f : \ZZ  \rightarrow  \ZZ$ is one which
is a linear function on each congruence class $(\bmod L$) for 
some finite $L$.
The author shows it is undecidable whether a given function has 
$f^{(k)} (1) =  0$
for some $k  \geq  1$, and also whether a given function has the property:
for each $n \geq  1$, some $f^{(k)} (n) =  0$.
He also shows
that the $3x+1$ conjecture is equivalent to a
certain functional equation having only the trivial solution
over the set of all
power-series $R(z) =  \sum_{n=}^\infty  a_n z^n$
with all $a_i =  0$ or 1.
The functional equation is
$$
3z^3 R( z^3 ) -3z^9  R(z^6) - R(z^2) - 
R ( \omega z^2 ) -R ( \omega^2  z^2 ) = 0
$$
where $\omega =  \exp ( \frac{2 \pi i}{3} )$.

\item
Robert N. Buttsworth and Keith R. Matthews (1990), 
{\em On some Markov matrices arising from the generalized Collatz mapping,}
Acta Arithmetica {\bf 55} (1990), 43--57.
(MR 92a:11016). \\
\newline
\hspace*{.25in}This paper studies maps 
$T (x) = \frac {m_i  x  -  r_i}{d}$ for
$x  \equiv  i (\bmod ~d)$, where
$r_i  \equiv  im_i (\bmod ~d)$.
In the case where g.c.d. $(m_0 , \ldots , m_{d-1} ,  d) =  1$ 
it gives information about the structure
of $T$-ergodic sets $(\bmod~m)$ as $m$ varies.
A set $S \subseteq  \ZZ$ is
{\em T-ergodic} $(\bmod~m)$
if it is a union of $k$ congruence classes $(\bmod~m)$,
$S =  C_1 \cup \ldots  \cup C_k$,
such that $T(S)  \subseteq S$ and
there is an $n$ such that $C_j  \cap  T^{(n)} (C_i ) \neq \phi$
holds for all $i$ and $j$.
It characterizes them in many cases.
As an example, for
$$
T(x) =  \left\{
\begin{array}{cl}
\df{x}{2} & \mbox{if}~~x \equiv 0~~(\bmod~2 )~, \\
~~~ \\
\df{5x-3}{2} & \mbox{if}~~x \equiv 1~~(\bmod~3) ~,
\end{array}
\right.
$$
the ergodic set $(\bmod ~m)$ is
unique and is
$\{   n : n   \in \ZZ$ and
$(n, m, 15) =  1 \}$, i.e. it is one of
$\ZZ ,  \ZZ -3,\ZZ  -  5 \ZZ$ or $\ZZ  -  3 \ZZ  -  5 \ZZ$
as $m$ varies.
An example is given having infinitely many different ergodic sets
$(\bmod ~m)$ as $m$ varies.

\item
Charles C. Cadogan (1984), 
{\em A note on the $3x +  1$ problem,}
Caribbean J. Math. {\bf 3} No. 2  (1984), 69--72.
(MR 87a:11013). \\
\newline
\hspace*{.25in}
Using an observation of S. Znam, this paper 
shows that to prove the $3x +  1$
Conjecture it suffices to check it for all $n  \equiv  1~~(\bmod ~4)$.
This result complements   the obvious fact that to prove the $3x  +  1$ Conjecture
it suffices to check it for all $n \equiv  3~ (\bmod ~4)$.
Korec and Znam (1987) obtained  other results in this spirit, for odd
moduli.

\item
Charles C. Cadogan (1991),
{\em Some observations on the $3x+1$ problem},
Proc. Sixth Caribbean Conference on Combinatorics \& Computing,
University of the West Indies: 
St. Augustine Trinidad ( C. C. Cadogan, Ed.)
Jan. 1991, 84--91.\\
\newline
\hspace*{.25in}
Cadogan (1984) reduced the $3x+1$ problem to the study of its
iterations on numbers $A_1=\{n: n \equiv 1~~(\bmod~4)\}$. Here the author
notes in particular  the subclass $A=\{ 1, 5, 21, 85,...\}$ where $x_{i+1} = 1+ 4x_i$ of $A_1$.
He considers the successive odd numbers occurring in the iteration. 
He forms a two-dimensional table partitioning all odd integers in which $A_1$ is the first row
and the first column of the table is called the anchor set (see Cadogan (1996) for more details) .
He observes that  Cadogan (1984) showed the iteration on higher rows
successively moves down rows to the first row, but from  row $A_1$ is flung back to 
higher rows, except for the subclass $A$, which remains on  the first row.
He comments that after
each revisit to row $A_1$ "the path may become increasingly
unpredictable." He concludes, concerning further work: "The target set $A_1$ is being vigorously
investigated."

\item
Charles  C. Cadogan (1996),
{\em Exploring the $3x+1$ problem I.},
Caribbean J. Math. Comput. Sci. {\bf 6} (1996), 10--18.
(MR 2001k:11032) \\
\newline
\hspace*{.25in}
This paper studies itertion of the Collatz function $C(x)$, which is here denoted $f(x)$,
and it includes most of the results of Cadogan (1991). 
The author  gives various criteria under  which trajectories will coalesce.
He partitions the odd integers $2\NN+1= \cup_{k=1}^{\infty} R_k$, in which
 $R_m= \{ n\equiv 2^m -1 ~(\bmod~2^{m+1})\}$.
It enumerates their elements $R_{m, j} = (2^m -1)+ (j-1) 2^{m+1}.$
and views these in an infinite two-dimensional array in which the $j$-th column $C_j$ consists
of the numbers $\{ R_{m, j}:  m \ge 1\}$.  Cadogan (1984)  showed that if
$x \in R_m$ for some $m \ge 2$, then $f^{(2)}(x) \in R_{m-1}$, thus after $2m$ iterations
one reaches an element of $R_1$.  Lemma 3.3 here observes that consecutive 
elements in columns are related by
by $R_{m+1, j} = 1+ 2 R_{m, j}.$  
For the
first set $R_1= \{ n \equiv 1~(\bmod~4)\}$, Theorem 4.1 observes that if $y=4x+1$ then $f^{(3)}(y)= f(x)$. 
The author creates chains $\{x_n: n\ge 1\}$ related by $x_{n+1}= 4 x_n +1$
and calls these $S$-related elements. Theorem 4.2 then observes that the trajectories
of $S$-related elements coalesce. 
%The paper concludes  with a recipe which 
 %given an odd integer $n$ will locate the unique pair $(m,j)$ for which $R_{m,j}=n$. 

\hspace*{.25in}
\item
Marc Chamberland (1996),
{\em A Continuous Extension of the $3x+1$ Problem to the Real Line,}
Dynamics of Continuous, Discrete and Impulsive Dynamical Systems,
{\bf 2} (1996), 495--509. (MR 97f:39031). \\
\newline
\hspace*{.25in}
This paper studies the iterates on $\RR$ of the function
\begin{eqnarray*}
f(x) & = & \df{\pi x}{2} \left( \cos \frac{\pi x}{2}\right)^2 + 
\df{3x+1}{2} \left( \sin \df{\pi x}{2}\right)^2 \\
  & = & x + \df{1}{4} - \left( \df{x}{2} + \df{1}{4} \right) \cos \pi x ~.
\end{eqnarray*}
which interpolates the $3x+1$ function $T(\cdot)$.
A fact crucial to the analysis is that $f$ has negative Schwartzian derivative
$Sf = \frac{f'''}{f'} - \frac{3}{2} \left( \frac{f''}{f'} \right)^2$
on $\RR^+$.
On the interval $[0, \mu_1 )$, where $\mu_1 = 0.27773 \ldots~$ 
all iterates of $f$ contract to a fixed point 0. 
Here $\mu_n$ denotes the $n-th$ positive fixed point of $f$.
The interval $[\mu_1 , \mu_3 ]$ is invariant under $f$,
where $\mu_3 = 2.44570 \ldots$ and this interval includes the trivial
cycle $A_1 = \{ 1,2 \}$.
On this interval almost every point is attracted to one of two
attracting cycles, which are $A_1$ and 
$A_2 = \{ 1.19253 \ldots ,~ 2.13865 \ldots \}$.
There is also an uncountable set of measure 0 on which the dynamics is
``chaotic.''
On the interval $[\mu_3 , \infty )$
the set of $x$ that do not eventually iterate to a point in 
$[ \mu_1 , \mu_3 ]$
is conjectured to be of measure zero.
The point $\mu_3$ is proved to be a ``homoclinic point,''
in the sense that for any $\epsilon > 0$ the iterates of
$[ \mu_3 , \mu_3 + \epsilon )$
cover the whole interval
$( \mu_1 , \infty )$.
It is shown that any nontrivial cycle of the $3x+1$ function on the 
positive integers would be an attracting
periodic orbit of $f$.

\item
Busiso P. Chisala (1994), 
{\em Cycles in Collatz Sequences, }
Publ. Math. Debrecen {\bf 45} (1994), 35--39.
(MR 95h:11019). \\
\newline
\hspace*{.25in}
The author
shows that for any $m$-cycle of the Collatz map on positive
{\em rationals}, the least element is at least as large as
$(2^{ [ m \theta ] /m} - 3 )^{-1}$, where
$\theta = \log_2 3$.
Using this result, he derives a lower bound for $3x+1$ cycle lengths 
based on the continued
fraction of $\theta = [ 1, a_1 , a_2 , a_3 , \ldots ]$, in which the $n$-th
convergent is $\frac{p_n}{q_n}$, and the intermediate convergent
denominator $q_n^i $ is $ iq_{n+1} + q_n$ for $0 \leq i < a_{n+1}$.
If the $3x+1$ conjecture is true for $1 \leq n \leq N$, and
$N \geq ( 2^{C(i,k)} -3)^{-1}$, where 
$C(i,k) = \frac{\lceil q_k^i \theta \rceil}{q_k^i}$, then there 
are no nontrivial cycles of the $3x+1$ function on
$\ZZ^+$ containing less than $q_k^{i+1}$ {\em odd} terms.
Using the known bound $N = 2^{40} \doteq 1.2 \times 10^{12}$,
the author shows that there are at least $q_{15} = 10~787~915$ 
odd terms in any cycle
of the $3x+1$ function on $\ZZ^+$.
\newline
[Compare these results with those of Eliahou (1993).]

\item
Vasik Chvatal, David A. Klarner and Donald E. Knuth (1972)
{\em Selected combinatorial research problems,}
Stanford Computer Science Dept. Technical Report STAN-CS-72-292,
June 1972, 31pages. \\
\newline
\hspace*{.25in}
This report contains thirty-seven research problems, the first 16 of which
are due to Klarner, the next 9 to Chvatal, and the remaining 11 to Knuth.
This list contains two problems about iterating affine maps.
Problem 1 asks whether the set of all positive integers reachable from $1$
using the maps $x \mapsto 2x+1$ and $x \mapsto 3x+1$ can be 
partitioned into a disjoint union  of  infinite arithmetic progressions. 
Problem 14 considers for nonnegative integers $(m_1, ..., m_r)$ the
set $S = \langle m_1x_1+ \cdots + m_r x_r: 1\rangle $ which is the smallest set
of natural numbers containing $1$ and which is closed under the operation
of adjoining $m_1x_1 +\cdots + m_r x_r$ whenever $x_i$ are in the set.
It states that Klarner has shown that $S$ is a finite union of arithmetic progressions
provided that $(i)$ $r \ge 2$, $(ii)$ the greatest common divisor $(m_1, ..., m)r)=1$,
and, $(iii)$ the greatest common divisor $(m_1+m_2 + \cdots +m_r, \prod_{i} m_i) = 1.$
It asks if the same conclusion holds if hypothesis (iii) is dropped.

{\em Note.}  Problem 1 was solved in the affirmative in Coppersmith (1975).
Problem 14 relates to the theory developed in Klarner and Rado (1974),
and was solved affirmatively in Hoffman and Klarner (1978), (1979). 

\item
Dean Clark (1995), 
{\em Second-Order Difference Equations Related to the 
Collatz $3n+1$ Conjecture, }
J. Difference Equations \& Appl., {\bf 1} (1995), 73--85.
(MR 96e:11031). \\
\newline
\hspace*{.25in}
The paper studies the integer-valued recurrence
$\frac{x_{n+1} + x_n}{2}$ if
$x_{n+1} + x_n$ is even, and
$x_n = \frac{b|x_{n+1} -x_n | +1}{2}$ if
$x_{n+1} + x_n$ is odd, for $b \geq 1$ an odd integer.
For $b =  1,3,5$ all recurrence sequences stabilize at some fixed point 
depending on $x_1 $
and $x_2$,
provided that $x_1 = x_2 \equiv \frac{b+1}{2}~ (\bmod ~b)$.
For $b  \geq  7$ there exist unbounded trajectories, and
periodic trajectories of period $\geq 2$.
In the ``convergent'' cases $b =  3$ or 5 the iterates exhibit an
interesting phenomenon, which the author calls
 {\em digital convergence,}
where the low order digits in base $b$ of $x_n$ successively
stabilize before the high order bits stabilize.

\item
Dean Clark and James T. Lewis (1995),
{\em A Collatz-Type Difference Equation, }
Proc. Twenty-sixth Internationsal Conference on Combinatorics,
Graph Theory and Computing (Boca Raton 1995),
Congr. Numer. {\bf 111} (1995), 129-135. (MR 98b:11008). \\
\newline
\hspace*{.25in}
This paper studies the difference equation
$$
x_n = \left\{
\begin{array}{cl}
\df{x_{n-1} + x_{n-2}}{2} & \mbox{if}~~x_{n-1} + x_{n-2}~\mbox{ is even}~, \\
~~~ \\
x_{n-1} - x_{n-2} & \mbox{if}~~x_{n-1} + x_{n-2}~ \mbox{is~odd} ~.
\end{array}
\right.
$$
with integer initial conditions $(x_0, x_1).$
It suffices to treat the case that $gcd(x_0 , x_1 )  = 1 $.
For such initial conditions the recurrence is shown to always converge to 
one of the 1-cycles 1 or $-1$
or to the 6-cycle $\{ 3,2,-1, -3, -2, -1 \}$.
%PAGEBREAK
%\newpage
\item
Dean Clark and James T. Lewis (1998),
{\em Symmetric solutions to a Collatz-like
system of  Difference Equations, }
Proc. Twenty-ninth Internationsal Conference on Combinatorics,
Graph Theory and Computing (Baton Rouge 1998),
Congr. Numer. {\bf 131} (1998), 101-114. (see MR 99i:00021). \\
\newline
\hspace*{.25in}
This paper studies the first order system of nonlinear difference equations
\begin{eqnarray*}
x_{n+1} & = & \lfloor \frac{x_n+y_n}{2} \rfloor\\
y_{n+1} &=& y_n - x_n,
\end{eqnarray*}
where  $\lfloor . \rfloor$ is the floor function (greatest integer function). Let 
$T(x, y) = ( \lfloor \frac{x+y}{2}\rfloor, y-x)$ be a map of the plane, noting that 
$T(x_n, y_n) = (x_{n+1}, y_{n+1})$. The function
$T$ is an invertible map of the plane, with inverse $S= T^{-1}$ given
by  $S(x,y)=( \lceil x - \frac{y}{2} \rceil, \lceil x+ \frac{y}{2} \rceil)$, using
the ceiling function $\lceil .\rceil$. 
One obtains an associated linear map of the plane,  by not imposing the floor function
above., i.e. $\tilde{T}(x, y) = ( \frac{x+y}{2}, -x+y)$.  The map $\tilde{T}$  is invertible, 
and  for arbitrary real initial conditions $(x_0, y_0)$ the full orbit 
$\{ \tilde{T}^{(k)}(x_0, y_0): -\infty < k < \infty\}$
 is bounded, with all points on it being confined to an invariant ellipse.
 The effect of the floor function is to  perturb this linear
dynamics. The authors focus on the question of whether all orbits having
integer initial conditions $(x_0, y_0)$ remain bounded; however they don't
resolve this question. 
Note that integer initial conditions imply the full orbit is integral;
then  invertibility implies
that bounded orbits of this type must be periodic.

The difference equation given above for integer initial conditions $(x_0, y_0)$
can  be transformed to a second order nonlinear recurrence by 
eliminating the variable $x_n$, obtaining. 
$$
y_{n+1} :=\left\{
\begin{array}{cl}
\frac{3y_{n}+1}{2} - y_{n-1}  & \mbox{if} ~~y_n  \equiv 1~~ (\bmod ~ 2) \\
~~~ \\
\frac{3y_{n}}{2} - y_{n-1}  & \mbox{if}~~ y_n \equiv 0~~ (\bmod ~2 )~.
\end{array}
\right.
$$
with integer initial conditions $(y_0, y_1)$.  They note a resemblance of
this recurrence in form
to the $3x+1$ problem, and 
view  boundedness of orbits as a (vague) analogue of the $3x+1$ Conjecture.
Experimentally they observe that all integer orbits appear
to be periodic, but the period of such orbits varies erratically with the
initial conditions. For example the starting
condition $(64,0)$ for $T$ has period 87, but that of $(65,0)$ has period 930.
They give a criterion (Theorem 1) for an integer orbit to be unbounded, but conjecture it is never
satisfied.

The paper also studies properties of periodic orbits imposed by some symmetry operators. 
They introduce the  operator  $Q(x,y):= ( \lfloor -x + \frac{y}{2})\rfloor, y)$,
observe it is an involution $Q^2=I$ satisfying $(TQ)^2=I$ and $S= QTQ^{-1}$.
They also introduce a second
symmetry operator $U(x,y): = (-x, 1-y)$ which is 
is an involution $U^2=I$ that commutes with $T$. These operators are used
to imply some symmetry
properties of periodic orbits, with respect to the line $x=y$. They also derive the
result (Theorem 6): the sum of the terms $y_k$ that are even integers in a complete period
of a periodic orbit is divisible by $4$; the sum of all the $y_k$ over a cycle is
strictly positive and equals the number of odd $y_k$ that appear in the cycle.
 
 The authors do not explore what happens to orbits of
 $T$ for general real  initial conditions; should all orbits remain bounded in
 this more general situation?

\item
Thomas Cloney, Eric C. Goles and G\'{e}rard Y. Vichniac (1987), 
{\em The $3x  +  1$ Problem: a Quasi-Cellular Automaton,}
Complex Systems {\bf 1}(1987), 349--360.
(MR 88d:68080). \\
\newline
\hspace*{.25in}
The paper presents computer graphics pictures of 
binary expansions of $ \{  T^{(i)}(m): i= 1, 2, \ldots  \}$ 
for ``random'' large $m$, using black and white pixels to represent 1 resp.~0
$(\bmod~2)$.  It discusses patterns seen in
these pictures.
There are no theorems.

\item
Lothar  Collatz (1986),
{\em On the Motivation and Origin of the $(3n +  1)$-Problem, }
J. of Qufu Normal University, Natural Science Edition
[Qufu shi fan da xue xue bao. Zi ran ke xue ban]
{\bf 12} (1986) No. 3, 9--11 (Chinese, transcribed by Zhi-Ping Ren).\\
%[German translation: \"Uber den Ursprung des $(3n +  1)$-Problems,
%by Zhang-Zheng Yu 1991, U.\ Hamburg] \\
\newline
\hspace*{.25in}
Lothar Collatz describes his interest since 1928 
in iteration problems
represented using associated graphs and hypergraphs.
He describes the structure of such graphs for several different problems.
He states that he invented the $3x  +  1$ problem and publicized it
in many talks.
He says:
``Because I couldn't solve it I never published anything.
In 1952 when I came to Hamburg I told it to my colleague 
Prof.\ Dr.\ Helmut Hasse.
He was very interested in it.
He circulated the problem in seminars and in other countries.''

{\em Note.} Lothar Collatz was part of the DMV delegation to the
1950 International Congress of Mathematicians in Cambridge, Massachusetts. 
There Kakutani and Ulam were invited speakers. Collatz reportedly
described the problem at this time to Kakutani and others in private conversations,
cf. Trigg et al (1976). 

\item
John H. Conway (1972),
{\em Unpredicatable Iterations},
In: Proc. 1972 Number Theory Conference, University of Colorado,
Boulder, CO. 1972, pp. 49--52. (MR 52 \#13717). \\
\newline
\hspace*{.25in}
This paper states the $3x+1$ problem, and shows
that a more general function iteration problem similar in form to
the $3x+1$ problem  is computationally undecidable.

The paper considers functions $f: \ZZ \to \ZZ$ for which there
exists a finite modulus $N$ and rational numbers
$\{a_j : 0 \le j \le N-1 \}$ such that 
$$
g(n) = a_j n ~~~\mbox{if}~~ n \equiv j~(\bmod~N).
$$
In order that the map take integers to integers it is
necessary that the denominator of $a_j$ divide $gcd(j, N)$.
The computationally undecidable question becomes: Given an
$f(\cdot)$ in this class and an input value $n = 2^k$ decide
whether or not some iterate of $n$ is a power of $2$.
More precisely, he shows that for any recursive 
function $f(n)$ there exists a choice of $g(\cdot)$ such
that for each $n$ there holds
$2^{f(n)} = g^{k}(2^n)$ for some $k \ge 1$ and
this is the smallest value of $k$ for which the
iterate is a power of $2$.
It follows that there is no decision procedure to 
recognize if the iteration of such a function, starting
from input $n=2^j$, will ever encounter another power of
$2$, particularly whether it will encounter the value $1$.
%For if it could, then it could recognize all recursive
%functions.

The proof uses an encoding of computations in the exponents
$e_p$ of a multiplicative factorization
$n = 2^{e_2} \cdot 3^{e_3} \cdot ..$,  in which only a fixed
finite number of exponents $(e_2, e_3, ..., e_{p_r})$
control the computation, corresponding to the primes dividing
the numerators and denominators of all $a_j$.
The computation is based on  a machine model with
a finite number registers storing integers of arbitrary size,
(the exponents $(e_2, e_3, ..., e_{p_r})$)
and there is a finite state controller.
These are 
called {\em Minsky machines} in the literature,
and are described in  Chapter 11 of 
M. Minsky, {\em Computation: Finite and Infinite Machines},
Prentice-Hall: Englewood Cliffs, NJ 1967 (especially Sec. 11.1).

Conway (1987) later formalized  this computational model 
as FRACTRAN,
and also constructed a 
universal function $f(\cdot)$. See Burckel (1994) for other 
undecidability results.

\item
John H. Conway (1987),
{\em FRACTRAN- A Simple Universal Computing Language for Arithmetic,}
in: {\em Open Problems in Communication and Computation}
(T. M. Cover and B. Gopinath, Eds.), Springer-Verlag: New York 1987,
pp. 3--27. (MR 89c:94003). \\
\newline
\hspace*{.25in}
FRACTRAN is a method of universal computation based on 
Conway's (1972)  earlier analysis in ``Unpredictable Iterations.''
Successive computations are done by multiplying the current
value of the computation, a positive  integer, by one of a finite list
of fractions, according to a definite rule which guarantees
that the resulting value is still an integer. A  FRACTRAN program
iterates a function $g(.)$ of the form 
$$
g(m) := \frac{p_r}{q_r}m   ~~~~\mbox{if}   ~~m \equiv r ~ (\bmod~N),  ~~0 \le r \le N-1, 
$$
and each fraction $\frac{p_r}{q_r}$  is positive with 
 denominator $q_r$ dividing $gcd(N, r)$,
so that the function takes positive integers as positive integers. 
If the input integer is of the form
 $m=2^n$ then the FRACTRAN program is said to {\em halt} at the first value
encountered which is again a power of $2$, call it  $2^{f(n)}$. The output $f(n)= \ast$ is
undefined if the program never halts. A FRACTRAN
program is regarded as computing the
partial recursive function $\{ f(n) ~:~ n \in \ZZ_{>0}\}$.
FRACTRAN programs can compute any partial recursive function.
The paper gives a number of  examples of FRACTRAN programs, e.g. 
for computing the decimal digits of $\pi$, and for computing the
successive primes. The prime producing algorithm was  described
earlier  in Guy (1983b).

Section 11 of the paper  discusses generalizations of the 3x+1 problem 
encoded as FRACTRAN programs, including a fixed such function for
which the halting problem is undecidable. 

\item
Don Coppersmith (1975),
{\em The complement of certain recursively defined sets},
J. Combinatorial Theory, Series A {\bf 18} (1975), No. 3, 243--251. 
(MR 51 $\#$5477) \\
\newline
\hspace*{.25in}
This paper studies sets  of nonnegative integers 
$S =\langle a_1x+ b_1, ..., a_k x+ b_k: c_1,..., c_m \rangle$
generated using a family of affine functions.
$x \mapsto a_i x + b_i$,   starting from the set of seeds
$A := \{ c_1, ..., c_m\}$. Here one proceeds  by producing new integers by 
the maps $x \mapsto a_i x + b_i$, in which all $a_i \ge 2, b_i \ge 0$ and
$c_i \ge 1$. He calls such sets {\em $RD$-sets}. Such a set is called {\em good}
if its complement can be given as a disjoint union of arithmetic progressions. 
He reduces the problem to the case of a single seed $m=1$, and shows 
that a set is good if and only if its complement can be covered with 
infinite arithmetic progressions. A (possibly negative) integer is a  {\em feedback
element}  if it is
a fixed point of some sequence of iterates of the maps above; such elements
need not be part of the $RD$-set. Theorem 1 says that for a fixed
set of maps,  if either there are no feedback elements, or  if there 
are  but no image of any of them under repeated iteration of the maps
becomes positive, then for any seed $c$ the associated
$RD$-set is good. Theorem 2 then gives a sufficient condition on a set of maps
for there to exist at least one seed giving a bad $RD$-set. In particular, Corollary 2b
says that  if  $a \ge 2, m\ge 0$ and $b \ge1$ 
then the set $S(c) := < ax + m(a-1), ax+ m(a-1) + b: c>$
is bad for some $c \ge 1$. Theorem 3 gives a stronger sufficient condition, very
complicated to state,  for
 a set of operators to have some seed giving a bad $RD$-set.  The author
 asserts that this sufficient condition is ``almost necessary."

 %\medskip
 {\em Note.} Klarner (1972)  asked if it was true that 
 $S := \langle 2x+1, 3x+1: 1 \rangle$ is a good set.
 That it is a good set follows using Theorem 1,  since it is easy to show by induction 
 that all  operators obtained by composition necessarily have the form $ax+ b$
 with $a >b \ge 1$, so such operators have a fixed point  $x$ satisfying $-1 < x < 0$. Thus
 this set of operators has no feedback elements. Theorem 1 imples that 
 $S(c):= \langle  2x+1, 3x+1: c \rangle$ is good for all $c \ge 1$.

\item
H. S. M. Coxeter (1971),
{\em Cyclic Sequences and Frieze Patterns,
(The Fourth Felix Behrend Memorial Lecture)}, 
Vinculum {\bf 8} (1971),  4--7.\\
\newline
\hspace*{.25in}
This lecture was given at the University of Melbourne in
1970. In this written version of  the lecture,  Coxeter discusses various integer
sequences, including the Lyness iteration $u_{n+1} = \frac{1+u_n}{u_{n-1}}$,
which has  the orbit  $(1,1, 2,3,2)$ as one solution. He observes that
a general solution to the Lyness iteration  can be produced by a frieze pattern
with some indeterminates.
He then introduces  the $3x+1$ iteration as "a more recent piece of
mathematical gossip." He states the $3x+1$ conjecture and then says:
"I am tempted to follow the example of Paul Erd\H{o}s, who offers prizes
for the solutions of certain problems. It the 
above conjecture is true, I will gladly offer a prize of fifty dollars to the 
first person who send me a proof that I can understand. If it is false, I offer
a prize of a hundred dollars to the first one who can establish a counterexample.
I must warn you not to try this in your heads or on the back of an old envelope, 
because the result has been tested with an electronic computer
for all $x_1 \le 500, 000.$"

{\em Note.}
Based on this talk, Coxeter  was credited in Ogilvy (1972) with proposing the $3x+1$ problem. 
For more on Frieze patterns, see 
H. S. M. Coxeter , {\em Frieze patterns}, Acta Arithmetica {\bf 18} (1971), 297--310,
and 
J. H. Conway and H. S. M. Coxeter {\em Triangulated polygons and Frieze patterns,}
Math. Gazette {\bf  57} (1973) no. 400, 87-94; no 401, 175--183. 
Vinculum is the journal of the Mathematical Association of Victoria (Melbourne,
Australia).

\item
Richard E. Crandall (1978),
{\em On the ``$3x + 1$'' problem},
Math. Comp. {\bf 32} (1978), 1281--1292.
(MR 58 \#494). \\
\newline
\hspace*{.25in}
This paper studies iteration of the  ``3x+1'' map
and more generally the ``$qx +r$'' map
$$
T_{q,r}(x) = \left\{
\begin{array}{cl}
\frac{qx+r}{2} & \mbox{if} ~~x \equiv 1~~ (\bmod~2) ~. \\
~~~ \\
\frac{x}{2} & \mbox{if}  ~~x \equiv 0~~ (\bmod~2)~.
\end{array}
\right.
$$
in which $q>1$ and $r \ge 1$ are both odd integers. He
actually considers iteration of the  map $C_{q,r}(\cdot)$ acting
on the domain of positive odd integers, given by
$$
C_{q,r}(x) = \frac{qx+r}{2^{e_2(qx+r)}},
$$
where $e_2(x)$ denotes the highest power of $2$ dividing $x$.

Most results of the paper concern the map $C_{3,1}(\cdot)$
corresponding to the $3x+1$ map.
He first presents a heuristic probabilistic argument why iterates
of $C_{3,1}(\cdot)$  should decrease at an exponential rate, 
based on this he formulates a conjecture  that the
number of steps $H(x)$ starting from $x$
needed to reach $1$ under iteration of $C_{3,1}(\cdot)$
should be approximately
$H(x) \approx \frac{\log x}{\log \frac{16}{9}}$
for most integers.
He proves that the number of odd integers $n$ taking exactly $h$
steps to reach $1$ is at least $\frac{1}{h!} (\log_2 x)^h$.
He deduces  that the function $\pi_1(x)$ which counts
the number of odd integers below $x$ that eventually reach
$1$ under iteration of $C_{3,1}(\cdot)$  has
$\pi_1(x) > x^c$ for a positive constant $c$. (He does not
compute its value, but his proof seems to give  $c= 0.05$.)
He shows there are no cycles of length less than 17985
aside from the trivial cycle, using approximations to
$\log_2 3$. 

Concerning the ``$qx+r$'' problem, he formulates the
conjecture that, aside from $(q, r) = (3,1)$, 
every map $C_{q,r}(\cdot)$ has
at least one orbit that never visits $1$.
He proves that this conjecture is true whenever $r \ge 3$, and
in the remaining case $r=1$ he proves it for $q = 5$, $q=181$ and $q=1093$.
For the first two cases he exhibits a periodic orbit not
containing $1$, while for $q=1093$ he uses
the fact that there are no numbers of height $2$ above $1$,
based on  the congruence $2^{q-1} \equiv 1~ (\bmod~ q^2)$.
(This last argument would apply as well to  $q= 3511$.)
He argues the conjecture is true in the remaining cases because
a heuristic probabilistic argument suggests that for each $q \ge 5$ the
``$qx+1$'' problem should have a divergent trajectory.

\item
J. Leslie Davison (1977),
{\em Some Comments on an Iteration Problem},
Proc. 6-th Manitoba Conf. On Numerical Mathematics,
and Computing (Univ. of Manitoba-Winnipeg 1976),
Congressus Numerantium XVIII, Utilitas Math.:
Winnipeg, Manitoba 1977,  pp. 55--59.
(MR 58 \#31773). \\
\newline
\hspace*{.25in}
The author considers iteration of the map $f: \ZZ^{+} \to \ZZ^{+}$ given
by 
$$
f(n) = \left\{
\begin{array}{cl}
\frac{3n+1}{2}  & \mbox{if} ~~n \equiv 1~~ (\bmod ~ 2), n> 1 \\
~~~ \\
\frac{n}{2} & \mbox{if}~~ n \equiv 0~~ (\bmod ~2 )~ \\
~~~\\
1  & \mbox{if}~~ n=1,
\end{array}
\right.
$$
This is essentially the $3x+1$ function, except for $n=1$. 
He calls a sequence of iterates  of this map  a {\em circuit}
if it starts with an odd number $n$, produces a sequence of odd
numbers followed by a sequence of even numbers, ending at
an odd number $n^{*}$. A circuit is a {\em cycle} if
$n= n^{*}.$  

Based on computer evidence, he conjectures that the
number of circuits required during the iteration of a
number $n$ to $1$ is at most $K \log n$, for some
absolute constant $K$. He presents a probabilistic
heuristic argument in support of this conjecture.

He asks whether a circult can ever be a cycle.
He shows that this
question can be formulated as the exponential
Diophantine equation: There exists
a circuit that is a cycle if and
only if there exist positive integers
$(k, l, h)$ satisfying
$$ 
(2^{k+l} - 3^k)h = 2^l - 1.
$$
(Here there are $k$ odd numbers and $l$ even numbers
in the cycle, and 
$h= \frac{n+1}{2^k}$ where $n$ is the smallest
odd number in the cycle.) 
The trivial cycle $\{1,2\}$ of the $3x+1$ map
corresponds to the solution $(k, l, h) = (1, 1, 1)$,
Davison states he has been unable to find any other solutions.
He notes that the $5x+1$ problem has a circuit that is a 
cycle. 

Steiner (1978) subsequently showed that $(1,1,1)$ is the only 
positive solution
to the exponential Diophantine equation above.

\item
Philippe Devienne, Patrick Leb\`egue, Jean-Chrisophe Routier (1993), 
{\em Halting Problem of One Binary Horn Clause is Undecidable,}
{\em Proceedings of STACS 1993,}
Lecture Notes in Computer Science No. {\bf 665},
Springer--Verlag 1993, pp. 48--57.
(MR 95e:03114). \\
\newline
\hspace*{.25in}
The halting problem for derivations using a 
single binary Horn clause for reductions is shown
to be undecidable,
by encoding Conway's  undecidability result on iterating periodically linear
functions having no constant terms, cf. Conway (1972). 
In contrast, the problem of whether or not ground can be 
reached using reductions by
a single binary Horn clause is decidable.
[M. Schmidt-Schauss, Theor. Comp. Sci.
{\bf 59}
(1988), 287--296.]

\item
James M. Dolan, Albert  F. Gilman and Shan Manickam (1987), 
{\em A generalization of Everett's result on the Collatz $3x  +  1$ problem,}
Adv. Appl. Math. {\bf 8} (1987), 405--409. (MR 89a:11018). \\
\newline
\hspace*{.25in}
This paper shows that for any $k  \geq  1$, the set of 
$m  \in  \ZZ + $ having $k$
distinct
iterates $T^{(i)} (m)  <  m$ has density one.

\item
Richard Dunn (1973),
{\em On Ulam's Problem,}
Department of Computer Science, 
University of Colorado, Boulder, Technical Report CU-CS-011-73, 15pp.\\
\newline
\hspace*{.25in}
This report gives early computer experiments on the $3x+1$ problem.
The computation numerically verifies the $3x+1$ conjecture on a CDC 6400
computer up to 22,882,247. Dunn also  calculates the densities $F(k)$ defined
in equation (2.16) of  Lagarias (1985) for $k \le 21$.

\item
Peter Eisele and Karl-Peter Hadeler (1990),
{\em Game of Cards, Dynamical Systems, and a Characterization of the
Floor and Ceiling Functions,}
Amer. Math. Monthly {\bf 97} (1990), 466--477. 
(MR 91h:58086). \\
\newline
\hspace*{.25in}
This paper studies  iteration of the mappings 
$f(x) =  a+ \lceil \frac{x}{b} \rceil$ on $\ZZ$
where $a,b$ are positive integers.
These are periodical linear functions $(\bmod~b)$.
For $b  \geq 2$, every
trajectory becomes eventually constant or reaches a cycle of order 2.

\item
Shalom Eliahou (1993), 
{\em The $3x +  1$ problem: New Lower Bounds on Nontrivial Cycle Lengths,} 
Discrete Math., {\bf 118}(1993), 45--56.
(MR 94h:11017). \\
\newline
\hspace*{.25in}
The author shows that any nontrivial cycle 
on $\ZZ + $ of the $3x  +  1$ function
$T(x)$ has period $p =  301994A  +  17087915B  +  85137581C$ with
$A, B, C$ nonnegative integers where $B  \geq  1$,
and at least one of $A$ or $C$ is zero.
Hence the minimal possible period length is at least 17087915.
The method uses the continued fraction expansion of $\log_2  3$, 
and the truth of the $3x  +  1$ Conjecture for all $n  <  2^{40}$.
The paper includes   a table of 
partial quotients and convergents to the continued fraction of $\log_2 3$.

\item
Peter D. T. A. Elliott (1985),
{\em Arithmetic Functions and Integer Products,}
Springer-Verlag, New York 1985. (MR 86j:11095) \\
\newline
\hspace*{.25in}
An {\em additive function}
is a function with domain $\ZZ^+ $, which satisfies $f(ab) =  f(a) + f(b)$
if $(a,b) =  1$.
In Chapters\ 1--3 Elliott studies additive functions having the property
that $| f(an +  b)  -  f(An  +  B) |  \leq  c_0$
for all $n  \geq  n_0$,
for fixed positive integers $a,b,A,B$ with
$\det \left[ \begin{array}{c}
a~~b \\
A~~B
\end{array} \right] \neq 0 $, and deduces that
$|f(n) |  \leq  c_1 ( \log n)^3$.
For the special case $\left[ \begin{array}{c}
a~~b \\
A~~B
\end{array}
\right] = \left[
\begin{array}{c}
1~~1 \\
1~~0
\end{array}
\right]$
an earlier argument of Wirsching yields a bound
$|f (n) |  \leq  c_2 ( \log n )$.
On page 19 Elliott indicates that the analogue of Wirsching's argument
for $\left[ \begin{array}{c}
a~~b \\
A~~B
\end{array}
\right] = \left[
\begin{array}{c}
3~~1 \\
1~~0
\end{array}
\right]$
leads to the $3x  +  1$ function, and implies that
$|f(n) |  \leq  c_3  \sigma_\infty (n)$.
A strong form of the $3x  +  1$ Conjecture claims that
$\sigma_\infty (n)  \leq  c_4  \log n$, see Lagarias and Weiss (1992).
Elliott proves elsewhere by other arguments that in fact 
$| f(n) | \leq c_4  \log  n$ holds.
[P. D. T. A. Elliott, J. Number Theory {\bf 16} (1983), 285--310.]

\item
Paul Erd\H{o}s and R. L. Graham (1979),
{\em Old and new problems and results in combinatorial number theory:
van der Waerden's theorem and related topics,}
Enseign. Math. {\bf 25} (1979), no. 3-4, 325--344.   (MR 81f:10005( \medskip

%\newline
%\hspace*{.25in}
This list of many problems includes problems raised in Klarner and Rado (1974),
Hoffman (1976), and Hoffman and Klarner (1978), (1979) on the smallest set of nonnegative
integers obtained from a given set $A$ under iteration of a finite set $R$ of
functions $\rho(x_1, ..., x_r)= m_0 + m_1 x_1 + \cdots + m_r x_r$, with
nonnegative integer $m_i$ for $i \ge 1$. 
Denoting this set $< R: A>$, one can ask for the size and structure of this
set. In the case of one variable functions $R= \{ a_1 x + b_1, ..., a_r x+ b_r\}$
Erd\H{o}s showed (see Klarner and Rado (1974)) that if $\sum \frac{1}{a_i}<1$,
then the set has density $0$. The case when $\sum \frac{1}{a_i} =1$
is pointed to as a source of unresolved problems. Erd\H{o}s had proposed as
a prize problem: For $R = \{ 2x+1, 3x+1, 6x+1\}$  and $A= \{ 1\}$ is the set $<R: A>$
of positive density? This was answered in the negative by D. J. Crampin and A. J. W. Hilton 
(unpublished), as summarized  in Klarner (1982) and Klarner (1988).

\item
Paul Erd\H{o}s and R. L. Graham (1980),
{\em Old and new problems and results in combinatorial number theory}
Monographie No. 28 de L'Enseignement Math\'{e}matique, 
Kundig: Geneva 1980. \\
\newline
\hspace*{.25in}
This book includes Erd\H{o}s and Graham (1979) as one chapter. 
Thus it  includes problems raised in Klarner and Rado (1974),
Hoffman (1976), and Hoffman and Klarner (1978), (1979) on the smallest set of nonnegative
integers obtained from a given set $A$ under iteration of a finite set $R$ of
functions $\rho(x_1, ..., x_r)= m_0 + m_1 x_1 + \cdots + m_r x_r$, with
nonnegative integer $m_i$ for $i \ge 1$. 
Denoting this set $< R: A>$, one can ask for the size and structure of this
set. In the case of one variable functions $R= \{ a_1 x + b_1, ..., a_r x+ b_r\}$
Erd\H{o}s showed (see Klarner and Rado (1974)) that if $\sum \frac{1}{a_i}<1$,
then the set has density $0$. The case when $\sum \frac{1}{a_i} =1$
is pointed to as a source of unresolved problems. Erd\H{o}s had proposed as
a prize problem: For $R = \{ 2x+1, 3x+1, 6x+1\}$  and $A= \{ 1\}$ is the set $<R: A>$
of positive density? This was answered in the negative  by D. J. Crampin and A. J. W. Hilton 
(unpublished), as summarized  in Klarner (1982) and Klarner (1988).

\item
C. J. Everett (1977),
{\em Iteration of the number theoretic function $f(2n)=n, f(2n+1)=3n+2$},
Advances in Math. {\bf 25} (1977), 42--45.
(MR 56\#15552). \\
\newline
\hspace*{.25in}
This is one of the first research papers 
specifically on the $3x+1$ function.
Note that  $f(\cdot)$ is the $3x+1$-function $T(\cdot)$.
The author shows that the set of positive integers $n$ having some iterate
$T^{(k)}(n) < n $ has natural density one. 
The result was obtained independently and
contemporaneously by Terras (1976).

\item
Carolyn Farruggia, Michael Lawrence and Brian Waterhouse (1996),
{\em The elimination of a family of periodic parity vectors
in the $3x+1$ problem,}
Pi Mu Epsilon J. {\bf 10} (1996), 275--280. \\
\newline
\hspace*{.25in}
This paper shows that the parity vector $1 0^k$ is
not the parity vector of any integral periodic orbit
of the $3x+1$ mapping whenever $k \ge 2$.
(For $k=2$ the orbit with parity vector $10$ is the integral orbit
$\{ 1, 2\}$. )

\item
Marc R. Feix,  Amador  Muriel, Danilo Merlini, and Remiglio Tartini (1995),
{\em The $(3x +  1)/2$ Problem: A Statistical Approach,}
in: {\em Stochastic Processes, Physics and Geometry~II,}
Locarno 1991. (Eds: S. Albeverio, U. Cattaneo, D. Merlini)
World Scientific, 1995, pp.~289--300. \\
\newline
\hspace*{.25in}
This paper formulates heuristic stochastic models
imitating various behaviors of the $3x  +  1$ function, and compares
them to some data on the $3x+1$ function. In Sect. 2
it describes a random walk model imitating "average" behavior of 
forward iterates of the
$3x+1$ function. In Sect. 3 it examines trees of inverse iterates of this function,
and predicts that the number of leaves at level $k$ of the tree should
grow approximately like $A (\frac{4}{3})^k$ as $k \to \infty$. In
Sect. 4 it describes computer methods for rapid testing of the
$3x+1$ Conjecture. In Sect. 5 it  briefly considers
related functions. In particular, it considers the $3x-1$ function and
the function
$$
\tilde{T}(x) = \left\{
\begin{array}{cl}
\df{x}{3}  & ~~\mbox{if} ~~x \equiv 0~~ (\bmod ~ 3) \\
~~~ \\
\df{2x+1}{3}  &~~ \mbox{if} ~~x \equiv 1~~ (\bmod ~ 3) \\
~~~ \\

\df{7x+1}{3} &~~ \mbox{if}~~ x \equiv 2~~ (\bmod ~3 )~.
\end{array}
\right.
$$
Computer experiments show that  the trajectories of $\tilde{T}(x)$ for $1\le n \le 200,000$ 
all reach the fixed point $\{1\}$. 
 
{\em Note.}  Lagarias (1992) and Applegate
 and Lagarias (1995c) study more detailed stochastic models 
 analogous  to those given here in Sects. 2 and 3,  respectively.
 %The number $A$ above seems  more accurately  modelled by a random variable, 
 %rather than  a constant.
 %with a rigorous analysis. 

\item
Marc R. Feix, Amador Muriel and Jean-Louis Rouet (1994), 
{\em Statistical Properties of an Iterated Arithmetic Mapping, }
J. Stat. Phys. {\bf 76} (1994), 725--741.
(MR 96b:11021). \\
\newline
\hspace*{.25in}
This paper interprets the iteration of the $3x+1$ map as
exhibiting a ``forgetting'' mechanism concerning the iterates
$(\bmod ~2^k )$, i.e. after $k$ iterations starting from
elements it draws from a fixed residue class $(\bmod ~2^k )$,
the iterate $T^k (n)$ is uniformly distributed
$(\bmod ~2^k )$.
It proves that certain associated $2^k \times 2^k$ matrices $M_k$
has $(M_k )^k =  J_{2^k}$ where $J_{2^k}$ is the
doubly-stochastic $2^k \times 2^k$ matrix having all
entries equal to $2^{-k}$.

\item
Piero Filipponi (1991), 
{\em On the $3n+1$ Problem: Something Old, Something New,}
Rendiconti di Mathematica, Serie VII, Roma {\bf 11} (1991), 85--103.
(MR 92i:11031). \\
\newline
\hspace*{.25in}
This paper derives  by elementary methods various facts about coalescences 
of trajectories and
divergent trajectories.
For example, the smallest counterexample $n_0$ to the $3x+1$
Conjecture, if one exists, must have $n_0  \equiv  7, 15, 27, 31$,
$39, 43, 63, 75, 79, 91 (\bmod ~96)$.
[The final Theorem 16 has a gap in its proof, because formula (5.11) 
is not justified.]

\item
Leopold Flatto  (1992) ,
{\em $Z$-numbers and $\beta$-transformations,}
in: {\em Symbolic dynamics and its applications (New Haven, CT, 1991)},
Contemp. Math. Vol. 135, American Math. Soc., Providence, RI 1992, 181--201.(MR94c:11065).\\
\newline
\hspace*{.25in}
This paper concerns the $Z$-number problem of Mahler (1968).
A real number $x$ is a {\em $Z$-number } if
$0 \le \{ \{ x (\frac{3}{2})^n\}\} < \frac{1}{2}$ holds for all $n \ge 0$, where $\{\{ x\}\}$ denotes  the
fractional part of $x$. Mahler  showed that there is at most one $Z$-number in
each unit interval $[n, n+1)$, for positive integer $n$, and bounded the number
of such $1 \le n \le X$ that can have  a $Z$-number by $X^{0.7}.$
This paper applies  the $\beta$-transformation of W. Parry [Acta Math.
Acad. Sci. Hungar. {\bf 11} (1960), 401-416] to get an improved upper  bound on  the number
of $n$ for which a $Z$-number exists. It studies symbolic
dyanmics of this transformation for  $\beta= \frac{3}{2}$, and deduces that
the number of $1\le n \le X$ such that there is a $Z$-number in $[n, n+1)$
is at most $X^{\theta}$ with $\theta = \log_2(3/2)\approx 0.59$. 
The paper also obtains related
results for more general $Z$-numbers associated to fractions $\frac{p}{q}$
having $q< p < q^2$.

\item
Zachary M.  Franco (1990), 
{\em Diophantine Approximation and the $qx  +  1$ Problem,}
Ph.D. Thesis, Univ. of Calif. at Berkeley 1990.
(H. Helson, Advisor). \\
\newline
\hspace*{.25in}
This thesis considers iteration of the $qx  +  1$ function defined by
$C_q (x) =  \frac{qx  +  1}{2^{{\rm ord}_2 (qx +1)}}$, where
$2^{{\rm ord}_2 (y)}  || y$,
and both $q$ and $x$ are odd integers.
The first part of the thesis studies a conjecture of Crandall (1978),
and the results appear in Franco and Pomerance (1995).
The second part of the thesis gives a method to determine for a
fixed $q$ whether there are any orbits of period 2, i.e. solutions of $C_q^{(2)} (x) =  x$,
and it shows that for $|q | < 10^{11}$, only
$q =  \pm 1$, $\pm 3$, $5, -11 ,  -91$, and 181 have such orbits.
The method uses an inequality of F. Beukers, [Acta Arith.
{\bf 38}
(1981) 389--410].

\item
Zachary Franco and Carl Pomerance (1995),
{\em On a Conjecture of Crandall Concerning the $qx  +  1$ Problem,}
Math. Comp. {\bf 64} (1995), 1333--1336.
(MR9 5j:11019). \\
\newline
\hspace*{.25in}
This paper considers iterates of the $qx  +  1$ function 
$C_ q (x) = \frac{qx  +  1}{2^{{\rm ord}_2(qx  +  1)}}$, 
where $2^{{\rm ord}_2 (y)}  || y$ and
both $q$ and $x$ are odd integers.
Crandall (1978) conjectured that for each odd $q  \geq  5$ there is some
$n  >  1$
such that the orbit $\{ C_q^{(k)} (n) : k \geq  0 \}$
does not contain 1, and proved it for $q =  5,  181$ and 1093.
This paper shows that $\{ q :$ Crandall's conjecture is true for $q \}$
has asymptotic density\ 1, by showing the stronger result that the set
$\{ q : C^{(2)} (m) \neq 1$ for all $m  \in  \ZZ \}$
has asymptotic density one.

\item
Michael Lawrence Fredman (1972)
{\em Growth properties of a class of recursively defined functions,}
Ph. D. Thesis, Stanford University, June 1972, 81 pages\\
\newline
\hspace*{.25in}
 Let $g(n)$ be a given function.
 This  thesis discusses solutions of the general recurrence 
$M(0)= g(0)$,
 $$M(n+1) := g(n+1) + \min_{0 \le k \le n} \left( \alpha M(k) + \beta M(n-k) \right), $$
 in which $\alpha, \beta > 0$. 
It has three chapters and a conclusion.
 The first two chapters concern the special case where $g(n)=n$,
Chapter 3 considers more general cases. 
The author how the quantity $M(n)$ has an interpretation in terms of minimum total
weight of weighted binary trees having $n$ nodes. For the analysis in the $h(n)=n$ case.
Theorem 2.1 sets $D(n) = M(n) -M(n-1)$  and sets $h(x) = \sum_{\{j: D(j)) \le x\}} 1$,
and shows that $h(x)$ satisfies 
 $h(x)=1$ for $0 \le x < 1$ and the functional equation
$$
h(x) = 1+ h(\frac{x-1}{\alpha}) + h(\frac{x-1}{\beta}).
$$
It determines the growth rates of $h(x)$ and $M(x)$ in many circumstances.
We describe here results only for the case  $\min( \alpha, \beta) >1$.
Let $\gamma$ be the unique positive solution to $\alpha^{-\gamma} +\beta^{-\gamma} =1$.
It is shown that the function $h(x)$ has order of magnitude $x^{\gamma}$ while 
$M(x)$ has order of magnitude $x^{1+ \frac{1}{\gamma}}$. 
Theorem 2.3.2 states that  $\lim_{x \to \infty} h(x) x^{-\gamma}$ exists if and only
if $\lim_{x \to \infty} M(x) x^{-1-\frac{1}{\gamma}}$ exists. It is shown the limits always
exist if $\frac{\log \alpha}{\log \beta} $ is irrational, but in general do not exist 
when $\frac{\log \alpha}{\log \beta} $ is rational. Chapter 3 obtains less precise growth
rate information for a wide class of driving functions $g(x)$. 
Some of the proofs use complex analysis and Tauberian theorems for Dirichlet series.
The conclusion of the thesis states 
applications and open problems. One given application is to answer a question
raised by Klarner, see  Klarner and Rado (1974). 
Fredman shows that the set $S$
of integers obtained starting from $1$ and  iterating the affine maps $x \mapsto 2x+1$,
$x \mapsto 3x+1$ has density $0$, and in fact the number of such integers below $x$ is at
most $O( x^{\gamma}$ where $2^{-\gamma} + 3^{-\gamma} =1$.
 (Here $\gamma \approx 0.78788$.) This result  follows from an upper bound on 
 growth of $h(x)$
 above when $\alpha=2, \beta=3,  g(n)=n$.
 It uses the fact that the function  $h(x)$ has an interpretation as counting  the number
of elements $\le x$ in the multiset 
$\tilde{S} := \cup_{ j=0}^\infty S_j $ generated by initial element $S_0= \{1\} $ and
inductively letting $S_{j+1}$ being the image of $S_j$ under 
iteration of the two affine maps $x \mapsto \alpha x +1, \, x \mapsto \beta x +1$
(counting elements  in $\tilde{S}$ with the multiplicity they occur).

{\em Note.} 
This thesis also appeared as Stanford Computer Science Technical Report
STAN-CS-72-296. 
Some results of this thesis were subsequently published in
 Fredman and  Knuth (1974).

\item
Michael Lawrence Fredman  and Donald E. Knuth (1974)
{\em Recurrence relations based on minimization,}
J. Math. Anal. Appl. {\bf 48} (1974), 534--559 (MR 57\#12364). \\
\newline
\hspace*{.25in}
This paper studies the asymptotics of solutions of the general recurrence
$M(0)= g(0)$,
 $$M(n+1) := g(n+1) + \min_{0 \le k \le n} \left( \alpha M(k) + \beta M(n-k) \right), $$
for various choices of $\alpha, \beta, g(n)$. They denote this $M_{g \alpha \beta}(n)$
In \S2-\S4 they   treat the case $g(n)=n$, where they develop an interpretation of
this quantity in terms of weighted binary trees.; $M_{g\alpha \beta}$ is the minimum
total weight of any rooted binary tree with $n$ nodes. The weight of  a nodes in a 
finite binary tree $T$ is
given by assigning the root node $\sigma= \epsilon$
the weight  $w(\emptyset)=1$, and then inductively defining
$w(L\sigma)= 1+ \alpha w(\sigma), w(R\sigma) = 1+ \beta w(\sigma);$
for example $w(LRR) = 1 + \alpha + \alpha \beta + \alpha \beta^2.$ The total weight function $\sM(T)$
of a tree is the sum of  the weights of all nodes in it. They set
$M(n) := \min_{T: |T|=n}  \sM(T),$ and show this 
quantity is $M_{g \alpha \beta} (n)$ for $g(n)=n$. 
In \S4 they analyze the asymptotic behavior  of $M(n)$ in the case $\min(\alpha, \beta) >1$. 
They let $H(x)=h(x)+1$ where $h(x)$
counts the number of node weights $w(\sigma) \le x$ and observe that
it satisfies $H(x)=1$ for $0 \le x < 1$ and the functional difference equation
$$
H(x) = H(\frac{x-1}{\alpha}) + H(\frac{x-1}{\beta}).
$$
They note that $H(x)$ and $M(x)$ are 
related using the fact that a depth $n$ node has  weight $w(\sigma_n)\le x$
if and only if $H(x) > n$. Now consider the case  $\alpha, \beta > 1$,
and let $\gamma$ be
the unique positive solution to $\alpha^{-\gamma} + \beta^{-\gamma}=1.$
Lemma 4.1 shows that $H(x)$ is on the order of $x^{\gamma}$, and
$M(x)$ is on the order of $x^{1+ \frac{1}{\gamma}},$ so that  one can write
$H(x) = c(x) x^{\gamma}$, $M(x)= C(x) x^{1+ \frac{1}{\gamma}}$, where
$c(x)$ and $C(x)$ are positive bounded functions. 
The asymptotic
behaviors of $H(x)$ and $M(x)$ now depends on properties of the 
positive real number $\frac{\log \alpha}{\log \beta}$. 
 Theorem 4.1 shows
that when  $\frac{\log \alpha}{\beta}$ is rational, 
the function
$C(x)$ is usually an oscillatory function having  no limiting value at $\infty$.
Theorem 4.3  shows that
when $\frac{\log \alpha}{\beta}$ is irrational, $C(x)$ has a positive limiting value as $x \to \infty$
so that  $M(x) \sim C x^{1+ \frac{1}{\gamma}}$. (Similar results hold for $c(x)$; this
was explicitly shown in  Fredman (1972).) This paper also analyzes various cases where
$0 < \alpha,  \beta <1$ and the growth rate is polynomial. 
Some proofs use generating functions,
others use complex
analysis and Tauberian theorems.
See also Pippenger (1993) for proofs of some results by more elementary methods.

\item
David Gale (1991), 
{\em Mathematical Entertainments: More Mysteries,} 
Mathematical Intelligencer {\bf 13}, No. 3, (1991), 54--55. \\
\newline
\hspace*{.25in}
This paper discusses the possible undecidability of the 
$3x +  1$ Conjecture, and also whether the orbit
containing 8 of the original Collatz function
$$
U(n) = \left\{
\begin{array}{ll}
\df{3}{2} n & \mbox{if}~n \equiv 0~~(\bmod~2) \\
~~~ \\
\df{3}{4} n + \df{1}{4} & \mbox{if}~n \equiv 1~~(\bmod~4 ) \\
~~~ \\
\df{3}{4} n -  \df{1}{4} & \mbox{if}~n \equiv 3~ (\bmod ~4)
\end{array}
\right.
$$
is infinite.

\item
Guo-Gang Gao (1993),
{\em On consecutive numbers of the same height in the Collatz problem,}
Discrete Math., {\bf 112} (1993), 261--267.
(MR 94i:11018). \\
\newline
\hspace*{.25in}
This paper proves
 that if there exists one $k$-tuple of consecutive integers all having
the same height and same total stopping time, then there exists infinitely
many such $k$-tuples.
(He attributes this result to P. Penning.)
There is a 35654-tuple starting from $2^{500}  +  1$.
He conjectures that the set $\{  n : C^{(k)} (n) =  C^{(k)}(n  +  1)$ for some
$k  \leq  \log_2 n \}$ has natural density one, and proves that it has a
natural density which is at least 0.389.

\item
Manuel V. P.  Garcia and Fabio A. Tal (1999),
{\em A note on the generalized $3n + 1$ problem,}
Acta Arith. {\bf 90}, No. 3 (1999), 245--250.
(MR 2000i:11019). \\
\newline
\hspace*{.25in}
This paper studies the generalized $3x + 1$ function,
defined for $m > d \geq 2$ by
$$
H(x) = \left\{
\begin{array}{cl}
\df{x}{d} & \mbox{if}~~x \equiv 0 ~~ (\bmod ~d)~, \\
~~\\
\df{mx - \pi (ma)}{2} & \mbox{if} ~~x \equiv a ~~ (\bmod ~d),  ~a \neq 0,
\end{array}
\right.
$$
where $\pi(x)$ denotes projection (mod d) onto a
fixed complete set of residues (mod d). The
{\em Banach density} of a set $B \subset \ZZ^+$ is
$$\rho_b(B) = \limsup_{n \to \infty} (max_{a \in \ZZ^+} 
\frac{\sharp(B \cap~\{a, a+1, ... ,a + n - 1\})}{n}).$$
The Banach density is always defined and is at least as large
as the natural density of the set $B$ , if it exists.
Call two integers $m_1$ and $m_2$ {\em equivalent} if
there is some positive integer $k$ such that 
$H^{(k)}(m_1) = H^{(k)}(m_2)$. The authors assume that
$m < d^{d/{d - 1}}$, a hypothesis which implies
that almost all integers have some iterate which is smaller,
and which includes the $3x +1$ function as a special case.
They prove that if $\sP$ is any
complete set of representatives of equivalence classes of $\ZZ^+$
then the Banach density of  $\sP$ is zero. As a corollary they
conclude  that
the Banach density of the orbit of any integer $n$ under 
such a map $H$
is zero. In particular, the Banach density of any divergent
trajectory for such a map is zero.

\item
Martin Gardner (1972),
{\em Mathematical Games},
Scientific American {\bf 226} No. 6, (June 1972), 114--118. \\
\newline
\hspace*{.25in}
This article is one
of the first places the $3x+1$ problem is stated
in print. Gardner attributes the problem to a
technical report issued 
 by M. Beeler, R. Gosper, and R. Schroeppel, HAKMEM, Memo 239,
Artificial Intelligence Laboratory, M.I.T., 1972, p. 64. 
The $3x+1$ problem certainly predates this memo of Beeler et al.,
which is a collection of problems. 

\item
Lynn E. Garner (1981),
{\em On the Collatz $3n+1$ algorithm},
Proc. Amer. Math. Soc. {\bf 82} (1981), 19--22.
(MR 82j:10090). \\
\newline
\hspace*{.25in}
The {\em coefficient stopping time} $\kappa(n)$
introduced by Terras (1976) is the least
iterate $k$ such that
$T^{(k)}(n) = \alpha(n) n + \beta(n)$, with $\alpha(n) < 1$.
Here $\alpha(n) = \frac{3^{a(n)}}{2^n}$
where $a(n)$ is the number of iterates $T^{(j)}(n) \equiv 1 (\bmod~2)$
with $0 \le j < k$. One has  $\kappa(n) \le \sigma(n)$,
where $\sigma(n)$ is the stopping time of $n$, and the 
Coefficient Stopping Time Conjecture of Terras (1976) asserts that
$\kappa(n) = \sigma(n)$ for all $n \ge 2$. 
This paper proves
that $\kappa(n) < 105,000$ implies that $\kappa(n) \le \sigma(n)$.
The proof methods used are those of  Terras (1976),
who proved the conjecture holds for $\kappa(n) < 2593.$
They invove the use of the 
continued fraction expansion of $\log_2 3$
and the truth of the $3x+1$ Conjecture for $n < 2.0 \times 10^9$.

\item
Lynn E. Garner (1985), 
{\em On heights in the Collatz $3n  +  1$ problem, }
Discrete Math. {\bf 55} (1985), 57--64.
(MR 86j:11005). \\
\newline
\hspace*{.25in}
This paper shows that infinitely many pairs of consecutive integers have
equal (finite) heights and equal total stopping times.
To do this he studies how trajectories of consecutive integers
can coalesce.
Given two consecutive integers $m$, $m  +  1$ having
$T^{(i)} (m) \neq T^{(i)} (m  +  1)$ for $i  <  k$ and
$T^{(k)} (m) =  T^{(k)} (m  +  1)$, associate to them the pair
$( \bv ,  \bv' )$
of $0  -  1$ vectors of length $k$ encoding the parity of $T^{(i)} (m)$
(resp. $T^{(i)} (m  +  1)$) for
$0  \leq  i \leq k  -  1$.
Call the set $\sA$ of pairs
$( \bv ,  \bv' )$ obtained this way {\em admissible pairs}.
Garner exhibits collections $\sB$ and $\sS$ of pairs of equal-length
$0  -  1$ vectors $( \bb ,  \bb')$ and $( \bs ,  \bs' )$ called
{\em blocks}
and
{\em strings},
respectively, which have the properties:
If $( \bv ,  \bv' )  \in  \sA$ and
$( \bb ,  \bb')  \in  \sB$
then the concatenated pair $( \bb \bv ,  \bb'\bv' )  \in \sA$,
and if $( \bs ,  \bs' )  \in  \sS$ then
$( \bs \bv' ,  \bs' \bv )  \in  \sA$.
Since $( 001, 100)  \in  \sA$, $(10, 01)  \in  \sB$
and $(000011,  101000)  \in \sS$, the set $\sA$ is infinite.
He conjectures that: (1)\ a majority of all positive integers
have the same height as an adjacent integer (2) arbitrarily long runs of
integers of the same height occur.

\item
Wolfgang Gasch\"{u}tz (1982),
{\em Linear abgeschlossene Zahlenmengen I.} [Linearly closed number sets I.]
J. Reine Angew. Math. {\bf 330} (1982), 143--158. (MR 83m:10095).  \\
\newline
\hspace*{.25in}
This paper studies subsets $S$ of $\ZZ$ or $\NN$ closed
under iteration of an affine map  $f(x_1, \cdots, x_r) = w_0 + w_1x_1 +w_2x_2+ \cdots + w_r x_r$ with
integer $w_0, w_1, ..., w_r$. 
He shows for a general function one
can reduce analysis to the case $w_0=1$
and initial seed value $0$. A polynomial $f(x_1, ... , x_r) \in \NN[x_1, ..., x_r]$
is {\em controlled} (``gebremst") if $f(x_1, ..., x_r) \preccurlyeq  (x_1+ \cdots + x_r) f(x_1, ..., x_r) + f(0,0,..., 0),$
where $f \preccurlyeq g$ means each coefficient of $f$ is no larger than the corresponding
coefficient of $g$. It is {\em $m$-controlled} if $f(0,0,  \cdots, 0) \le m.$ Let $F_{r, m}$ denote
the set of $m$-controlled polynomials in $r$ variables, and let
$F_{r, m} (w_1, ..., w_r) =\{ f(w_1, ..., w_r) : f \in F_{r,m}\}$.
Theorem 4.1 asserts that for $f := 1+ w_1 x_1 +\cdots + w_rx_r$ the smallest set
containing $0$ and closed under iteration of $f$ is $F_{r,1}(w_1, w_2, ..., w_r)$. He also shows that if
the greatest common divisor of $(w_1, ..., w_r)=1$ then (i) if $w_i$ are nonnegative then 
 $F_{r, m}(w_1, .., w_r) = \NN$ for all large enough $m$; (ii) if some $w_i$ is negative
and some other $w_i$ nonzero then for  $F_{r, m}(w_1, .., w_r) = \ZZ$
for large enough $m$.
He applies these results to obtain a  characterization of  those functions of two variables
 $f(x_1, x_2) = w_0 + w_1 x_1 + w_2 x_2$
such that 
 the smallest set of integers
containing $0$ and closed under its action consists of all nonnegative integers $\NN$.
They are exactly those functions with $w_0=1$,  $w_1, w_2 \ge 0$, 
with greatest  common divisor $(w_1, w_2)=1$, having the extra property that
iteration $(\bmod~ w_1w_2)$ visits  all residue classes $(\bmod~ w_1w_2)$.
He also obtains a criterion  for the smallest set to be $\ZZ$.
He notes that similar results were  obtained earlier by Hoffman and Klarner (1978), (1979). 

%\medskip
{\em Note.} This study was motivated by the author's
 earlier work on single-word criteria for subgroups of abelian groups: 
W. Gasch\"{u}tz, {\em Untergruppenkriterien f\"{u}r abelsche Gruppen}, Math. Z. {\bf 146} (1976), 89--99.

\item
H. Glaser and Hans-Georg Weigand (1989), 
{\em Das-ULAM Problem-Computergest\"utze Entdeckungen,}
DdM (Didaktik der Mathematik) {\bf 17}, No. 2 (1989), 114--134. \\
\newline
\hspace*{.25in}
This paper views  the $3x+1$ problem as iterating the Collatz function,
and views it as an algorithmic problem between mathematics
and computer science. It views study of this problem as useful
as training in exploration of mathematical ideas. It formulates 
exploration as a series of questions to ask
about it, and answers some of them. 
Some of these concern properties of the trees of inverse iterates of
the Collatz function starting from a given number. It proves
branching properties of the trees via congruence properties
modulo powers of $3$. 
This paper also discusses programming the $3x+1$ iteration in
the programming languages
Pascal and  LOGO.
%[I have not seen this paper.]

\item
Gaston Gonnet (1991), 
{\em Computations on the $3n+1$ Conjecture,} 
MAPLE Technical Newsletter  {\bf 0}, No. 6, Fall 1991.\\
\newline
\hspace*{.25in}
This paper describes
 how to write computer code to efficiently compute $3x+1$ function
iterates for very large $x$ using MAPLE.
It displays a computer plot of the total stopping function for $n  <  4000$, 
revealing an
interesting structure of well-spaced clusters of points.

\item
Richard K. Guy (1981),
{\em Unsolved Problems in Number Theory},
Springer-Verlag, New York 1981. 
[Second edition: 1994. Third Edition: 2004.] \\
\newline
\hspace*{.25in}
Problem E16 discusses  the $3x+1$ Problem.
Problem E17 discusses permutation sequences, includes Collatz's original permutation,
see Klamkin (1963). 
Problem E18 discusses Mahler's Z-numbers, see Mathler (1968).
 Problem E36 (in the second edition) discusses Klarner-Rado sequences,
see Klarner and Rado (1974).

{\em Note.} Richard Guy [private conversation] informed me that he first
heard of the problem in the early 1960's from his son Michael Guy,
who was a student at Cambridge University and friends with John Conway.
John Conway [private conversation] confirms that he heard of the problem 
and worked on it as a Cambridge undergraduate (BA 1959). 
\item
Richard K. Guy (1983a)
{\em Don't try to solve these problems!},
Amer. Math. Monthly {\bf 90} (1983), 35--41. \\
\newline
\hspace*{.25in}
The article gives some brief history of work
on the $3x+1$ problem. It mentions at second hand
a statement of P. Erd\H{o}s
regarding the $3x+1$ problem: ``Mathematics is not yet ripe enough 
for such questions.'' 

The $3x+1$ problem is stated as Problem 2. 
Problem 3 concerns cycles in the  original Collatz problem,
for which see Klamkin (1963). 
Problem 4 asks the question due to Klarner (1982): Let S be the smallest set
of positive integers containing $1$ which is
closed under $x \mapsto 2x$, $x \mapsto 3x+2$, $x  \mapsto 6x+3$.
Does this set have a positive lower density?

\item
Richard K. Guy (1983b),
{\em Conway's prime producing machine},
Math. Magazine {\bf 56} (1983), 26--33.
(MR 84j:10008). \\
\newline
\hspace*{.25in}
This paper gives a function  $g(\cdot)$ of
the type in Conway (1972) having the 
follwing property. If $p_j$ denotes the
$j$-th prime, given in increasing order,
then starting from the value $n= 2^{p_j}$ 
and iterating under $g(\cdot)$, the first power
of $2$ that is encountered in the iteration
is $ 2^{p_{j+1}}$. 

He shows that the associated register machine
uses only four registers. See also the paper
Conway (1987) on FRACTRAN. 

\item
Richard K. Guy (1986), 
{\em John Isbell's Game of Beanstalk and John Conway's 
Game of Beans Don't Talk,} 
Math. Magazine {\bf 59} (1986), 259--269.
(MR 88c:90163). \\
\newline
\hspace*{.25in}
John Isbell's game of Beanstalk has two players alternately make
moves using
the rule
$$
n_{i+1} =  \left\{
\begin{array}{cl}
\df{n_1}{2} & \mbox{if} ~n_i \equiv 0~~(\bmod~2 ) ~, \\
~~~ \\
3n_i \pm 1 &  \mbox{if}~n_i \equiv 1~~(\bmod~2 ) ~,
\end{array}
\right.
$$
where they have a choice if $n_i$ is odd.
The winner is the player who moves to 1.
In Conway's game the second rule becomes
$\df{3n  \pm 1}{2^\ast }$,
where $2^\ast$ is the highest power of 2 that divides the numerator.
It is unknown whether or not there are positions from which
neither player can force a win.
If there are then the $3x  +  1$ problem must have a nontrivial cycle or a
divergent trajectory.

\item
Lorenz Halbeisen and Norbert Hungerb\"uhler (1997), 
{\em Optimal bounds for the length of rational Collatz cycles,} 
Acta Arithmetica {\bf 78} (1997), 227-239. (MR 98g:11025). \\
\newline
\hspace*{.25in}
The paper presents ``optimal'' upper bounds for the
size of a minimal element in a rational cycle of length $k$ for the
$3x+1$ function.
These estimates improve on Eliahou (1993) but currently do not
lead to a better linear bound for
nontrivial cycle length than the value 17,087,915 obtained by Eliahou.
They show that if the $3x+1$ Conjecture is verified for
$1 \leq n \leq 212,366,032,807,211$, which is
about $2.1 \times 10^{14}$, then the lower bound
on cycle length jumps to 102,225,496.

[The $3x+1$ conjecture now  verified up to
$1.9 \times 10^{17}$ by E. Roosendaal,
see the comment on Oliveira e Silva (1999).]

\item
Gisbert Hasenjager (1990), 
{\em Hasse's Syracuse-Problem und die Rolle der Basen,} in:
{\em Mathesis rationis. Festschrift f\"{u}r Heinrich Schepers,}
(A. Heinekamp, W. Lenzen, M. Schneider, Eds.) Nodus Publications: M\"{u}nster 1990
(ISBN\ 3-89323-229-X), 329--336. \\
\newline
\hspace*{.25in}
This paper, written for a philosopher's anniversary, 
comments on the complexity of the
$3x+1$ problem compared to its simple appearance.
It suggests looking for patterns in the iterates written in base 3 or base 27,
rather than in base 2.
It contains heuristic speculations and no theorems.

\item
Helmut Hasse (1975),
{\em Unsolved Problems in Elementary Number Theory},
Lectures at University of Maine (Orono), Spring 1975.
Mimeographed notes. \\
\newline
\hspace*{.25in}
Hasse  discusses the $3x+1$ problem on pp. 23--33.
He calls it the Syracuse (or Kakutani) algorithm.
He asserts that A. Fraenkel checked it for $n < 10^{50}$,
which is not the case (private communication with A. Fraenkel).
He states that Thompson has proved the Finite Cycles
Conjecture, but this seems not to be the case, as no
subsequent publication has appeared. 

He suggests a generalization $( \bmod~m)$ for $m > d \ge 2$
and in which the map is 
$$
T_d(x) := 
\frac{mx+f(r)}{d}   \mbox{if} ~~x \equiv r~~ (\bmod ~ d) 
$$
in which $f(r) \equiv -mr ~(\bmod~d)$. He gives a 
probabiliistic argument suggesting  that all orbits are
eventually periodic, when $m=d+1$. 

Hasse's circulation of the problem motivated
some of the first publications on it.  He proposed
a class of generalized $3x+1$ maps studied
in M\"{o}ller (1978) and Heppner (1978).

\item
Brian Hayes (1984), 
{\em Computer recreations: The ups and downs of hailstone numbers},
Scientific American {\bf 250}, No. 1 (January 1984), 10--16. \\
\newline
\hspace*{.25in}
The author introduces the $3x  +  1$ problem to a general audience under yet
another name --- hailstone numbers.

\item
Ernst Heppner (1978),
{\em Eine Bemerkung zum Hasse-Syracuse Algorithmus},
Archiv. Math. {\bf 31} (1978), 317--320.
(MR 80d:10007) \\
\newline
\hspace*{.25in}
This paper studies iteration of generalized $3+1$ maps
that belong to a class
formulated by H. Hasse. This class consists of 
maps depending on parameters of the form
$$
T(n) = T_{m, d, R}(n) := \left\{
\begin{array}{cl}
\df{mn+r_j}{d}  & \mbox{if} ~~n \equiv j~~ (\bmod ~ d),~1 \le j \le d-1 \\
~~~ \\
\df{x}{d} & \mbox{if}~~ n \equiv 0~~ (\bmod ~d )~.
\end{array}
\right.
$$
in which the parameters $(d, m)$ satisfy
$d \ge 2$, $gcd(m, d) =1$,
and  the set $R= \{ r_j: 1 \le j \le d -1\}$ has each
$r_j \equiv -mj~(\bmod~d)$.  
The qualitative behavior of iterates of these maps 
are shown to depend on the
relative sizes of $m$ and $d$.

Heppner proves that if $m < d^{d/{d-1}}$ then almost all
iterates get smaller, in the following quantitative sense:
There exist positive real numbers $\delta_1, \delta_2$ such
that for any $x > d$, and 
$N = \lfloor \frac{\log x}{\log d}\rfloor$,
there holds
$$
\#\{ n \le x:~ T^{(N)}(n) \ge n x^{-\delta_1} \} = O ( x^{1 - \delta_2}).
$$
He also proves that if $m > d^{d/{d-1}}$ then
almost all iterates get larger, in the following quantitative sense:
There exist positive real numbers $\delta_3, \delta_4$ such
that for any $x > d$, and 
$N = \lfloor \frac{\log x}{\log d}\rfloor$,
there holds
$$
\#\{ n \le x:~ T^{(N)}(n) \le n x^{\delta_3} \} = O ( x^{1 - \delta_4}).
$$
In these results the constants $\delta_j$ depend on $d$ and
$m$ only while the implied constant in the $O$-symbols depends
on $m, d$ and $R$.
This results improve on those of M\"{o}ller (1978).

\item
Dean G. Hoffman and David A. Klarner (1978),
{\em  Sets of integers closed under affine operators- the closure
of finite sets,}
Pacific Journal of Mathematics {\bf 78} (1978), No. 2, 337--344. 
(MR 80i:10075) \\
\newline
\hspace*{.25in}
This work extends work of Klarner and Rado (1974),
concerning sets of integers generated by  iteration of  a single
multi-variable affine function. 
It considers iteration of the affine map $f(x_1, ..., x_m) = m_1x_1 + \cdots + m_r x_r + c$
assuming that : (i) $r \ge 2$, (ii) each $m_i \ne 0$
and the greatest common divisor $(m_1, m_2, ..., m_r) =1$.
The author supposes that $T$ is a set of (not necessarily positive)
 integers that is closed under iteration
of $f$ in the sense that if $t_1, ..., t_r \in T$ then $f(t_1, ..., t_r) \in T$.
The main result, Theorem 12, states that if in addition 
$f(t, t, ..., t) > t$ holds for all $t \in T$, then the following two statements are
equivalent: (1) $T$ is a finite union of infinite arithmetic progressions, and
(2) T is generated by some finite set $A$ under iteration of the map $f$,
using $A$ as  a ``seed." 

\item
Dean G. Hoffman and David A. Klarner (1979),
{\em  Sets of integers closed under affine operators- the finite basis theorems,}
Pacific Journal of Mathematics {\bf 83} (1979), No. 1, 135--144. 
(MR 83e:10080) \\
\newline
\hspace*{.25in}
This paper strengthens the results in Hoffman and Klarner (1978)
concerning sets of integers generated by iteration of a single multi-variable
affine function. 
The same hypotheses (i), (ii) are imposed on the function $f$ as in
that paper. 
Without any further hypothesis, the authors now conclude that if $T$ is any
set closed under iteration by $f$, then $T$ is generated by a finite set
of  ``seeds" $A$, so we may write $T =\langle  f: A \rangle$ in the notation of 
Klarner and Rado (1974). The second main result (Theorem 13)
is the conclusion that either $T$ is a one-element set, or 
a finite union of one-sided infinite arithmetic progessions, bounded below,
or a finite union of two-sided infinite arithmetic progressions.

%\item
% Hong, B. Y.  (1986),
%{\em  $3x+1$},
%J. of Hubei Normal University, Natural Science Edition,
%[Hubei shi fan xue yuan xue bao. Zi ran ke xue ban] (1986), No. 1, 1--5.\\
%\newline
%\hspace*{.25in}
%[I have not seen this paper.]
\item
Hong, Bo Yang (1986),
{\em  About $3X+1$ problem} (Chinese),
J. of Hubei Normal University, Natural Science Edition,
[Hubei shi fan xue yuan xue bao. Zi ran ke xue ban] (1986), No. 1, 1--5.\medskip

%\newline
%\hspace*{.25in}
Theorem 1 shows that if there is a positive odd number such that the $3x+1$ 
conjecture fails for it, then there are infinitely many such odd numbers, and
 the smallest such number must belong to  one of the congruence classes
$7, 15, 27, 31~(\bmod \,32)$. Theorem 2 shows that if there is a positive
odd number such that the $3x+1$ conjecture fails for it, then there is such
an odd number that is $53~(\bmod \, 64)$.

\item
Stephen D. Isard and Harold M. Zwicky (1970),
{\em Three open questions in the theory of one-symbol Smullyan systems,}
SIGACT News, Issue No. 7,  1970, 11-19.\\
\newline
\hspace*{.25in}
Smullyan systems are described in R.  Smullyan, {\em First Order Logic,}
[Ergebnisse  der Math. Vol. 43, Springer-Verlag, NY 1968
(Corrected Reprint: Dover, NY 1995)] In this paper
open question 2 concerns the set  of integers that are generated by 
certain Smullyan systems with one symbol $x$.
% with rules 
%\begin{eqnarray*}
%| ~~~~&& \\
%x ~~~&\to& xxx\\
%xx ~| & \to & x
%\end{eqnarray*}
If we let $n$ label a string of $n$ $x$'s, then the system allows the two 
string rewriting operations
$f(n)= 3n$ and $g(4m+3 )= 2m+1$.  The problem asks
whether it is true that,  starting from
$n_0= 1$, we can reach every positive  number $n \equiv 1~(\bmod~3)$
by a sequence of such string rewritings. The authors note that
one can reach every such number except possibly those with $n \equiv 80~(\bmod~81),$ 
and that the first few numbers in this congruence class are reachable.
In  the reverse direction we are allowed instead to apply either of the rules 
$F(3n)=n$ or $G(2m+1)= 4m+3$, and we then wish to get from
an arbitrary $n \equiv 1~(\bmod ~3)$ to $1$. 
An undecidability result given at the end of the paper, using Minsky machines,
has some features in common with Conway (1972).

{\em Note.} This symbol rewriting problem
involves many-valued functions in both the forward and backward directions,
which are both linear on congruence classes to some modulus.
 Trigg et. al (1976)  cited this problem in a discussion of
the $3x+1$ problem, as part of its prehistory. 
The $3x+1$ problem has similar feature (cf. Everett (1977))
with the difference that   the $3x+1$ function is single-valued in the forwards direction and
many-valued only in the backwards direction.
 Michel (1993) considers
other rewriting rules which are functions in one direction.

\item
Frazer Jarvis (1989), 
{\em 13, 31 and the $3x  +  1$ problem, }
Eureka {\bf 49} (1989), 22--25. \\
\newline
\hspace*{.25in}
This paper studies the function $g(n) = h(n  +  1)  -  h(n)$,
where $h(n)$ is the height of $n$,
and observes empirically that $g(n)$ appears unusually often to be
representable as $13 x  +  31 y$ with small values of $x$ and $y$.
It offers a heuristic explanation of this observation in terms of 
Diophantine approximations to
$\log_6 2$.
Several open problems are proposed, mostly concerning $h(n)$.

\item
John  A. Joseph (1998), 
{\em A chaotic extension of the Collatz function to $\ZZ_2 [i]$, }
Fibonacci Quarterly {\bf 36} (1998), 309--317. (MR 99f:11026). \\
\newline
\hspace*{.25in}
This paper studies the function~1
$\sF : \ZZ_2 [i] \rightarrow \ZZ_2 [i]$ given by
$$
\sF ( \alpha ) = \left\{
\begin{array}{cl}
\df{\alpha}{2} & \mbox{if} ~ \alpha \in [0]~, \\
\df{3\alpha+1}{2} & \mbox{if}~ \alpha \in [1]~, \\
\df{3 \alpha+i}{2} & \mbox{if}~ \alpha \in [i]~,
\\
\df{3 \alpha +1 +i}{2} & \mbox{if}~ \alpha \in [1+ i]~,
\end{array}
\right.
$$
where $[ \alpha]$ denotes the equivalence class of $\alpha$
in $\ZZ_2 [i]/2 \ZZ_2 [i]$.
The author proves that the map $\tilde{T}$ is chaotic in the
sense of Devaney
[{\em A First Course in Dynamical Sysems: Theory and Experiment},
Addison-Wesley 1992].
He shows that $\tilde{T}$ is not conjugate to $T \times T$ via a
$\ZZ_2$-module isomorphism, but is topologically
conjugate to $T \times T$.
This is shown using an analogue $\tilde{Q}_\infty$ of the 
$3x+1$ conjugacy map $Q_\infty$ studied in
Lagarias (1985), Theorem~L and in
Bernstein and Lagarias (1996).

\item
Frantisek Kascak (1992),
{\em Small universal one-state linear operator algorithms},
in: Proc. MFCS '92, Lecture Notes in Computer Science No. 629,
Springer-Verlag: New York 1992, pp. 327--335. [MR1255147]\\
\newline
\hspace*{.25in}
A {\em one-state linear operator algorithm} (OLOA) with modulus $m$,
is specified by data $(a_r, b_r, c_r)$ a triple
of integers, for each residue class $r~0\bmod~m)$.
Given an input integer $x$, the OLOA  does the following in
one step. It finds
 $x \equiv r ~(\bmod~m)$ and based on the value of $r$, it 
  does the following. 
If $c_r=0$, the machine halts and outputs $x$, calling it
a final number, 
If $c_r \ne 0$ it computes 
to $(a_r x + b_r ~\mbox{DIV} ~c_r= \frac{a_r x + b_r}{c_r}$, and
outputs this value, provided
this output is a nonnegative  integer. If $\frac{a_r x + b_r}{c_r}$ is a 
non-integer  or is negative, the the number $x$ is called a terminal number,
and the machine stops. 
An OLOA $L$ computes  a partially defined function $m \to f_{L}(m)$ as follows. 
For  initial input value  $m$, the OLOA is iterated, and halts with  $f_{L}(m)=x$ if
$x$ is a final number. Otherwise $f_{L}(m)= \uparrow$ is viewed as undefined,  
and this occurs if the computation 
either reaches a terminal number or else runs for an infinite number of 
steps without stopping.

The author observes that the $3x+1$ problem can be encoded as
an OLOA $L$ with modulus $m=30.$ The increase of modulus of the
$3x+1$ function from $m=2$ to $m=30$ is made in order to 
 encode the $3x+1$ iteration arriving at  value $1$ as a halting state in the sense above. 
 The $3x+1$ function input $n$ is encoded as the integer $m= 5\cdot 2^{n-1}$ 
 to the OLOA. It is not known if the function $f_L$ computed by this
 OLOA has a recursive domain $D(f_L)$, but  if the $3x+1$ conjecture is true,
 then $D(F_L) = \NN$ will be recursive.
 
 The main result of the paper is that there exists an OLOA $L$ with modulus $m=396$
which is a universal OLOA, i.e. there exists a unary recursive function $d$
such that $d f_{L}$ is a universal unary recursive function.  Here a   unary partial
recursive function $u$  is
{\em universal} it there exists 
some binary recursive function $c$ such that $F_x(y)= u c(x,y)$
where $F_x(y)$ is an encoding of the $x$-th function in an encoding of all
unary partial recursive function, evaluated at input $y$.
In particular this machine $L$  has an unsolvable halting problem, i.e. 
the function $f_L$ has domain $D(f_L)$ which is not recursive. 
The encoding of the construction uses some ideas from Minsky machines, 
as in Conway (1972). The construction of Conway would give a universal OLOA with
a much larger modulus.

\item
Louis H. Kauffman (1995),
{\em Arithmetic in the form,}
Cybernetics and Systems {\bf 26} (1995), 1--57. [Zbl 0827.03033]\\
\newline
\hspace*{.25in}
This interesting paper describes how to do arithmetic in certain symbolic
logical systems developed by  G. Spencer Brown, {\em Laws of Form},
George Unwin \& Brown: London 1969. Spencer-Brown's  formal system starts from a primitive
notion of the additive void $0$ and the multiplicative void $1$,
and uses division of space by boundaries to create numbers,
which are defined using
certain transformation rules defining equivalent 
expressions. Kauffman argues that these rules give operations more primitive than addition
and multiplication. He develops  a formal natural number  arithmetic in this context,
specifying replacement rules for rearrangements of symbolic expressions,
and gives a proof of consistency of value (Theorem, p. 27).

In Appendix A he presents a second symbolic system, called {\em string arithmetic}.
He discusses the $3x+1$ problem in this context. This system has 
three kinds of  symbols $\ast, \langle, \rangle$, where in representing integers the angles occur only in 
matching pairs. The integer $N$ is represented in unary as a string
of $N$ asterisks. This integer
 has alternate expressions, in which  the angles encode multiplication by $2$.
 Addition is given by concatenation of expressions. 
The transformation rules between equivalent expressions are, if
$W$ is any string of symbols
\begin{eqnarray*}
\ast \ast &\longleftrightarrow & \langle \ast \rangle \\
\rangle \langle &\longleftrightarrow &  (\mbox{blank}) \\
\ast W &\longleftrightarrow &  W \ast
\end{eqnarray*}
A natural number can be defined as any string $W$ equivalent under these rules to some
string of asterisks. 
For examples, these rules give as equivalent representations of $N=4$:
$$
 \ast \ast \ast \ast = \langle \ast \rangle \langle \ast \rangle = \langle \ast \ast \rangle = 
\langle \langle \ast \rangle \rangle.
$$
He observes that the $3x+1$ problem is  particularly simple
to formulate in this system. He gives transformation rules defining the two steps in
the $3x+1$ function. 
If $N = \langle W\rangle$, then it is even, and then $\frac{N}{2}= W$, while if
$N =  \langle W\rangle \ast$,  then it is odd and $\frac{3N+1}{2} = \langle W \ast \rangle W$.
These give the transformation rules
$$
\begin{array}{ccccc}
(T1): &~~~~~& \langle W \rangle & \mapsto &W\\
~~~~~~~~~~\\
(T2): & ~~~~~&\langle W \rangle \ast & \mapsto & \langle W \ast \rangle W
\end{array}
$$
The $3x+1$ Conjecture is equivalent to the assertion:
starting from a string of  asterisks of any length, and using the
transformation rules plus the rules (T1) and (T2), one can reach the string
of one asterisk $\ast$. Kauffman says: "A more mystical reason for writing the
Collatz [problem] in string arithmetic is the hope that there is some subtle pattern
right in the notation of string arithmetic that will show the secrets of the iteration."

\item
David Kay (1972),
{\em An Algorithm for Reducing the
Size of an Integer,}  Pi Mu Epsilon Journal {\bf 4}  (1972), 338.\\
\newline
\hspace*{.25in}
This short note  proposes the $3x+1$ problem as a possible undergraduate research project.
The  Collatz map is presented, being denoted  $k(n)$, and it is noted that 
it is an unsolved problem to show that all iterates on the positive integers
go to $1$. It states this has been verified for all integers up to a fairly
large bound.

The proposed research project asks if one can find  integers $p, q, r$ with
$p, q \ge 2$,  having the property that a result analogous
to the $3x+1$ Conjecture  can be rigorously proved for the function 
$$
k(n) := \left\{
\begin{array}{cl}
\df{n}{p} & \mbox{if} ~~n \equiv 0~~ (\bmod~p) ~. \\
~~~ \\
qn+r& \mbox{if}  ~~x \not\equiv 0~~ (\bmod~p)~.
\end{array}
\right.
$$

{\em Note.} This project was listed as proposed by the Editor of Pi Mu Epsilon;
the Editor was  David Kay.

\item
Timothy P. Keller (1999),
{\em  Finite cycles of certain periodically linear functions},
Missouri J. Math. Sci. {\bf 11} (1999), no. 3, 152--157.
(MR 1717767) \\
\newline
\hspace*{.25in}
This paper is motivated by the original Collatz function
 $f(3n)=2n,$ $f(3n+1)=4n+1$, $f(3n+2)=4n+3$, 
which is a permutation of $\ZZ$,
see Klamkin (1963). The class of
periodically linear functions which are permutations.
were characterized by Venturini (1997). 
The author defines a {\em permutation of type V} to
be a periodically linear function $g_p$ defined for an
odd integer $p>1$ by 
 $g_p(pn+r) = (p + (-1)^r)n + r$
for $0 \le r \le p-1$. The permutation $g_3$  is
conjugate to the original Collatz function. This paper shows that 
for each $L \ge 1$ any  permutation
of type $V$ has only finitely many cycles of period $L$. 
The author conjectures that for each odd $p$ there is a
constant $s_0(p)$ such that if $L \ge s_0(p)$ is the minimal period
of a  periodic orbit of $g_p$, then $L$ is a denominator of a
convergent of the continued fraction
expansion of $\gamma_p :=\frac{\log (p-1) - \log p}{\log p - \log (p-1)}.$
For $p=7$ these denominators start 
$\Lambda(7) =\{ 1, 2, 13, 28, 265, 293, ...\}$. He finds that $g_7$
has an orbit of period 265 with starting value $n_1=1621$,
and an orbit of period 293 with starting value $n_2=293.$ 
This question was raised for the original Collatz function
by Shanks (1965).

\item
Murray S. Klamkin (1963), {\em Problem $63-13^{*}$},
SIAM Review {\bf 5} (1963), 275--276. \\
\newline
\hspace*{.25in}
He states the problem: ``Consider the infinite permutation 
$$
P \equiv  \left( \begin{array}{lllllll}
1 & 2 & 3 & 4 & 5 & 6 & ...\\
1 & 3 & 2 & 5 & 7 & 4 & ,,,
\end{array} \right) 
$$
taking $n \mapsto f(n)$ where 
$f: \NN^{+} \to \NN^{+}$ 
is given by  $f(3n)=2n, f(3n-1)=4n-1,
f(3n-2)=4n-3.$  
We now write $P$ as a product of cycles
$$
P \equiv (1)~ (2, 3)~ (4, 5, 7, 9, 6)~ (8, 11, 15, ...
$$
It is conjectured that the cycle  $(8, 11, 15, ...)$
is infinite. Other problems concerning $P$ are:

(a) Does the permutation $P$ consist of finitely many cycles?

(b) Are there any more finite cycles than those indicated? ``

This function was the original function proposed by
L. Collatz in his private notes in 1932. 
See Shanks (1965) and Atkin (1966) for comments on this problem.
This problem remains unsolved concerning the orbit of $n=8$
and part (a).
Concerning (b) one more cycle was found, of
period $12$ with smallest element $n=144$.
Atkin (1966) presents a heuristic argument suggesting
there are finitely many cycles.

{\em Note.} In 1963 M. Klamkin  proposed  
another problem, jointly with A. L. Tritter, concerning the orbit structure
of  a different
infinite permutation of the integers
[{\em Problem 5109},
Amer. Math. Monthly {\bf 70} (1963), 572--573].
For this integer permutation all orbits are cycles.
A solution to this problem was given
by G. Bergman, Amer. Math. Monthly {\bf 71} (1964), 569--570.\\

\item
David A. Klarner  (1981),
{\em An algorithm to determine when certain sets have $0$ density,}
Journal of Algorithms {\bf 2} (1981), 31--43. (MR 84h:10076) \\
\newline
\hspace*{.25in}
This  paper studies sets of integers that are closed under the iteration of
certain one variable affine maps. It considers the case when
all maps $f_i(x)= m x + a_i$ are expanding maps with the same ratio $m \ge 2$.
The problem is reduced to the study of sets $S(c) = \langle f_1, ..., f_k: c\rangle$ which
denotes the closure under iteration of of these maps starting from a
 single seed element $c$.
The {\em density} $\delta(S)$ of a sequence of nonnegative integers is
the lower asymptotic density
$$
\delta(S) := \liminf_{n \to \infty} \frac{1}{n+1} |S \cap \{0, 1, ..., n\}|.
$$
If $k<m$ then each  sequence $S(c)$ automatically has zero density.
If $k=m$ and the semigroup of affine maps having generators  $A:=\{ f_i: 1 \le i \le m\}$
has a nontrivial relation, then all sequences $S(c)$ have zero density,
while if this semigroup is free on $k$ generators, then they have positive density.
The paper gives an algorithm to determine whether a sequence $S(c)$ has
zero density that works for all $k$. It is based on the fact that the
number $n_t:=|A^t/ \sim|$ of distinct affine functions obtained by composition
exactly $t$ times
satisfies a linear homogeneous difference equation with constant coefficients. 
It follows that the generating function $F(z) := \sum_{t=0}^{\infty} n_t z^t$
is a rational function $F(z)= \frac{P(z)}{Q(z)}$
for relatively prime polynomials $P(z)$, $Q(z)$. The problem is reduced
to testing whether  $Q(\frac{1}{m})=0$ holds.

\item
David A. Klarner  (1982),
{\em A sufficient condition for certain semigroups to be free,}
Journal of Algebra {\bf 74} (1982), 140--148. (MR 83e:10081) \\
\newline
\hspace*{.25in}
This  paper studies sets of integers that are closed under the iteration of
certain one variable affine maps $f_i(x) = a_i x + b_i$, where $a_i, b_i$
are integers with $a_i \ge 2$. It determines sufficient conditions for the
semigroup generated by the maps $f_i$ under composition to be a free
semigroup. 
As motivation, it mentions the Erd\"{o}s problem asking whether the
set of integers generated starting from seed $S= \{ 1\}$ by the maps
$f_1(x) = 2x+1, \, f_2(x)= 3x+1$ and $f_3(x) = 6x+1$ is of positive density.
This problem was solved by D. J. Crampin and A. J. W. Hilton (unpublished)
who observed the density must be zero using the
fact that  the semigroup generated by these functions is not free. 
The paper Klarner (1981) gives an effective algorithm for the case when
all $a_i$ are equal to an integer $m \ge 2$. 
This paper first observes that a
necessary condition for freeness is that
$
\frac{1}{a_1} + \frac{1}{a_2} + \cdots + \frac{1}{a_k} \le 1.
$
It next supposes the functions are ordered so that the quantities
$p_j = \frac{b_j}{a_j -1}$  satisfy
$ p_1 \le  p_2 \le \cdots \le p_j.$
Here it notes  that strict inequalities are needed for
freeness because  an equality  $p_j=p_{j+1}$  implies 
that the generators $f_j$ and $f_{j+1}$  commute, giving a nontrivial relation.
A sufficient condition for freeness  
is  that on each collection $L(i_1, ..., i_k)$
 of compositions using exactly   $i_j$ generators of type $j$
the linear ordering of all these functions evaluated at $0$ coincides with the lexicographic
ordering  induced on functional composition,
ordering the  generators as  $f_1 < f_2 < \cdots < f_k$.
Theorem 2.1 characterizes
the latter condition, showing that   
these two orderings coincide if and only if,
 the $p_j$ satisfy the two conditions
(i) $p_1< p_2 < \cdots < p_k$
and, (ii) 
$\frac{p_k + b_j}{a_j} \le \frac{p_1+ b_{j+1}}{a_{j+1}}$ holds for $1 \le j \le k-1$. 
It then deduces some examples where the semigroups are free. 
These include $f_1(x) = 2x+ b_1, \, f_2(x) = 3x+ b_2, \, f_3(x) = 6x+ b_3)$
for  the six cases $(b_1, b_2, b_3) = (0, 3, 10), \, (0, 2, 3), \, (2, 0, 15), \, (1, 0, 2), \, (2, 1, 0), \,$
and $(1, 6, 0)$. The second of these cases gives the functions
$2x, 3x+2, 6x+3$ posed as unsolved Problem 4 in Guy (1983a).

{\em Note.}   The interest of D. J. Crampin and A. J. W. Hilton in such
questions may have arisen because they used iteration of (vector-valued) affine functions
to show the existence of $n \times n$ Latin squares orthogonal to their transpose of
all sufficiently large $n$. Roughly speaking they used constructions that took Latin squares
of one size (with extra features) and used them to build Latin rectangles of larger sizes,
whose size was an affine function of the size of the orignal squares 
They could thus produce suitable Latin squares of sizes generated by various
vector-valued affine functions, and needed to show the resulting sizes 
 included every suffciently large integer.
A large computer calculation was required, and is described in: D. J. Crampin and A. J. W. Hilton, 
{\em Remarks on Sade's disproof of the Euler conjecture with an application
to Latin squares orthogonal to their transpose,} J. Comb. Theory, Series A,
{\bf 18} (1975), 47--59. 

\item
David A. Klarner  (1988),
{\em $m$-recognizability of sets closed under certain affine functions,} 
Discrete Applied Mathematics {\bf 21} (1988), no. 3, 207--214. 
(MR 90m:68075) \\
\newline
\hspace*{.25in}
This  paper studies sets of integers $T$ that are closed under the iteration of
certain one variable affine maps. It assumes that all maps have the
special form $f_i(x) = m^{e_i} x + a_i, $ all with the same $m \ge 2$,
with all $e_i \ge 1, a_i \ge 0$. The author writes $T = \langle  A: S \rangle$, where
$A $ denotes the finite set of maps and $S$ a finite set of  ``seeds". 
The main idea of the paper is to study the base $m$-representations of
the integers and to show these are described by languages accepted
by a finite automaton; such sets are called here {\em $m$-recognizable.}
This is the content of Theorem 1. This is exhibited
on the example $f_1(x) = 3x, \, f_2(x)= 3x+1, \, f_3(x)= 3x+4.$
Theorem 2 asserts that if $S$ is an $m$-recognizable set, then so
is $T= \langle A: S \rangle$. Theorem 3 asserts that if $T$ is $m$-recognizable so is
the translated set $T+h$.

At the end of the paper more general cases are discussed. He remarks that for 
 $f_1(x) = 2x+1, \, f_2(x) = 3x+1, \, f_3(x) = 6x+1$ with $S= \{ 1\}$, 
suggested by Erd\H{o}s, 
the set $T = <A: S>$ has no discernible structure,  and that Erd\H{o}s
showed that  $T$ contains no infinite arithmetic progression.
He mentions that  Crampin and A. W. S. Hilton (unpublished) 
showed  that this sequence has
(lower asymptotic) density zero, 
answering a question of Erd\H{o}s. 
They used  the existence of  the  nontrivial relation under composition 
$$
f_1 \circ f_1 \circ f_2 = f_3 \circ f_1= 12x+7,
$$
of these affine functions, so that 
$A$ is not a free semigroup under composition in
this case. Klarner also  notes that for the functions considered in this paper there is
an effective algorithm to test whether the resulting set has density zero,
generalizing results in Klarner (1981).

{\em Note.} This paper shows  these sets are {\em $m$-automatic sequences}, in
the terminology of 
J.-P. Allouche and J. O. Shallit, {\em Automatic Sequences},  theory,
applications, generalizatins. Cambridge University Press, Cambridge 2003.

\item
David A. Klarner  and Karel Post (1992),
{\em Some fascinating integer sequences,} 
A collection of contributions in honour of Jack van Lint. 
Discrete  Mathematics {\bf 106/107} (1992),  303--309. 
(MR 93i:11031) \\
\newline
\hspace*{.25in}
This  paper studies sets of integers $T$ that are closed under the iteration of
certain multi-variable affine maps $\alpha(x_1, \cdots, x_r) = m_0 + m_1x_1 + \cdots + m_rx_r$.
It is related to Klarner (1988), which 
considered one-variable functions.
Here the authors consider a single function of two variable. 
%It studies the lattice of all sets $T$ that are closed under iteration of all functions
%in a finite family $A:= \{ \alpha_i: 1 \le i\le\}$.
%It studies a very special case, which uses a single function of two variables. 
They fix an  integer $m \ge 2$, and study sets of integers closed under action of the 
two-variable function
$$
\gamma_m(x, y) := m x + my +1.
$$
Let  $\langle \gamma_m:0 \, \rangle$ denote the set generated starting from the element $\{0\}$.
using iteration of this function, which may have a complicated structure.
They  set $G^{(1)} =\,  \langle \gamma_m: 0 \, \rangle$ and consider  the hierarchy of sumsets
$G^{(k+1)} = G^{(1)} + G^{(k)}.$ One finds that
$G^{(2k)} = \{ 0\} \cup \left( \cup_{i=1}^{2k} (m G^{(2i)} +i)\right)$.
They introduce  the affine linear functions  $\mu_i(x) = mx+i$, and prove  that
$M^{(2k)} := \langle \mu_{k+1}, \mu_{k+2}, \cdots , \mu_{2k}: \{0, 1, ..., 2k\}\, \rangle $ has
$M^{(2k)} \subseteq G^{(2k)}.$  They deduce that $M^{(2m-2)} = \NN$
and hence in Theorem 1 that $G^{(2k)} = \NN$ for all $k \ge m-1$. They deduce from this 
that $\langle  \gamma_m : 0 \rangle $ is an $m$-recognizable set in the sense of 
Klarner (1988).
They assert that for all  $m$ that the sets $G^{(1)}$ and
$G^{(2k)}$ all have positive limiting natural densities. They carry this
out  for the case $m=5$, and find that $(G^{(1)}, G^{(2)}, G^{(4)}, G^{(6)}, G^{(8)})$
have natural densities $(\frac{1}{40}, \frac{1}{8}, \frac{1}{2}, \frac{7}{8}, 1)$,
respectively. They then derive a finite automaton for generating the set $\langle \gamma_5: 0 \rangle$.

\item
David A. Klarner and Richard Rado (1973),
{\em Linear combinations of sets of consecutive integers,}
American Math. Monthly  {\bf 80} (1973), No. 9, 985--989.
(MR  48 $\#$8378) \\
\newline
\hspace*{.25in}
This paper arose from questions concerning the iteration of integer-valued
affine maps, detailed in Klarner and Rado (1974). 
This paper proves results on additive unions of sets of 
consecutive positive integers, containing certain larger sets of consecutive integers.
It proves a  vector-valued form of
such a result. The results allow improvement
of certain results  in Klarner and Rado (1974), implying for 
example that when $m,n$ are positive integers having
greatest common divisor $(m, n)=1$ the set $S= \langle 1+ mx+ny: 0 \rangle$ of integers generated
starting from $0$ by iteration of the function $\rho(x,y)= 1+ mx+ny$ contains all but
finitely many positive integers in each arithmetic progressions $(\bmod~ mn)$
that it can allowably reach.

\item
David A. Klarner and Richard Rado (1974),
{\em Arithmetic properties of certain recursively defined sets,}
Pacific J. Math. {\bf 53} (1974), No. 2, 445--463. 
(MR 50  $\#$9784) \\
\newline
\hspace*{.25in}
This  paper studies 
the smallest set of nonnegative
integers obtained from a given set $A$ under iteration of a finite set $R$ of affine
maps 
$$\rho(x_1, ..., x_r)= m_0 + m_1 x_1 + \cdots + m_r x_r,
$$
in which all  $m_i$ are integers, and $m_i \ge 0$ for $i \ge 1$.
They denote  this set $\langle R: A \rangle$. One can now ask questions
concerning  the size and structure of this
set.  This paper gives some sufficient conditions for such a set to be a finite
union of arithmetic progressions, which are called {\em per-sets}; such sets  
necessarily have positive
density. Theorem 4 shows that if $A$ is already a per-set,
then closure under a map having greatest common divisor $(m_1, ..., m_r)=1$
will give a per-set. 
Theorem 5 gives a general condition for a set $<R: A>$ to be closed
under multiplication. 
Conjecture 1 states  that if  $R$ includes
a function having  greatest common divisor $(m_1, ..., m_r)=1$, and $A= \{1\}$ 
then $<R: A>$ is a per-set. 
(The paper also announces a subsequent proof by Klarner of this conjecture; this was 
carried out  in Hoffman and Klarner (1978), (1979).) 
Conjecture 1 is proved here in some cases, showing (Theorem 11) that 
$\langle 2x+ny: 1 \rangle$ is a per-set for
all odd integers $n \ge 1$. Conjecture 2 asserts  that  for all $m,n \ge 1$ the set $<mx+ny: 1>$ 
contains a non-empty per-set. 
 In the case of one variable functions $R= \{ a_1 x + b_1, ..., a_r x+ b_r\}$
it includes a theorem of Erd\H{o}s (Theorem 8) showing that
 if $\sum \frac{1}{a_i}<1$, then the set has density $0$. The authors mention
 numerical study of the set $S= \langle  2x+1, 3x+1: 1 \rangle $, which has density $0$.
 Here Klarner (1972) had conjectured that the complement $\NN \smallsetminus S$  could
 be written as a disjoint union of infinite arithmetic progressions; this was
  proved by Coppersmith (1975).

{\em Note.} This work originally appeared as a series
of Stanford Computer Science Dept. Technical Reports in 1972,
numbered: STAN-CS-72-269.
% STAN-CS-72-274, 
Related  sequels by Klarner are  STAN-CS-72-275, STAN-CS-73-338.
These results are superseded by Hoffman and Klarner (1978), (1979).

\item
Ivan Korec (1992),
{\em The $3x  +  1$ Problem, Generalized Pascal Triangles, 
and Cellular Automata,}
Math. Slovaca, {\bf 42} (1992), 547--563.
(MR 94g:11019). \\
\newline
\hspace*{.25in}
This paper shows that the iterates of the Collatz function $C(x)$ can
actually be encoded in a simple way by a one-dimensional nearest-neighbor
cellular automaton with 7 states.
The automaton encodes the iterates of the function $C(x)$
written in base 6.
The encoding is possible
because the map $x  \rightarrow  3x +  1$ in base 6 does not have 
carries propagate.
(Compare the $C(x)$-iterates of $x =  26$ and $x =  27$ in base 6.)
The $3x  +  1$ Conjecture is reformulated in terms of special 
structural properties
of the languages output by such cellular automata.

\item
Ivan Korec (1994), 
{\em A Density Estimate for the $3x+1$ Problem, }
Math. Slovaca {\bf 44} (1994), 85--89.
(MR 95h:11022). \\
\newline
\hspace*{.25in}
Let $S_\beta  = \{ n :{\rm some} ~ T^{(k)} (n) < n^\beta \}$.
This paper shows that for any
$\beta  > \frac {\log ~3}{\log ~4} \doteq .7925$
the set $S_\beta$ has density one.
The proof follows the method of E.\ Heppner [Arch. Math.
{\bf 31}
(1978) 317-320].

\item
Ivan  Korec and Stefan Znam (1987), 
{\em A Note on the $3x  +  1$ Problem,}
Amer. Math. Monthly {\bf 94} (1987), 771--772.
(MR 90g:11023). \\
\newline
\hspace*{.25in}
This paper shows
 that to prove the $3x  +  1$ Conjecture it suffices to verify it for
the set of all numbers $m  \equiv  a$
$(\bmod ~p^n)$, for any
fixed $n  \geq  1$, provided that 2 is a
primitive root $ (\bmod ~p$)
and $(a, p) =  1$.
This set has density $p^{-n}$.

\item
Ilia Krasikov (1989), 
{\em How many numbers satisfy the $3x  +  1$ Conjecture?,}
Internatl. J. Math. \& Math. Sci. {\bf 12} (1989), 791--796.
(MR 91c:11013). \\
\newline
\hspace*{.25in}
This paper
shows that the number of integers $ \leq  x$ for which the $3x  +  1$
function has an iterate that is 1 is at least $x^{3/7}$.
More generally, if 
$\theta_a (x) =  \{ n \leq  x : ~  {\rm some}~ T^{(k)} (n) =  a \}$, 
then he shows that, for $a \not\equiv 0~~(\bmod~3)$, 
$\theta_a (x)$ contains at least
$x^{3/7}$ elements, for large enough $x$.
For each fixed $k \ge 2$
this  paper derives a system  of difference inequalities 
based on information $(\bmod~ 3^k)$. The bound $x^{3/7}$
was obtained using $k=2$, and by using larger
values of $k$ better exponents can be obtained.
This was done in  Applegate and Lagarias (1995b)
and Krasikov and Lagarias (2003).

\item
James R. Kuttler (1994),
{\em On the $3x+1$ Problem,} 
Adv. Appl. Math. {\bf 15} (1994), 183--185.
(Zbl. 803.11018.) \\
\newline
\hspace*{.25in}
The author states the oft-discovered fact that
$T^{(k)} (2^k n-1 ) = 3^k n-1$ and derives his main result
that if $r$ runs over all odd integers $1 \leq r  \leq 2^k - 1$
then $T^{(k)} (2^k n+r) = 3^p n+s$,
in which $p \in \{ 1, 2 , \ldots , k \}$ and
$1 \leq s \leq 3^p$,
and each value of $p$ occurs exactly $\left(  {k-1 \atop p-1} \right) $
times.
Thus the density of integers with $T^{(k)} (n)  >  n$ is exactly
$\frac{\alpha}{2^k}$, where $\alpha =  \sum_{p > \theta k}   {k-1 \choose p-1}  $,
and $\theta =  \log_2 ~3$.
This counts the number of inflating vectors in Theorem~C in Lagarias (1985).

\item
Jeffrey C. Lagarias (1985), 
{\em The $3x  +  1$ problem and its generalizations,}
Amer. Math. Monthly {\bf 92} (1985), 3--23.
(MR 86i:11043). \\
\newline
\hspace*{.25in}
This paper is a
 survey with extensive bibliography of known results on $3x  +  1$
problem and related problems up to 1984.
It also contains improvements of previous results and some new results, 
including in particular Theorems D, E, F, L and M.

\hspace*{.25in}
Theorem O has several misprints.
The method of Conway (1972) gives, for any partial recursive function 
$f$, a periodic piecewise linear function $g$, with the property:
iterating $g$ with the starting value $2^n$ will never be a power of 2 if
$f(n)$ is undefined, and will eventually reach a power of 2 if 
$f(n)$ is defined, and the first such
power of 2 will be $2^{f(n)}$.
Thus, in parts (ii) and (iii) of Theorem O, occurrences of $n$ 
must be replaced by $2^n$. On page 15, the thirteenth partial quotient of
$\log_2 (3)$ should be
$q_{13}= 190537$.

An updated version of this paper is accessible on the World-Wide Web: \\
{\tt http://www.cecm.sfu.ca/organics/papers/lagarias/index.html}.
 %or
%{\tt http://www.mathsoft.com} (unsolved problems menu)

\item
Jeffrey C. Lagarias (1990), 
{\em The set of rational cycles for the $3x  +  1$ problem,} 
Acta Arithmetica {\bf 56} (1990), 33--53.
(MR 91i:11024). \\
\newline
\hspace*{.25in}
This paper studies the sets of those integer cycles of
$$
T_k (x) =  \left\{
\begin{array}{cl}
\df{3x+k}{2} & \mbox{if}~~x \equiv 1~~(\bmod~2 ) , \\
~~~ \\
\df{x}{2} & \mbox{if}~~x \equiv 0~~(\bmod~2 ) ~,
\end{array}
\right.
$$
for positive $k  \equiv  \pm 1~~(\bmod~6)$ which have $(x, k) =  1$.
These correspond to rational cycles $\frac{x}{k}$ of the
$3x  +  1$ function $T$.
It conjectures that every $T_k$ has such an integer cycle.
It shows that infinitely many $k$ have at least $k^{1 -  \epsilon}$
distinct such cycles of period at most $\log~ k$, and infinitely many $k$
have no such cycles having period length less than
$k^{1/3}$.
Estimates are given for the counting function $C(k, y)$
counting the number of such cycles of $T_k$
of period $ \leq  y$, for all
$k  \leq  x$ with $y =  \beta  \log x$.
In particular $C(k,  1.01k)  \leq  5k ( \log~  k)^5$.

\item
Jeffrey C. Lagarias (1999),
{\em How random are 3x + 1 function iterates?,}
in: {\em The Mathemagician and Pied Puzzler: 
A Collection in Tribute to Martin Gardner},
A. K. Peters, Ltd.: Natick, Mass. 1999, pp. 253--266. \\
\newline
\hspace*{.25in}
This paper briefly summarizes results on extreme trajectories of 3x+1 
iterates, including those in  Applegate and Lagarias (1995a), (1995b),(1995c)
and Lagarias and Weiss (1992). It also presents some large integers $n$ 
found  by V. Vyssotsky whose trajectories take  $c~ log~ n$ steps to
iterate to $1$ for various
$c > 35$. For example $n= 37~66497~18609~59140~59576~52867~40059$ has 
$\sigma_\infty(n)=2565$ and $\gamma (n)=35.2789$. 
It also mentions some unsolved problems.

The examples above were  
the largest values of $\gamma$ for $3x+1$ iterates known
at the time, but are now superseded by
examples of Roosendaal (2004+) achieving
$\gamma= 36.716.$

\item
Jeffrey C. Lagarias, Horacio A. Porta and Kenneth B. Stolarsky (1993),
{\em Asymmetric Tent Map Expansions I.\ Eventually Periodic Points, }
J. London Math. Soc., {\bf 47} (1993), 542--556. 
(MR 94h:58139). \\
\newline
\hspace*{.25in}
This paper studies the set of eventually
periodic points ${\rm Per} ( T_\alpha )$ 
of the asymmetric tent map
$$
T_\alpha (x) =  \left\{
\begin{array}{ll}
\alpha x & \mbox{if $0 < x \leq \df{1}{\alpha}$} \\
~~~ \\
\df{\alpha}{\alpha-1} (1-x) & \mbox{if $\df{1}{\alpha} \leq x < 1$} ~,
\end{array}
\right.
$$
where $\alpha  >  1$ is real.
It shows that ${\rm Per} ( T_\alpha ) =  \QQ ( \alpha )  \cap  [0, 1]$ 
for those
$\alpha$ such that both $\alpha$ and $\frac{\alpha}{\alpha -1}$ are 
Pisot numbers. It finds $11$ such numbers, of degree up to four,
 and proves that  the set of all such numbers is finite.

It conjectures that the property
${\rm Per} ( T_\alpha ) =  \QQ ( \alpha )  \cap  [0, 1]$  
holds for certain other $\alpha$,
including the real root of $x^5  -  x^3  -  1 = 0$.
The problem of proving that 
${\rm Per} ( T_\alpha ) = \QQ ( \alpha )  \cap [ 0, 1]$ 
in these cases appears analogous to the problem of proving that the
$3x  +  1$ function has no divergent trajectories.

[C. J. Smyth, {\em There are only
eleven special Pisot numbers}, (Bull. London Math. Soc. {\bf 31} (1989) 1--5)
proved that the set of $11$ numbers found above is the complete set.]

\item
Jeffey C. Lagarias, Horacio A. Porta and Kenneth B. Stolarsky (1994),
{\em Asymmetric Tent Map Expansions II.\ Purely Periodic Points, }
Illinois J. Math. {\bf 38} (1994), 574--588. 
(MR 96b:58093, Zbl 809:11042). \\
\newline
\hspace*{.25in}
This paper continues Lagarias, Porta and Stolarsky (1993).
It studies the set 
${\rm Fix}(T_{\alpha})$ of purely periodic
points of the asymmetric tent map $T_{\alpha}(\cdot)$
and the set ${\rm Per}_0(T_{\alpha})$ with 
terminating $T_{\alpha}$-expansion, 
in those cases wwhen $\alpha$ and
$\frac{\alpha}{1 - \alpha}$ are simultaneously Pisot
numbers. It shows that 
${\rm Fix}(T_{\alpha}) \subseteq \{ \gamma \in \QQ(\alpha)  \mbox{~and~}
\sigma(\alpha) \in A_{\alpha}^\sigma \mbox{~for~all~embeddings}~
\sigma: \QQ(\alpha) \to \CC \mbox{~with~} \sigma(\alpha) \ne \alpha \},$
in which each set $ A_{\alpha}^\sigma$ is a compact set in $\CC$
that is the attractor of a certain hyperbolic iterated function
system. It shows that equality holds in this
inclusion in some cases, and not in others.
Some related results for ${\rm Per}_0(T_{\alpha})$ are established.

\item
Jeffrey C. Lagarias and Alan Weiss (1992),
{\em The $3x  +  1$ Problem: Two Stochastic Models, }
Annals of Applied Probability {\bf 2} (1992), 229--261.
(MR 92k:60159). \\
\newline
\hspace*{.25in}
This paper studies two stochastic models that mimic the ``pseudorandom''
behavior of the $3x  +  1$ function. The models are branching
random walks, and the analysis uses the theory of large
deviations.

For the models  the
average number of steps to get to $1$ is $\alpha_0 \log n$,
where
$\alpha_0 = \left(\frac{1}{2} \log \frac{3}{4}\right)^{-1}\approx 6.9$
For both models it is shown that there is a constant
$c_0 =  41.677. . .$ such that with probability one
for any $\epsilon  >  0$ only finitely many $m$
take $(c_0  +  \epsilon )  \log m$ iterations to take the value 1,
while infinitely many $m$ take at least 
$(c_0  -  \epsilon )  \log  m$ iterations
to do so.
This prediction is shown to be consistent with empirical data 
for the $3x  +  1$ function.

The paper also studies the maximum excursion
$$
t(n) : = \max_{k \geq 1} T^{(k)} (n)
$$
and conjectures that $t(n) < n^{2+o(1)}$ as $n \to \infty$.
An analog of this conjecture is proved for one
stochastic model.
The conjecture $t(n) \leq n^{2 + o(1)}$ is consistent with empirical
data of Leavens and Vermeulen (1992) for
$n < 10^{12}$.

\item
Gary T. Leavens and Mike Vermeulen (1992),
{\em $3x +1$ Search Programs, }
Computers \& Mathematics, with Applications 
{\bf 24}, No. 11,(1992), 79--99. 
(MR 93k:68047). \\
\newline
\hspace*{.25in}
This paper describes  methods for computing $3x+1$ function iterates,
and gives results of extensive computations done on a
distributed network of workstations, taking an estimated 10\ CPU-years in
total.
The $3x+1$ Conjecture is verified for all $n < 5.6 \times 10^{13}$.
Presents statistics on various types of extremal trajectories of
$3x+1$ iterates in this range.
Gives a detailed discussion of techniques used for program optimization.

[This  bound for verifying the
$3x+1$ conjecure is now superseded by that of Oliveira e Silva (1999).]

\item
George M. Leigh (1986), 
{ \em A Markov process underlying the generalized Syracuse algorithm,} 
Acta Arithmetica {\bf 46} (1986), 125--143.
(MR 87i:11099) \\
\newline
\hspace*{.25in}
This paper considers mappings $T(x) = \frac  {m_i  x  -  r_i }{d}$ if
$x  \equiv  i~ (\bmod  ~d)$, where all
$r_i  \equiv  im_i~ (\bmod  ~d$). This work is motivated by earlier
work of Matthews and Watts (1984), (1985), and introduces
significant new ideas. 

Given such a mapping $T$
and an auxiliary modulus $m$ it introduces two Markov chains, denoted $\{X_n\}$ and $\{Y_n\}$,
%(depending on $(T, m)$) 
whose behavior encodes information on the iterates of $T$ (modulo $md^k$) for all
$k \ge 1$. In general, both Markov chains have a countable number of states; however
in many interesting cases both these chains are   finite state Markov chains. The states
of each chain are labelled by certain congruence classes 
$B(x_j, M_j):= \{ x: x \equiv x_j (\bmod~M_j)\}$ with $M_j | m d^k$ for some $k \ge 1$.
(In particular $B(x_j, M_j) \subset B(x_k, M_k)$ may occur.) Any finite path in such a  chain can be realized by a sequence of iterates
$x_j= T^{j}(x_0)$ for some  $x_0 \in \ZZ$ satisfying the congruences specified by the states. 
The author shows the two chains contain equivalent information; the chain $\{ X_n\}$ has
in general fewer states than the chain $\{Y_n\}$ but the latter is more suitable for
theoretical analysis, in particular one chain is finite if and only if the other one is.  
If such a  chain has a positive ergodic class of states, then the limiting frequency distribution of
state occupation of a path in the chain in this class exists, and from this
the corresponding frequencies of iterates  
 in a given congruence class (mod $m$) can be calculated (Theorem 1-  Theorem 4).  
 
 For applications to the map $T$ on the integers, the
  author 's guiding conjecture is: {\em any condition that occurs with probability zero in the
 Markov chain model does not occur in divergent trajectories on $T$.}
Thus the author conjectures  that in cases where the Markov chain is finite,  there are 
limiting densities that
divergent trajectories for $T$ (should they exist) spend in each residue class 
$j$ $(\bmod  ~m)$, obtained from frequency distribution in ergodic states
of the chain. 

The author gets a complete analysis whenever the associated Markov
chains are finite.  Furthermore Theorem 7
shows that any map $T$ satisfying the condition 
$gcd(m_i, d^2) = gcd(m_i, d)$ for $0 \le i \le d-1$,
has, for every modulus $m$,
 both Markov chains $\{X_n\}$ and $\{Y_n\}$ 
being  finite chains.  All  maps studied in  the earlier work of Matthews and Watts (1984, 1985)
satisfy this condition, and the author  recovers nearly all of their results in 
these cases.

In Sect. 5 the author suggests that  when infinite Markov chains occur, that one 
approximate them using  
a series of larger and larger finite Markov chains obtained by truncation.
Conjecture 2  predicts limiting frequencies of iterates (mod $m$) for the map $T$,
even in these cases.
The paper concludes with several worked examples. 

This Markov chain approach was extended further by  Venturini (1992).

\item
Simon Letherman, Dierk Schleicher and Reg Wood (1999),
{\em On the $3X+1$ problem and holomorphic dynamics,}
Experimental Math. {\bf 8}, No. 3 (1999), 241--251. (MR 2000g:37049.) \\
\newline
\hspace*{.25in}
This paper studies the class of entire functions
$$
f_h(z) := \frac{z}{2} + (z + \frac{1}{2})(\frac{1 + \cos \pi z}{2})
 + \frac{1}{\pi}(\frac{1}{2} - \cos \pi z) \sin \pi z
 + h(z) (\sin \pi z)^2,
$$
in which $h(z)$ is an arbitrary entire function. Each function in this
class reproduces the $3x+1$-function on the integers, and the set of
integers is contained in the set of critical points of the function.
The simplest such function takes $h(z)$ to be identically zero, and is
denoted $f_0(z)$. The authors study the iteration properties of this
map in the complex plane $\CC.$ They show that $\ZZ$ is contained in
the Fatou set of $f_h(z)$. There is a classification of the connected
components of
the Fatou set of an entire function into six categories: (1) (periodic)
immediate basins of (super-) attracting periodic points, (2) (periodic)
immediate basins of attraction of rationally indifferent periodic points,
(3) (periodic) Siegel disks, (4) periodic domains at infinity (Baker domains), (5) preperiodic 
components of any of the above, and (6) wandering domains.
The following results all apply to the function $f_0(z)$, and some of them
are proved for more general $f_h(z)$.
The Fatou component of any integer must be in the basin of attraction
of a superattracting periodic point or be in a wandering domain; the authors
conjecture the latter does not happen.
The existence
of a divergent trajectory is shown equivalent to $f_0(z)$ containing
a wandering domain containing some
integer. No two integers are in the same
Fatou component, except possibly $-1$ and $-2$. The
real axis contains points of the Julia set between any two integers,
except possibly $-1$ and $-2$.
It is known that  holomorphic dynamics of an entire function $f$ 
is controlled by its critical points. The critical points of $f_0(z)$ on
the real axis consist of the integers together with $-\frac{1}{2}$.
The authors would like to choose $h$ to reduce the number of other
critical points, to simplify the dynamics. However they show that
any map $f_h(z)$ must contain at least one more critical point in addition to
the critical points on the real axis (and presumably infinitely many.)
The authors compare and contrast their results with the real dynamics studied
in Chamberland (1996).

\item
Heinrich Lunkenheimer (1988),
{\em Eine kleine Untersuchung zu einem zahlentheoretischen Problem},
PM (Praxis der Mathematik in der Schule) {\bf 30} (1988), 4--9.\\
\newline
\hspace*{.25in}
The paper considers the Collatz problem, giving no prior history
or references.
It observes the coalescences of some orbits.
It lists some geometric series of integers which iterate to $1$.
\item
Kurt Mahler (1968),
{\em An unsolved problem on the powers of $3/2$},
J. Australian Math. Soc. {\bf 8} (1968), 313--321.
(MR 37 \#2694). \\
\newline
\hspace*{.25in}
A $Z$-number is a real number $\xi$ such that
$0 \le {{ (\frac{3}{2})^k \xi }} \le \frac{1}{2}$
holds for all $k \ge 1$, where ${{x}}$ denotes the fractional
part of $x$. Do $Z$-numbers exist?
The $Z$-number problem  was originally
proposed by Prof. Saburo Uchiyama (Tsukuba Univ.)
according to S. Ando (personal communication),
and was motivated by a connection with the
function $g(k)$ in Waring's problem, for which
see G. H. Hardy and E. M. Wright , {An Introduction to the
Theory of Numbers} (4-th edition), Oxford Univ. Press 1960,
Theorem 393 ff. and Stemmler (1964).

Mahler shows that
existence of $Z$-numbers relates to a question concerning
the iteration of the function
$$
g(x) = \left\{
\begin{array}{cl}
\df{3x+1}{2} & \mbox{if} ~~x \equiv 1~~ (\bmod~2) ~. \\
~~~ \\
\df{3x}{2} & \mbox{if}  ~~x \equiv 0~~ (\bmod~2)~.
\end{array}
\right.
$$
Mahler showed that a  $Z$-number exists in the
interval $[n, n+1)$ if and only if no iterate
$g^{(k)}(n) \equiv 3 ~(\bmod~4).$ He uses 
this relation to prove  that
the number of $Z$-numbers below $x$ is at most $O(x^{0.7}).$
He conjecures that no  $Z$-numbers exist, a problem which is
still unsolved. 

{\em Note.} Leopold Flatto (1992)   subsequently improved
Mahler's upper bound on $Z$-numbers below $x$ to $O(x^{\theta})$,
with $\theta = \log_2 \frac{3}{2} \approx 0.59$.
%[L. Flatto, {\em $Z$-numbers and $\beta$ transformations}, pp. 181--201 in:
%{\em Symbolic Dynamics and Applications (New Haven, CT 1991)},
%Contemp. Math. Vol. 135, AMS: Providence RI 1992)]

\item
Jerzy Marcinkowski (1999),
{\em Achilles, Turtle, and Undecidable Boundedness
Problems for Small DATALOG programs},
SIAM J. Comput. {\bf 29} (1999), 231--257. 
(MR 2002d:68035). \\
\newline
\hspace*{.25in}
DATALOG is the language of logic programs without
function symbols. A DATALOG program consists of a finite
set of Horn clauses in the language of first order
logic without equality and without functions.

The author introduces {\em Achilles-Turtle
machines}, which model the iteration of 
{\em Conway functions},  as introduced in Conway (1972).
These are functions $g: \NN \to \NN$ having the form
$$
g(n) = \frac{a_j}{q_j} n 
~~\mbox{~if~} n \equiv j (\bmod~p),~ 0 \le j \le p - 1,
$$
where for each $j$,  $a_j \ge 0,  q_j\ge 1$ are integers 
with  $q_j | gcd (j, p)$ and with $\frac{a_j}{q_j} \le p$.
He cites Devienne, Leb\'{e}gue and Routier (1993)
for the idea of relating Conway functions 
to Horn clauses. In Section 2.3 
he explicitly describes  an
Achilles-Turtle machine associated to computing the Collatz 
function. 

The paper proves that several questions concerning
the uniform boundedness of computations are undecidable.
These include uniform boundedness for ternary linear programs;
unfiorm boundedness for single recursive rule ternary
programs; and uniform boundedness of single rule programs.
These results give no information regarding the $3x+1$ function
itself.

\item
D\u{a}nu\c{t} Marcu (1991), 
{\em The powers of two and three,}
Discrete Math. {\bf 89} (1991), 211-212. \\
(MR 92h:11026). \\
\newline
\hspace*{.25in}
This paper obtains a similar result to Narkiewicz (1980),
with a similar proof, but obtains a slightly worse bound on the exceptional set 
$N(T)< 2.52T^{\theta}$ where $\theta= \frac{\log 2}{\log 3}.$\\

{\em Note.} This paper  cites  a paper of Gupta appearing in the same journal
as Narkiewicz's paper [Univ. Beograd Elektrotech. Fak.  Ser. Mat. Phys.] slightly
earlier in the same year,  but fails to cite Narkiewicz (1980). 
 For clarification, consult the Wikipedia article on D. Marcu.

\item
Maurice Margenstern and Yuri Matiyasevich (1999),
{\em A binomial representation of the $3X +1$ problem,}
Acta Arithmetica {\bf 91}, No. 4 (1999), 367--378.
(MR 2001g:11015). \\
\newline
\hspace*{.25in}
This paper encodes the $3x+1$ problem as a logical
problem using one universal quantifier and
existential quantifiers, with 
an arithmetical formula using polynomials and binomial
coefficients. The authors observe that the use of such expressions
in a language with binomial coefficients
often leads to shorter formulations than are possible
in a language just allowing polynomial equations and quantifiers.
They give three equivalent restatements of the $3x+1$ Conjecture
in terms of quantified binomial coefficient equations.

\item
Keith R. Matthews (1992), 
{\em Some Borel measures associated with the generalized Collatz mapping,} 
Colloq. Math. {\bf 53} (1992), 191--202.
(MR 93i:11090). \\
\newline
\hspace*{.25in}
This paper studies maps of the form
$$
T(x) = \frac{m_i  x-r_i }{d} ~~~ {\rm if} ~~ x  \equiv  i~(   \bmod ~d)~.
$$
These extend first to maps on the $d$-adic integers $\ZZ_d$,
and then further to maps on the polyadic integers $\hat{\ZZ}$.
Here $\hat{\ZZ}$ is the projective limit of the set of
homomorphisms $\phi_{n,m} : \ZZ  \rightarrow  \ZZ /m \ZZ$ where $m | n $.
The open sets $B(j,m) $ $: = \{ x \in \hat{\ZZ} : x \equiv j (  \bmod~m)\} $
put a topology on $\hat{\ZZ}$, which has a Haar measure 
$\sigma ( B( j,m)) = \frac{ 1}{m}$.
The paper proves conjectures of Buttsworth and Matthews (1990) 
on the structure
of all ergodic open sets $( \bmod~m)$.
In particular the ergodic sets link together to give finitely many
projective systems, each giving $T$-invariant measure on $\hat{\ZZ}$.
The paper gives examples where there are infinitely many ergodic sets.

\item
Keith R. Matthews and George M. Leigh (1987),
{\em A generalization of the Syracuse algorithm to $F_q [x]$,}
J. Number Theory {\bf 25} (1987), 274--278. (MR 88f:11116). \\
\newline
\hspace*{.25in}
This paper defines mappings analogous to the $3x  +  1$ function on 
polynomials over finite fields,
e.g. $T(f) =  \frac{f}{x}$ if $f  \equiv  0~ (\bmod  ~x)$ and
$\frac{(x  +  1)^3 f  +  1}{x}$ if $f  \equiv  1~~(\bmod ~x)$, over GF(2).
It proves that divergent trajectories exist for certain such maps.
These divergent trajectories
have a regular behavior.

\item
Keith R. Matthews and Anthony M. Watts (1984),
{\em A generalization of Hasse's generalization of the Syracuse algorithm,} 
Acta. Arithmetica {\bf 43} (1984), 167--175.
(MR 85i:11068). \\
\newline
\hspace*{.25in}
This paper studies functions
$T(x) = \frac{m_i  x  -  r_i}d$ for
$x  \equiv  i~ (\bmod  ~d)$, where all $m_i $ are positive integers and
$r_i  \equiv  im_i~ (\bmod  ~d)$.
It is shown that if $\{  T^{(k)} (m) : k \geq  0 \}$
is unbounded and uniformly distributed $(\bmod ~d)$
then $m_1 m_2 \cdots m_d  >  d^d$ and
$\lim_{k \rightarrow \infty} | T^{(k)} (m) |^{1/k} =  
\frac{1}{d} (m_1 \cdots m_d )^{1/d}$.
The function $T$ is
extended to a mapping on the $d$-adic integers and is shown 
to be strongly mixing,
hence ergodic, on $\ZZ_d$.
The trajectories $\{  T^{(k)} ( \omega : k  \geq  0 \}$ for
almost all $\omega  \in \ZZ_d$ are equidistributed
$(\bmod ~d^k)$ for all $k  \geq  1$.

\item
Keith R.  Matthews and Anthony M. Watts (1985), 
{\em A Markov approach to the generalized Syracuse algorithm, }
Acta Arithmetica {\bf 45} (1985), 29--42. 
(MR 87c:11071). \\
\newline
\hspace*{.25in}
This paper studies the functions 
$T(x) = \frac{m_i  x  -  r_i}d$ for
$x  \equiv  i~ (\bmod ~ d)$, where all $m_i $ are positive integers and
$r_i  \equiv  im_i~ (\bmod  ~d)$, which were considered in
Matthews and Watts\ (1984).
Given a modulus $m$ one associates to $T$ a row-stochastic matrix
$Q =  [q_{jk} ]$ in which $j, k$ index residue classes $(\bmod  ~m$) and
$q_{jk}$ equals $\frac{1}{md}$ times the number of residue classes
$(\bmod ~md$) which are
$ \equiv k ( \bmod~m)$ and whose image under $T$ is $ \equiv j  (  \bmod~m)$.
It gives sufficient conditions for the entries
of $Q^l$ to be the analogous probabilities associated to the 
iterated mapping $T^{(l)}$.
Matthews and Watts conjecture that if $\sS$ is an ergodic set 
of residues $(\bmod  ~m$) and
$( \alpha_i ;~ i  \in  \sS )$ is the corresponding stationary vector on $\sS$,
and if 
$A =  \Pi_{i  \in \sS} \left( \frac{m_i }{d} \right)^{\alpha_i }  <  1$ 
then all
trajectories of $T$ starting in $\sS$ will eventually be periodic, 
while if $A  >  1$ almost all trajectories starting in $\sS$ will
diverge.
Some numerical examples are given.

\item
Danilo Merlini and Nicoletta Sala (1999),
{\em On the Fibonacci's Attractor and the Long Orbits in the 
$3n+1$ Problem}, International Journal of Chaos Theory and
Applications {\bf 4}, No. 2-3 (1999), 75--84.\\
\newline
\hspace*{.25in}
This paper studies heuristic models for the Collatz tree of inverse iterates
of the Collatz function $C(x)$, which the authors call the "chalice".
They also consider the length of  the longest trajectories for the Collatz map. 
 The author's predictions are that
the number of leaves in the Collatz tree at depth $k$ should grow like
$Ac^k$ where $c= \frac{1}{2}(1+ \sqrt{\frac{7}{3}}) \approx 1.2637$
and $A=\frac{3c}{3c-2} \approx 0.6545$.
This estimate is shown to agree  closely with numerical data computed to depth $k=32$.
They predict that the longest orbits of the 
Collatz map should take no more than $67.1 \log n$ steps to reach $1$.
This is shown to agree with empirical data for the longest orbit found
up to ${10}^{10}$ by Leavens and Vermeulen (1992).

{\em Note.} These models have
some features related to that in Lagarias and Weiss (1992)
for the $3x+1$ map $T(n)$. The asymptotic constants they obtain differ from those 
In the  models of Lagarias and Weiss (1992) because the Collatz function $C(x)$
takes one extra iteration each time an odd number occurs, compared to 
for the  $3x+1$ mapping $T(x)$. The Applegate and Lagarias 
estimate  for the number of nodes
in the $k$-th level of a $3x+1$ tree is $(\frac{4}{3})^{k(1+o(1)}$,
and  their estimate for the maximal number of steps to reach $1$ for the
$3x+1$ map is about $41.677 \log n$.

\item
Karl Heinz Metzger (1995),
{\em Untersuchungen zum $(3n+1)$-Algorithmus.
Teil I: Das Problem der Zyklen,}
PM (Praxis der Mathematik in der Schule) {\bf 38} (1995), 255--257. \\
({\em Nachtrag zum Beweis de $(3n+1)$-Problems, }
PM  {\bf 39} (1996), 217.)\\
\newline
\hspace*{.25in}
The author studies cycles for the $3x+1$ map, and more generally
considers the $bx+1$ map for odd $b$.
 He exhibits a cycle for the $5x+1$ that does not contain $1$, namely
$n=13$.  He obtains a general formula (4.1) for a rational number $a$ to be in a  cycle
of the $3x+1$ map, as in Lagarias (1990).
At the end  of the paper  is a theorem asserting  that the only cycle of the Collatz map acting on
the positive integers is the trivial cycle $\{ 1, 2, 4\}$.
 The proof of this theorem has a gap, however. The author 
points this out in the addendum, and says he hopes to return to the question in
future papers. 

The author's approach  to proving that the trivial cycle is the
only cycle on the positive integers is as follows. He observes that the condition (4.1) for $c$
to be a rational cycle with $\nu$ odd elements can be rewritten as 
$$
c= \frac{A_{\nu} + x_0^{\ast}}{A_{\nu} + y_0^{\ast}}
$$
for certain integers $(x_0^{\ast}, y_0^{\ast})$, 
taking $A_{\nu}=2^{2\nu} - 3^{\nu}$.
%3^{\nu -1} + 3^{\nu-2} 2^2 + \cdots + 3^0 2^{2\nu -2} = 2^{2\nu} - 3^{\nu}$.
He then observes that this can be viewed as  a  
special case of a linear Diophantine equation in $(x,y)$, 
$$
c= \frac{A_{\nu} + x}{A_{\nu} + y}, ~~~~\mbox{with}~~A, c~\mbox{fixed}, 
$$
which has the  general solution $(x, y) = ( (c^2-1)A + \lambda c, (c-1)A + \lambda c)$, for
some integer $\lambda$. His proof   shows that the case where $(x_0^{\ast}, y_0^{\ast})$
has associated value  $\lambda=0$ leads to a contradiction. The gap in the proof is
that cases where $\lambda \ne 0$ are not ruled out.

\item
Karl Heinz Metzger (1999),
{\em Zyklenbestimmung beim $(bn+1)$-Algorithmus},
PM (Praxis der Mathematik in der Schule) {\bf 41} (1999), 25--27.\\
\newline
\hspace*{.25in}
For an odd number $b$ the Collatz version of the $bx+1$ function is 
$$
C_{b}(x) = 
\left\{
\begin{array}{cl}
bx+1 & \mbox{if}~ x \equiv 1~~ (\bmod ~2 ) ~, \\
~~~ \\
\df{x}{2} & \mbox{if} ~~x \equiv 0~~ (\bmod~2) ~.
\end{array}
\right.
$$
The author studies  cycles for  the $bn+1$ map.
Let the {\em size} of a cycle count the number of distinct odd integers it contains.
The author first observes that $1$ is in a size  one
cycle for the $bn+1$ problem if and only if $b=2^\nu -1$ 
for some $\nu \ge 1$.
For the $5x+1$ problem he  shows that on the positive integers
there is no cycle of size $1$, a unique cycle of size $2$,
having  smallest element $n=1$, and exactly two cycles of size  $3$,
having smallest elements $n=13$ and $n=17$, respectively. 

At the end of the paper the author 
states a theorem asserting  that the only cycle of the Collatz map acting on
the positive integers is the trivial cycle $\{ 1, 2, 4\}$.  This paper apparently
is intended to repair the faulty proof in Metzger (1995). However the proof of this
result still has a gap.

\item
Pascal Michel (1993), 
{\em Busy Beaver Competition and Collatz-like Problems,}
Archive Math. Logic {\bf 32} (1993), 351--367.
(MR 94f:03048). \\
\newline
\hspace*{.25in}
The Busy Beaver problem is to find that Turing machine
which, among all $k$-state Turing machines, when given the 
empty tape as input, eventually halts and produces the
largest number $\Sigma (k)$ of ones on the output tape.
The function $\Sigma (k)$ is well-known to be non-recursive.
This paper shows that the current Busy Beaver record-holder for 
5-state Turing
machine computes a Collatz-like function.
This machine $M_5$ of Marxen and Buntrock [Bulletin EATCS No.
{\bf 40}
(1990) 247--251] has $\Sigma (M_5 ) = 47,176,870$.
Michel shows that $M_5$ halts on all inputs if and only if iterating the
function $g( 3m) =  5m  +  6$,
$g(3m + 1) =  5m + 9$, $g(3m+2) =  \uparrow$ eventually halts at 
$\uparrow$ for all inputs.
Similar results are proved for several other 5-state Turing machines.
The result of Mahler [J. Aust. Math. Soc.
{\bf 8}
(1968) 313-321] on $Z$-numbers is restated in this framework:
If a $Z$-number exists then there is an input $m_0$
such that the iteration $g(2m) =  3m$,
$g(4m~+~ 1) =  6m~+~ 2$,
$g(4m~+~ 3) =  \uparrow$
never halts with $\uparrow$ when started from $m_0$.

\item
Filippo Mignosi (1995),
{\em On a Generalization of the $3x  +  1$ Problem,} 
J Number Theory, {\bf 55} (1995), 28--45.
(MR 96m:11016). \\
\newline
\hspace*{.25in}
For real $\beta  >  1$, define the function 
$T_\beta : \NN  \rightarrow  \NN$ by
$$
T_\beta (n) =  \left\{
\begin{array}{cl}
\lceil \beta n \rceil & \mbox{if}~n \equiv 1~~(\bmod~2 ) \\
~~~ \\
\df{n}{2} & \mbox{if $n \equiv 0~~(\bmod~2 ) ~.$}
\end{array}
\right.
$$
Then $\beta =  \frac{3}{2}$ gives the
$(3x  +  1)$-function.
Conjecture $C_\beta$ asserts that $T_\beta$ has finitely many cycles and
every $n  \in  \NN$ eventually enters a cycle under iteration of $T_\beta$.
The author shows that, for any fixed $0 < \epsilon < 1$, and for $\beta$
either transcendental or rational with an even denominator, if 
$1 < \beta < 2$,
then the set 
$S( \epsilon , \beta) = \{ n:$ some $T_\beta^{(k)} (n) < \epsilon n \}$
has natural density one, while if $\beta  >  2$
then $S^\ast ( \epsilon , \beta ) = \{  n :$ some 
$T_\beta^{(k)} (n)  >  \epsilon^{-1} n \}$
has natural density one.
For certain algebraic $\beta$ different behavior may occur,
and Conjecture $C_\beta$ can sometimes be settled.
In particular Conjecture $C_\beta$ is true for $\beta =  \sqrt 2$
and false for $\beta = \frac{1  +  \sqrt 5}{2}$.

\item
Herbert M\"{o}ller (1978),
{\em Uber Hasses Verallgemeinerung des Syracuse-Algorithmus (Kakutanis
Problem)}, Acta Arith. {\bf 34} (1978), No. 3, 219--226.
(MR 57 \#16246). \\
\newline
\hspace*{.25in}
This paper studies for parameters $d, m$ with $d \ge 2, m \ge 1$ and
$(m, d) =1$ the class of maps 
$H: \ZZ^{+} \backslash \ZZ^{+} \to    \ZZ^{+} \backslash d\ZZ^{+}$
having the  the form
$$
H(x) =\frac{1}{d^{a(x)}}( mx - r_j) ~~\mbox{~when~}~~rx \equiv j ~(\bmod~m),~
1 \le j \le d-1,
$$
in which $R(d) := \{ r_j :~ 1 \le j \le d-1 \}$
is any set of integers satisfying
$r_j \equiv mj (\bmod~ d)$ for $1 \le j \le d-1$,
and $d^{a(x)}$ is the maximum power of $d$ dividing
$mx - r_j$. The $3x+1$ problem corresponds to $d=2, m=3$
and $R(d) = \{ -1\}$. The paper shows that if
$$ m \le d^{d/(d-1)} $$
then the set of positive integers $n$ which have some iterate
$H^{(k)}(n) < n$,
has full natural density $1 - \frac{1}{d}$ in the set
$\ZZ^{+} \backslash d\ZZ^{+}.$ He conjectures that 
when $ m \le d^{d/(d-1)} $ the exceptional set of positive 
integers which don't satisfy the conditon is finite.

The author's result generalizes that of Terras(1976) and 
Everett (1977).
In a note added in proof the author asserts that the
proofs of Terras (1976) are faulty. Terras's proofs  seem essentially
correct to me, and in response Terras (1979) provided 
further details.

\item
Helmut A. M\"{u}ller (1991), 
{\em Das `$3n+1$' Problem,}
Mitteilungen der Math. Ges. Hamburg {\bf 12} (1991), 231--251. 
(MR 93c:11053). \\
\newline
\hspace*{.25in}
This paper presents  basic results on $3x+1$ problem, 
with some overlap of Lagarias (1985),
and it  presents  complete proofs
of the results it states.
it contains new observations on the $3x+1$ function $T$ 
viewed as acting on the
 2-adic integers $\ZZ_2$.
For $\alpha  \in  \ZZ_2$ the 2-adic valuation
$| \alpha |_2 =  2^{-j}$ if
$2^j || \alpha$.
M\"{u}ller observes that 
$T( \alpha ) = \sum_{n=0}^\infty (-1)^{n-1} (n+1)2^{n-2} {\alpha \choose n}$
where 
$ {\alpha \choose n}   = \frac{\alpha ( \alpha -1 ) \cdots(\alpha - n+1)}{n!}$.
The function $T ( \alpha )$ is locally constant but not analytic.
Define the function $Q_\infty : \ZZ_2  \rightarrow  \ZZ_2$ by
$Q_\infty ( \alpha ) =  \sum_{i=0}^\infty a_i 2^i$ where
$a_i =  0$ if $| T^{(i)} ( \alpha ) |_2 < 1$
and $a_i =  1$ if $| T^{(i)} ( \alpha ) |_2 =  1$.
Lagarias (1985) showed that this function is a measure-preserving 
homeomorphism
of $\ZZ_2$ to itself.
M\"{u}ller proves that $Q_\infty$ is nowhere differentiable.

\item
Helmut A. M\"uller (1994),
{\em \"Uber eine Klasse 2-adischer Funktionen im Zusammenhang 
mit dem ``$3x+1$''-Problem,}
Abh. Math. Sem. Univ. Hamburg {\bf 64} (1994), 293--302. \\
(MR 95e:11032). \\
\newline
\hspace*{.25in}
Let $\alpha = \sum_{j=0}^\infty a_{0,j} 2^j$ be a 
2-adic integer,
and let $T^{(k)} ( \alpha ) = \sum_{j=0}^\infty a_{k,j} 2^j$ denote its
$k$-th iterate.
M\"uller studies the functions $Q_j : \ZZ_2 \rightarrow \ZZ_2$ defined by
$Q_j ( \alpha ) = \sum_{k=0}^\infty a_{kj} 2^k$.
He proves that each function $Q_j$ is continuous and nowhere differentiable.
He proves that the function $f : \ZZ_2 \rightarrow \QQ_2$ given by
$f = \sum_{j=1}^N A_j Q_j ( \alpha )$ with constants $A_j \in \QQ_2$
is differentiable at a point $\alpha$ with $T^{(k)} ( \alpha ) = 0$ for
some $k \geq 0$ if and only if $2A_0 + A_1 = 0$ and 
$A_2 = A_3 = \ldots = A_N = 0$.

\item
W\l adys\l aw Narkiewicz (1980),
{\em A note on a paper of H. Gupta concerning powers
of 2 and 3,} Univ. Beograd. Publ. Elecktrotech. Fak. Ser. Mat. Fiz. No {\bf 678-715}
 (1980), 173-174.\\
\newline
\hspace*{.25in}
Erd\"os raised the question: ``Does there exist an integer
$m \neq 0,2,8$ such that $2^m$ is a sum of distinct integral powers of 3?''
This was motivated by work of Gupta [
Univ. Beograd. Publ. Elecktrotech. Fak. Ser. Mat. Fiz. Nos. {\bf 602--633}
 (1980), 151--158].
who checked numerically that the
only solutions were $m=0, 2$ and $8$, for $m \le 4734$.
This paper shows that if $N(T)$ denotes the number of such $m \leq T$
then $N(T) \leq 1.62 T^\theta$ where $\theta = \frac{\log 2}{\log 3}$.
%(A weaker result was obtained later in  Marcu (1991), who fails to cite this
%paper.)

\item
J\'"{u}rg Nivergelt (1975),
{\em Computers and Mathematics Education,}
Computers \& Mathematics, with Applications, {\bf 1} (1975), 121--132.\medskip

The paper argues there is a  strong connection between mathematics
education and computers. It points out that mathematics is an
experimental science. In an appendix the  $3x+1$ problem is formulated 
as an example for experimentation, and a number of its properties are derived. 
The author
formulates the $3x+1$ conjecture, says he does not know the answer,
but that it can be numerically studied.
For the Collatz function he gives a heuristic that predicts that the number of elements that
take exactly $s$ steps to iterate to $1$ should grow like $\alpha^s$,
with $\alpha = \frac{1}{2}\left(1 + \sqrt{\frac{7}{3}}\right) \approx 1.26376.$

\medskip
{\em Note.} A similar  discussion of the $3x+1$ problem also appears 
in pages 211--217 of the book: J\"{u}rg Nivergelt, J. Craig Farrar and Edward M. Reingold (1974),
{\em Computer Approaches to Mathematical Problems,}
Prentice-Hall, Inc.: Englewood Cliffs, NJ 1974.

%\item
%J\'{u}rg Nivergelt, J. Craig Farrar and Edward M. Reingold (1974),
%{\em Computer Approaches to Mathematical Problems,}
%Prentice-Hall, Inc.: Englewood Cliffs, NJ 1974.\\
%\newline
%\hspace*{.25in}
%This  textbook  in computer science included a discussion of
%verifying the $3x+1$ problem to some bound  as a programming 
%and data structures challenge. Reference is pages 211--217.

\item
C. Stanley Ogilvy (1972),
{\em Tomorrow's Math: unsolved problems for the amateur},
Second Edition,   Oxford University Press: New York 1972. \\
\newline
\hspace*{.25in}
The $3x+1$ problem is discussed on pages 103-104. He states:
"H. S. M. Coxeter, who proposed it in 1970, stated then that it
had been checked for all $N \le 500,000$. However if the conjecture
is true, which seems likely, a proof will have nothing to do with computers."
He notes that the analogous conjecture for the
$5x+1$ problem is  false, since there is a cycle
not reaching $1$.

{\em Note.} See Coxeter (1971), where Coxeter reports the
problem ``as a piece of mathematical gossip."

\item
Tom\'{a}s Oliveira e Silva (1999),
{\em Maximum Excursion and Stopping Time Record-Holders 
for the $3x+1$ Problem: Computational Results, }
Math. Comp. {\bf 68} No. 1 (1999), 371-384, 
(MR 2000g:11015). \\
\newline
\hspace*{.25in}
This paper reports on computations that verify the $3x+1$
conjecture for $n < 3 \cdot 2^{53} = 2.702 \times 10^{16}$.
It also reports the values of $n$ that are champions for the quantity
$\frac{t(n)}{n}$, where
$$
t(n) : =  \sup_{k \geq  1} T^{(k)} (n) ~.
$$
In this range all $t(n) \leq 8n^2$, which is consistent with the 
conjecture $t(n) \leq n^{2 + o(1)}$
of Lagarias and Weiss (1992).

{\em Note}: In 2004 he
implemented an improved version of this
algorithm. As of 2008 his computation verified the 3x+1 conjecture 
up to $ 19 \times 2^{58} \approx 5.476 \times 10^{18}$. See his webpage at: 
{\tt http://www.ieeta.pt/\~tos/}; email: {\tt tos@ieeta.pt}.  
This is the current record value for verifying the $3x+1$ conjecture.

\item
Elio Oliveri and Giuseppe Vella (1998),
{\em Alcune Questioni Correlate al Problema Del ``3D+1'',} 
Atti della Accademia di scienze lettere e arti di Palermo, Ser. V,
{\bf 28} (1997-98), 21--52. (Italian)\\
\newline
\hspace*{.25in} 
This paper studies the structure of forward and backward iterates
of the $3x+1$ problem, treating the iteration as proceeding from one
odd integer to the next odd iterate.
It obtains necessary conditions  for existence of a
nontrivial cycle in the $3x+1$ problem. It observes
that if $(D_1, ..., D_n)$ are the odd integers in such a cycle, then
$$
\prod_{i=1}^n \frac{3D_i+1}{D_i} = 2^{e_1+ \cdots+ e_n} = 2^k.
$$
If the $3x+1$ conjecture is verified up to $D$ then $\frac{3D+1}{D} > \frac{3D_i+1}{D_i}$, for all $i$, whence
$$
n \log_2 (\frac{3D+1}{2}) > k > n \log_2 3.
$$
They conclude that there must be an integer in the interval
$[n \log_2 (\frac{3D+1}{2}) , n \log_2 3]$, a condition which imposes a lower bound on $n$, 
the number of odd integers in the cycle. In fact one must have 
$k = 1+ \lfloor n \log_2 3 \rfloor,$ 
and $Max_{i}\{  D_{i} \}> (2^{\frac{k}{n}} -3)^{-1},$  
and $\frac{k}{n}$ is a good rational approximation to $\log_2 3$.
Using the continued fraction expansion of
$\log_2 3$, including the intermediate convergents, the authors conclude
that  if the $3x+1$ conjecture
is verified for all integers up to $D=2^{40}+1$,
then  one may conclude $n \ge 1078215.$

The paper also investigates 
the numbers $M(N,n)$ of parity-sequences of length $N$ which
contain $n$ odd values, at positions $k_1, k_2, ..., k_n$ with $k_n=N$, and such that
the partial sums $h_j=k_1+k_2 +\cdots+ k_j$ with all $h_j \le L_j := \lfloor j \log_2 3\rfloor$.
Here $\log_2 (3) \approx 1.585$. 
%They set $Q_{N,n}= \frac{1}{2^N} M(N,n)$.
These are the possible candidate
 initial parity sequences for the smallest positive number in a nontrivial cycle
of the $3x+1$ iteration of period $N$ or longer.

{\em Note: }  The method of the authors for getting a lower bound on
the length of nontrivial cycles is similar in
spirit to earlier methods, cf. Eliahou (1993) and Halbeisen and H\"{u}ngerbuhler (1997). 
 Later computations of the authors show that
if one can take $D= 2^{61}+1\approx 2.306 \times 10^{18}$, then one may conclude there
are no nontrivial cycles containing less than
$n=6,586,818,670$ odd integers.
The calculations of T. Olivera e Silva, as of 2008, have verified
the $3x+1$ conjecture to $19 \times 2^{58} \approx 5.476 \times 10^{18}$, so this improved
cycle bound result is now unconditional.

\item
Clifford  A. Pickover (1989), 
{\em Hailstone $3n  +  1$ Number graphs, }
J. Recreational Math. {\bf 21} (1989), 120--123. \\
\newline
\hspace*{.25in}
This paper gives two-dimensional graphical plots of $3x  +  1$ function 
iterates revealing ``several patterns and a diffuse background 
of chaotically-positioned dots.''

\item
Margherita Pierantoni and Vladan \'{C}ur\v{c}i\'{c} (1996),
{\em A transformation of iterated Collatz mappings,}
Z. Angew Math. Mech. {\em 76}, Suppl. 2 (1996), 641--642. 
(Zbl. 900.65373). \\
\newline
\hspace*{.25in}
A generalized Collatz
map has the form $T(x)= a_i x + b_i$, for $x \equiv i ~(\bmod~n)$,
in which  $a_i= \frac{\alpha_i}{n}, b_i = \frac{\beta_i}{n}$, with $\alpha_i, \beta_i$
integers satisfying $i\alpha_i + \beta_i \equiv 0~(\bmod~n).$ The authors
note there is  a continuous extension of this map to the real line given by
$$
\hat{T}(x) = \sum_{m=0}^{n-1} (a_m x + b_m)
 \left( \sum_{h=0}^{n-1} e^{\frac{2 \pi ih(m-x)}{n}}\right)
= \sum_{h=0}^{n-1} ( A_h x+ B_h) e^{-\frac{2\pi i hx}{n}},
$$
where $\{A_k\}, \{B_k\}$ are the discrete Fourier transforms of the $\{a_i\},$ res. $\{ b_i\}$.
This permits analysis of the iterations of the map using the 
discrete Fourier transform and its inverse. 

They specialize to the case $n=2$, where the data $(a_0, a_1, b_0, b_1)$
describing $T$ are half-integers, with $b_0$ and $a_1+b_1$ integers.
The  $3x+1$ function corresponds to
$(a_0, a_1, b_0, b_1)= (\frac{1}{2}, \frac{3}{2}, 0, \frac{1}{2})$.
One has  $A_0= \frac{1}{2}(a_0+a_1),
A_1= \frac{1}{2}(a_0-a_1), B_0= \frac{1}{2}(b_0+b_1), B_1= \frac{1}{2}(b_0-b_1).$
The authors note that  the recursion 
$$
x_{k+1}= \hat{T}(x_k)= (A_0 x_k+ B_0) + (A_1 x_k+ B_1) \cos( \pi x_k)
$$
can be transformed into  a two-variable system
 using the auxiliary variable $\xi_k= \cos (\pi x_k)$, as
\begin{eqnarray*}
x_{k+1}&=& (A_0x_k + B_0) + (A_1x_k + B_1) \xi_k\\
\xi_{k+1} & =& \cos \left( \pi(A_0x_k+B_0) + \pi(A_1x_k+ B_1) \xi_k\right).
\end{eqnarray*}
They give a formula for $x_{m}$ in terms of the
data $(x_0, \xi_0, \xi_1, ..., \xi_{m-1})$, namely
$$
x_m = \left(\prod_{j=0}^{m-1}(A_0+ A_1 \xi_j)\right)x_0 + (B_0+\xi_{m-1}B_1)+
\sum_{k=0}^{m-2}(B_0 + \xi_kB_1) \prod_{j=k+1}^{m-1}(A_0+ \xi_k A_1).
$$
Now they study the transformed system  in terms of the auxiliary variables $\xi_j$ .
When the starting values $(x_0, x_1)$ are integers, all subsequent values are integers,
and then the auxiliary variables $\xi_k= \pm 1$. Then the  recursion for $\xi_{k+1}$ above
can be simplified by trigonometric sum of angles formulas. In particular the obtain
recursions for the $\xi_j$ that are  independent of the $x_k's$ whenever  $a_0$ and $a_1$
are both integers, but not otherwise. Finally the authors observe  that  on integer orbits 
$\xi_m$ is a periodic function of $x_0$  (of period dividing $2^{m+1}$) which can
be interpolated using a Fourier series  in $\cos \frac{2\pi h}{2^{m+1}}, \sin \frac{2\pi h}{2^{m+1}}$,
using the inverse discrete Fourier transform. They give explicit  
interpolations for $\xi_m$ for the $3x+1$ function for $m=0, 1, 2.$

\item
Nicholas Pippenger (1993),
{\em An elementary approach to some analytic asymptotics,}
SIAM J. Math. Anal. {\bf 24} (1993), No. 5, 1361--1377. (MR 95d:26004). \medskip

%\newline
%\hspace*{.25in} 
This paper studies asymptotics of  recurrences of a type treated in Fredman and Knuth (1974).
It treats those  in \S5 and \S6 of their paper, in which  $g(n) =1$, and we set $h(x):=H(x)-1$.   
These concern  the recurrence
$M(0)=1$, 
$$M(n+1) = 1+ \min_{0 \le k \le n} \left( \alpha M(k) + \beta M(n-k) \right), $$
in the parameter range $\min( \alpha, \beta) > 1$. They reduced the problem  to study of  
the function $h(x)$ satisfying $h(x) = 0,$ $0 < x <1$ satisfying, for $1 \le x < \infty$, 
$$
h(x) = 1+ h( \frac{x}{\alpha}) +h(\frac{x}{\beta}).
$$
For $\alpha, \beta > 1$ let $\gamma$ be
the unique positive solution to $\alpha^{\gamma} + \beta^{-\gamma}=1.$
Fredman and Knuth showed that $h(x) \sim C x^{\gamma}$, when $\frac{\log \alpha}{\log \beta}$
is irrational, and $h(x) \sim D(x) x^{\gamma}$, where $D(x)$ is a periodic function
of the variable $\log x$, if $\frac{\log \alpha}{\log \beta}$ is rational.
Pippenger rederives these results 
using elementary arguments based on  a geometric interpretation of $h(x)$ as  a sum of
binomial coefficients in a triangular subregion of the Pascal triangle. 
He obtains  detailed information on the function $D(x)$ in the case that 
 $\frac{\log \alpha}{\log \beta}$ is rational.

\item
Susana Puddu (1986), 
{\em The Syracuse problem (Spanish), 5$^{\rm  th}$}
Latin American Colloq. on Algebra -- Santiago 1985,
Notas Soc. Math. Chile {\bf 5} (1986), 199--200.
(MR88c:11010). \\
\newline
\hspace*{.25in}
This note considers iterates of the Collatz function $C (x)$.
It shows every positive $m$ has some iterate $C^k (m)  \equiv  1~~(\bmod 4)$.
If $m  \equiv  3~ (\bmod  4)$ the smallest such $k$ must have
$C^k (m)  \equiv  5~ (\bmod~12)$.

\item
Qiu, Wei Xing (1997),
{\em Study on $``3x+1"$ problem} (Chinese),
J. Shanghai Univ. Nat. Sci. Ed.
[Shanghai  da xue xue bao. Zi ran ke xue ban] {\bf 3} (1997), No. 4, 462--464. \\
\newline
\hspace*{.25in} 
English Abstract:  ``The paper analyses the structure presented in the problem $``3x+1"$
and points out there is no cycle in the problem except that $x=1$."
%[The proof of the claimed result is faulty.]

\item
 Raymond Queneau (1963),
 {\em Note compl\'{e}mentaire sur la Sextine,}
 Subsidia Pataphysica, No.1 (1963), 79--80.\\
 \newline
\hspace*{.25in}
Raymond Queneau (1903-1976) was a French poet and novelist,
and a founding member in 1960 of the
French mathematical-literary group Oulipo (Ouvroir de litt\'{e}rature potentielle).
He was  a member of the mathematical society of France from 1948 on, and published
in 1972 a mathematical paper in additive number theory.
 His final essay in 1976 was titled: ``Les fondaments de la litt\'{e}rature d'apr\`{e}s David Hilbert"
(La Biblioth\`{e}que Oulipienne, No. 3); it set out an axiomatic foundation of  literature in
imitiation of Hilbert's {\em Foundations of Geometry}, replacing ``points", ``straight line", and ``plane"
with ``word" , ``sentence"  and ``paragraph", respectively, in some of Hilbert's axioms.
[English translation: The Foundations of Literature (after David Hilbert),
in:  R. Queneau, I. Calvino, P. Fournel,
J. Jouet, C. Berge and H. Mathews, {\em Oulipo Laboratory, Texts from the Biblioth\`{e}que Oulipienne},
Atlas Press: Bath 1995.] He deduces from his axioms:  
``THEOREM 7. {\em Between two words of a sentence there exists an infinity of other words.}" To 
explain this result,  he posits the existence of   ``imaginary words."

This short note is a comment on a preceding article by
A. Taverna, {\em Arnaut Daniel et la Spirale},
Subsidia Pataphysica, No. 1 (1969), 73--78.
Arnaut Daniel was a 12-th century troubadour who
composed poems in Occitan having  a particular rhyming
pattern, called a sestina. Dante admired him and honored
him in several works, including
the {\em Divine Comedy}, where Daniel  is depicted as doing penance in Purgatory
in {\em Purgatorio}. 
The rhyming pattern of the sestina had
six sextets with rhyme pattern involving a cyclic
permutation of order $6$, followed by a triplet.   The work
of Taverna observes that this  cyclic permutation can be represented using a spiral
pattern.   

Queneau considers the ``spiral permutation"  on  numbers $\{1, 2, ..., n\}$
which which takes $2p$ to $p$
and $2p+1$ to $n-p$.
He raises the question: For which $n$ is  a similar spiral permutations
in the symmetric group $S_{n}$ a cyclic permutation?  Call these allowable $n$. 
The example
of the sestina is the case $n=6$; this pattern he terms the sextine.  
He says it is easy to show  that numbers of the form $n= 2xy+x+y$, with $x, y \ge 1$, are not
allowable numbers; this excludes $n=4,7 , 10$ etc.  
He states that  the following $31$ integers $n\le 100$ are allowable:
$n=1, 2, 3, 5, 6, 9, 11, 14, 18, 23, 26,
29, 30, 33, 35, 39, 41, 50, 51, 53, 65, 69, 74, 81, 83, 86, 89, 90, 95,$ \\
$ 98, 99.$

This paper is in the  Oulipo spirit, considering literature obtainable when mathematical
restrictions are placed on its allowable form. Queneau also discussed the 
 topic of  ``Sextines"   in his 1965 essay on
the aims  of  the group Oulipo, 
``Litt\'{e}rature potentielle", published in
{\em B\^{a}tons, chiffres et lettres}, 2nd Edition, Gallimard: Paris 1965.
[English translation: Potential Literature, pp. 181--196 in
 R. Queneau, {\em Letters,
Numbers, Forms: Essays 1928--1970} (Jordan Stump, Translator), Univ. of
Illinois Press: Urbana and Chicago 2007.]

{\em Note.}  Queneau's question was later formulated as the behavior under
iteration of a $(3x+1)$-like function $\delta_n(x)$ on the range
$x \in \{1, 2, ..., n\}$, as observed in Roubaud (1969) and Bringer (1969).
Queneau later published a paper in additive number theory [J. Comb. Theory  A {\bf 12}
 (1972), 31--71], which is on a different topic, but
 mentions in passing   at least one $(3x+1)$-like function  (on page 63).

%{\em Note.} 
\item
Raymond Queneau (1972),
 {\em Sur les suites $s$-additives,}
J. of Combinatorial Theory, Series A, {\bf 12} (1972), 31-71. (MR 46 \# 1741) \\
\newline
\hspace*{.25in}
This is the detailed paper following the announcement
in Comptes Rendus Acad. Sci. Paris (A-B) {\bf 266} (1968), A957--A958.
This paper studies sequences of integers constructed by a
  ``greedy" algorithm, where the first $2s$ integers $0< u_1<u_2< ...,<u_{2s}$ are arbitrary, 
and thereafter each integer $a_t$ is the smallest integer that can be written in
exactly $s$ distinct ways as $S_{ij} = u_i +u_{i+1}+ \cdots + u_j$, for some
$1 \le i < j< t$. He shows that in order for the series not to terminate
the initial set must be a union of two arithmetic progressions of length $s$, 
having the same common difference, one being
$\{u, 2u, ..., su\}$ and the other $\{ v, v+u, ..., v+ (s-1)u\}$. Letting $S(s,u, v)$
denoting the sequence generated this way, he shows that for $s\ge 2$ and
$u \ge 3$, then the sequence is infinite, consisting of
$\{ v+nu: n \ge s\}$, together with the single term $2v+(2s-1)u$.  He then considers
cases with $s\ge 2$ and $u=1$ or $u=2$. Some sequences of the above form
are finite, and some are infinite; he presents some results and conjectures. For $s=0, 1$ all 
$s$-sequences are infinite, and may have a complicated structure. 

 On page 63 there appears a $3x+1$-like function $\sigma(2p)=p-1, \sigma(2p+1)= p$ of
the form considered in Queneau (1963).

\item
Lee Ratzan (1973),
{\em Some work on an Unsolved Palindromic Algorithm},
Pi Mu Epsilon Journal {\bf 5} (Fall 1973), 463--466. \\
\newline
\hspace*{.25in}
The $3x+1$ problem and a generalization were proposed
as an Undergraduate Research Project by David Kay (1972).
The author gives computer code to test the conjecture,
and used it to verify the $3x+1$ Conjecture up to 31,910.
He notices some patterns in two consecutive numbers 
having the same total stopping time and makes
conjectures when they occur. 

\item
Daniel A.  Rawsthorne (1985), 
{\em Imitation of an iteration, }
Math. Magazine {\bf 58} (1985), 172--176.
(MR 86i:40001). \\
\newline
\hspace*{.25in}
This paper proposes a multiplicative random walk that models 
imitating the ``average''
behavior of the $3x  +  1$ function and similar functions.
It compares the mean and standard deviation that this model predicts 
with empirical $3x  +  1$ function data,
and with data for several similar mappings,
and finds good agreement between them.

\item
H. J. J. te Riele (1983a),
{\em Problem 669,}
Nieuw Archief voor Wiskunde, Series IV, {\bf 1}  (1983), p. 80.\\
\newline
\hspace*{.25in}
This problem asks about the iteration of the function
$$
f(x) = \left\{
\begin{array}{cl}
\frac{n}{3}  & \mbox{if} ~~n \equiv 0~~ (\bmod ~ 3) \\
~~~ \\
\lfloor n \sqrt{3} \rfloor  & \mbox{if}~~ n \not\equiv 0~~ (\bmod ~3 )~.
\end{array}
\right.
$$
on the positive integers.  The value $n=1$ is a fixed point , 
 and all powers $n= 3^k$ eventually reach this fixed point.
 The problem asks to show that if two consecutive
 iterates are not congruent to $0~(\bmod~3)$ then the
 trajectory of this orbit thereafter grows monotonically,
 so diverges to $+\infty$; otherwise the orbit reaches $1$.

{\em Note.} Functions similar in form to $f(x)$, but with 
a condition $(\bmod~2)$,  were studied by
Mignosi (1995) and Brocco (1995).

\item
H. J. J. te Riele (1983b),
{\em Iteration of Number-theoretic functions},
Nieuw Archief voor Wiskunde, Series IV,  {\bf 1} (1983), 345--360.
[MR 85e:11003].
\newline
\hspace*{.25in}
The author surveys  work on a wide variety of number-theoretic functions
which take positive integers to positive integers, whose behavior under iteration is not understood.
This includes the $3x+1$ function (III.1) and the $qx+1$ function (III.2).
The author presents as Example 1  the function in te Riele (1983a), 
$$
f(x) = \left\{
\begin{array}{cl}
\frac{n}{3}  & \mbox{if} ~~n \equiv 0~~ (\bmod ~ 3) \\
~~~ \\
\lfloor n \sqrt{3} \rfloor  & \mbox{if}~~ n \not\equiv 0~~ (\bmod ~3 )~.
\end{array}
\right.
$$
He showed that this function has divergent trajectories, and
that all non-divergent orbits  converge to the fixed point $1$. He  
conjectures that almost all orbits on positive integers tend to $+ \infty$.
Experimentally only 459 values of $n < 10^5$ 
have orbits converging to $1$.

\item
Jacques Roubaud (1969) ,
{\em Un probl\`{e}me combinatoire pos\`{e} par la po\'{e}sie lyrique des troubadours,}
 Math\'{e}matiques et Sciences Humaines
[Mathematics and Social Science], {\bf 27} Autumn 1969, 5--12.\\
\newline
\hspace*{.25in}
This paper addresses  the question of  suitable rhyming schemes for poems,
suggested by the schemes used by medieval
troubadours, and raised in Queneau (1963). This leads to combinatorial questions concerning
the permutation structure of such rhyming schemes.  The author classifies the movement
of rhymes by permutation patterns. In section 5 he formulates  three questions
suggested by the rhyme pattern of the 
 sestina of Arnaut Daniel. One of these (problem ({\em Pa})) is the question of 
Queneau (1963), which he states as iteration of the $3x+1$-like function
$$
\delta_n(x) := \left\{ 
 \begin{array}{cl}
\df{x}{2} & ~~\mbox{if} ~~x  ~~\mbox{is ~even}\\
~&~\\
\df{2n+ 1-x}{2} & ~~\mbox{if} ~~x  ~~\mbox{is ~odd}\\
\end{array}
%\.
\right.
$$
restricted to the domain  $\{1, 2, ..., n\}$. This paper goes together with the analysis of 
this function done by his student   Monique Bringer (1969).  
In  Robaud (1993) he proposes  further generalizations of these problems.

{\em Note.} Jacques Roubaud is a mathematician and  a member of the
literary group Oulipo. In 1986 he wrote a spoof, an obiturary for N. Bourbaki
(see pages 73 and 115 in: M. Mashaal, {\em Bourbaki} (A. Pierrehumbert, Trans.), Amer. Math. Soc.,
Providence 2006).
He discussed the mathematical work of Raymond Queneau in 
J. Roubaud, La mathematique dans la methode de Raymond Queneau,
Critique: revue g\'{e}n\'{e}rale des publications fracaises et \'{e}trang\`{e}res 
{\bf 33} (1977), no. 359, 392--413.
[English Translation: Mathematics in the Method of Raymond Queneau, pp. 79--96
in: Warren F. Motte, Jr. (Editor and Translator), {\em Oulipo, A Primer of
Potential Literature}, Univ. of Nebraska Press: Lincoln, NB 1986.]

\item
Jacques Roubaud (1993),
{\em $N$-ine, autrement dit quenine (encore)},
{\em R\'{e}flexions historiques et combinatoires sur la $n$-ine,
autrement dit quenine.}
La Biblioth\'{e}que Oulipienne, num\'{e}ro 66, 
Rotographie \`{a} Montreuil (Seine-Saint-Denis), 
November 1993.\\
\newline
\hspace*{.25in} 
This paper considers work on rhyming patterns in poems suggested
by those of medieval troubadours. The $n$-ine is a generalization
of the rhyme pattern of the sestina of Arnaut Daniel, a 12-th century
troubadour, see Queneau (1963), Roubaud (1969). It considers a
``spiral permutation" of the symmetric group $S_{N}$, which he earlier gave
mathematically as 
$$
\delta_n(x) := \left\{ 
 \begin{array}{cl}
\df{x}{2} & ~~\mbox{if} ~~x  ~~\mbox{is ~even}\\
~&~\\
\df{2n+ 1-x}{2} & ~~\mbox{if} ~~x  ~~\mbox{is ~odd}\\
\end{array}
%\.
\right.
$$
restricted to the domain  $\{1, 2, ..., n\}$.

Roubaud summarizes the work of Bringer (1969) giving restrictions
on admissible $n$. She observed $p=2n+1$ must be prime and
that $2$ is a primitive root of $p$ is a sufficient condition for
admissibility. 
He advances the conjectures that 
the complete set of admissible  $n$ are those with  $p=2n+1$ a 
prime, such that either $ord_p(2)=2n$, so that $2$ is a primitive root, or else 
$ord_p(2)=n$.  (This conjecture later turned out to require a correction;
namely, to allow those $n$ with  $ord_p(2)=n$ only when $n \equiv 3 (\bmod~4)$.) 

Roubaud suggests various generalizations of the problem.
Bringer(1969) had studied the inverse permutation to $\delta_n$, given by
$$
d_n(x) := \left\{ 
 \begin{array}{cl}
2x & ~~\mbox{if} ~~1 \le x \le \frac{n}{2}\\
~&~\\
2n+1-2x & ~~\mbox{if} ~~\frac{n}{2} < x \le n\\
\end{array}
%\.
\right.
$$
In Section 5 Roubaud suggests a generalization of this, to the permutations
$$
d_n(x) := \left\{ 
 \begin{array}{cl}
3x & ~~\mbox{if} ~~1 \le x \le \frac{n}{3}\\
~&~\\
2n+1-3x & ~~\mbox{if} ~~\frac{n}{3} < x \le \frac{2n}{3}\\
~&~\\
3x- (2n+1) & ~~\mbox{if} ~~\frac{2n}{3} < x \le n.\\
\end{array}
%\.
\right.
$$
Here $n=8$ is a solution, and he composes a poem, ``Novembre", a $3$-octine,
to this rhyming pattern, in honor of Raymond Queneau. 
Call these solutions
$3$-admissible, and those of the original problem $2$-admissible. 
He gives the following list of $3$-admissible solutions for $n \le 200$ that are not $2$-admissible:
$n=8, 15, 21, 39, 44, 56, 63, 68, 111, 116, 128, 140, 165, 176, 200.$

He concludes with some other proposed
rhyming schemes involving spirals, including the ``spinine", also called  ``escargonine".

\item
Olivier Rozier (1990), 
{\em Demonstraton de l'absence de cycles d'une
certain forme pour le probl\'{e}me de Syracuse,}
Singularit\'{e} {\bf 1} no. 3 (1990), 9--12. \\
\newline
\hspace*{.25in}
This paper proves that there are no cycles except the trivial cycle 
whose iterates $(\bmod~2)$
repeat a pattern of the form $1^m 0^{m'}$.
Such cycles are called
{\em circuits}
by R. P. Steiner, who obtained this result in 1978.
Rozier's proof uses effective transcendence bounds similar to
Steiner's, but is simpler since he can quote the recent bound
$$
\left| \df{\log 3}{\log 2}  -  \frac{p}{q} \right|  >  q^{-15}  , 
$$
when $(p,  q) =  1$,
which appears in M. Waldschmidt, Equationes 
diophantiennes et Nombres Transcendents,
Revue Palais de la D\'ecouverte,
{\bf 17},
No. 144 (1987) 10--24.

\item
Olivier Rozier (1991), 
{\em Probleme de Syracuse:
``Majorations'' Elementaires des Cycles,} 
Singularit\'e {\bf 2}, no. 5 (1991), 8--11. \\
\newline
\hspace*{.25in}
This paper proves  that if the integers in a cycle of the 
$3x +  1$ function are
grouped into blocks of integers all having the same parity,
for any $k$ there are only a finite number of cycles having $ \leq  k$
blocks.
The proof uses the Waldschmidt result
$\left| \frac{\log 3}{\log 2}  -  \frac{p}{q} \right|  >  q^{-15}$.

\item
J\"{u}rgen W. Sander (1990), 
{\em On the $(3N  +  1)$-conjecture,}
Acta Arithmetica {\bf 55} (1990), 241--248.
(MR 91m:11052). \\
\newline
\hspace*{.25in}
The paper shows that the number of integers $\leq x$ for which
the $3x  +  1$ Conjecture is true is at least $x^{3/10}$, by extending the
approach of Crandall (1978).
(Krasikov (1989) obtains a better lower bound $x^{3/7}$ by another method.)
This paper also shows that if the $3x +  1$ Conjecture is true for all
$n  \equiv  \frac{1}{3} ( 2^{2k} - 1)~ (\bmod  ~2^{2k-1} )$,
for any fixed $k$,
then it is true in general.

\item
Benedict G. Seifert (1988), 
{\em On the arithmetic of cycles for the Collatz-Hasse 
(`Syracuse') conjectures,}
Discrete Math. {\bf 68} (1988), 293--298.
(MR 89a:11031). \\
\newline
\hspace*{.25in}
This paper gives
 criteria for cycles of $3x  +  1$ function to exist, and bounds the
smallest number in the cycle in terms of the length of the cycle.
Shows that if the $3x +  1$ Conjecture is true, then the only 
positive integral solution of
$2^l -3^r =1$ is $l=2, r=1$.

[All integer solutions of $2^l -3^r =  1$ have been found
unconditionally. This can be done  by a method of
St$\o$rmer, Nyt. Tidsskr. Math. B {\bf 19} (1908) 1--7.
It was done as a special case of various more general
results by S. Pillai, Bull. Calcutta Math. Soc. {\bf 37} (1945), 15--20 
(MR\# 7, 145i);
W. LeVeque, Amer. J. Math. {\bf 74} (1952), 325--331 (MR 13, 822f);
R. Hampel, Ann. Polon. Math. {\bf 3} (1956), 1--4 (MR 18, 561c); etc.]

\item
Jeffrey  O. Shallit and David W. Wilson (1991), 
{\em The ``$3x+1$'' Problem and Finite Automata, }
Bulletin of the EATCS
(European Association for Theoretical Computer Science),
No. 46, 1991, pp. 182--185. \\
\newline
\hspace*{.25in}
A set $S$ of positive integers is said to be
{\em 2-automatic}
if the binary representations of the integers in $S$ form a
regular language $L_S  \subseteq \{  0,1 \}^\ast$.
Let $S_i$ denote the set of integers $n$ which have some
$3x+1$ function iterate $T^{(j)} (n) =  1$, and whose 
$3x+1$ function iterates include exactly $i$ odd integers $\geq  3$.
The sets $S_i$ are proved to be 2-automatic for all $i  \geq  1$.

\item
Daniel Shanks (1965),
{\em Comments on Problem $63-13^{*}$},
SIAM Review {\bf 7} (1965), 284--286. \\
\newline
\hspace*{.25in}
This note gives comments on Problem $63-13^{*}$,
proposed by   Klamkin (1963). This 
problem concerns the
iteration of Collatz's original function, which
is a permutation of the integers.
He states that these problems date back at least
to 1950,  when L. Collatz mentioned them in
personal conversations at the  International
Math. Congress held  at Harvard University.
He gives the results of a computer search of the orbit of
$8$ for Collatz's original function, observing
that it reaches numbers larger than $10^{10}.$
 He observes there are known cycles of length 
$1$, $2$, $5$ and $12$, the last having
smallest element $n=144$. He observes that
this seems related to the fact that the
continued fraction expansion of $\log_2 3$
has initial convergents  having denominators
$1$, $2$, $5$, $12$, $41$, ...
However he does not know of any  cycle of period $41$. He
notes that it is not  known whether the only cycle lengths
that can occur must be denominators of such
partial quotients. Later Atkin (1966) proved
there exists no cycle of period $41$.

\item
Daniel Shanks (1975), 
{\em Problem $\#5$},
Western Number Theory Conference 1975, Problem List,
(R. K. Guy, Ed.). \\
\newline
\hspace*{.25in}
The problem concerns iteration of the Collatz function
$C(n)$. Let $l(n)$ count the number of distinct integers
appearing in the sequence of iterates $C^{(k)}(n)$, $k \ge 1$, 
assuming it eventually enters a cycle.
Thus $l(1) = 3, l(2)=3, l(3) = 8$, for example. Set
$S(N) = \sum_{n=1}^N l(n)$. The problem asks 
whether it is  true that
$$
S(N) = A N \log N + B N + o(N) ~~\mbox{as}~~ N \to \infty,
$$
where
$$ A = \frac{3}{2} \log \frac{4}{3} \approx 5.21409 \mbox{~and~}
B = A(1 - \log 2) \approx 1.59996.
$$

In order for $S(N)$ to
remain finite, there must be no divergent trajectories.
This problem formalizes the result of a heuristic
probabilistic calculation based on assuming the
$3x+1$ Conjecture to be true. 

\item
Ray P. Steiner (1978),
{\em A Theorem on the Syracuse Problem},
Proc. 7-th Manitoba Conference on Numerical Mathematics
and Computing (Univ. Manitoba-Winnipeg 1977),
Congressus Numerantium XX, Utilitas Math.: Winnipeg, 
Manitoba 1978, pp. 553--559.
(MR 80g:10003). \\
\newline
\hspace*{.25in}
This paper studies periodic orbits of the
$3x+1$ map, and a
problem raised by Davidson (1976).
 A sequence of iterates  $\{n_1, n_2, ... , n_p, n_{p+1}\}$ 
with $T(n_j) = n_{j+1}$ is called by Davidson (1977) a
{\em circuit} if it consists of a sequence of
odd integers $\{n_1, n_2,..., n_j\}$ followed by a sequence
of even integers $\{n_{j+1}, n_{j+2}, ..., n_p\}$,
with $n_{p+1} = T(n_p)$ an odd integer.
A circuit is a {\em cycle} if $n_{p+1} = n_1$.

This paper shows that  the only  circuit on
the positive integers that is a cycle is $\{1, 2\}$.
It uses the observation of Davison (1977) 
that these corresponds to positive solutions $(k, l, h)$ 
to the exponential Diophantine equation
$$ 
(2^{k+l} - 3^k)h = 2^l - 1.
$$
The paper  shows that the only solution to this equation
in positive integers is  $(k, l, h) = (1,1,1)$.
The proof  uses results from transcendence theory, 
Baker's method of linear forms in logarithms
(see A. Baker, {\em Transcendental Number Theory},
Cambridge Univ. Press 1975, p. 45.)

[One could prove similarly that there are exactly three
circuits that are cycles on the non-positive integers,
namely  $\{-1\}$, and $\{ -5, -7, -10\}$.
These correspond to the solutions to the exponential
Diophantine equation 
 $(k, l, h) =  (k, 0, 0)$ for any $k \ge 1$; and $ (2,1,-1)$,
respectively. A further solution  $(0, 2, 1)$
corresponds to the cycle $\{0\}$ which is not a circuit
by the definition above. ]

\item
Ray P. Steiner (1981a),
{\em On the ``$Qx+1$'' Problem, $Q$ odd},
Fibonacci Quarterly {\bf 19} (1981), 285--288.
(MR 84m:10007a) \\
\newline
\hspace*{.25in}
This paper studies the $Qx+1$-map
$$
h(n) = \left\{
\begin{array}{cl}
\frac{Qn+1}{2}  & \mbox{if} ~~n \equiv 1~~ (\bmod ~ 2), n > 1 \\
~~~ \\
\frac{n}{2} & \mbox{if}~~ n \equiv 0~~ (\bmod ~2 )~ \\
~~~\\
1  & \mbox{if}~~ n=1,
\end{array}
\right.
$$
when $Q$ is odd and $Q >3$.
It proves that the only circuit which is a cycle  when $Q= 5$ is
$\{13, 208\}$, that there is no circuit which is a cycle for  $Q=7.$
Baker's method is again used, as in Steiner (1978).

\item
Ray P. Steiner (1981b),
{\em On the ``$Qx+1$'' Problem, $Q$ odd II},
Fibonacci Quarterly {\bf 19} (1981), 293--296.
(MR 84m:10007b). \\
\newline
\hspace*{.25in}
This paper continues to study  the $Qx+1$-map
$$
h(n) = \left\{
\begin{array}{cl}
\frac{Qn+1}{2}  & \mbox{if} ~~n \equiv 1~~ (\bmod ~ 2), n > 1 \\
~~~ \\
\frac{n}{2} & \mbox{if}~~ n \equiv 0~~ (\bmod ~2 )~ \\
~~~ \\
1  & \mbox{if}~~ n=1,
\end{array}
\right.
$$
when $Q$ is odd and $Q >3$. It makes general remarks on
the case $Q >7$, and presents data on from the computation of
$\log_2 \frac{5}{2}$ and $\log_2 \frac{7}{2}$ used in the
proofs in part I.

\item
Rosemarie M. Stemmler (1964),
{\em The ideal Waring problem for exponents $401-200,000$,}
Math. Comp. {\bf 18} (1964), 144--146. (MR 28 \#3019)\\
\newline
\hspace*{.25in}
The ideal Waring theorem states  that for a given $k \ge 2$ each positive integer is
the sum of $2^k + \lfloor (\frac{3}{2}^k \rfloor$ non-negative $k$-th powers,
provided that the fractional part $\{\{ x\}\} := x - \lfloor x \rfloor$ of $(\frac{3}{2})^k$
satisfies
$$
0 \le \{\{ (\frac{3}{2})^k \}\} < 1 - (\frac{3}{4})^k.
$$
This paper checks that this inequality holds for $401 \le k \le 200000$.
This problem on the fractional parts of powers of $(\frac{3}{2})^k$ motivated
the work of  Mahler (1968) on $Z$-numbers.

\item
Kenneth S. Stolarsky (1998), 
{\em A prelude to the $3x+1$ problem,}
J. Difference Equations and  Applications {\bf 4} (1998), 451--461.
(MR 99k:11037). \\
\newline
\hspace*{.25in}
This paper studies a purportedly ``simpler''
analogue of the $3x+1$ function.
Let $\phi = \frac{1+ \sqrt 5}{2}$.
The function $f: \ZZ^+ \rightarrow \ZZ^+$
is given by
$$
\left\{
\begin{array}{ll}
f( \lfloor n \phi \rfloor ) = \lfloor n \phi^2 \rfloor +2~,  
& \mbox{all} ~~ n \geq 1 ~, \\
~~~ \\
f ( \lfloor n \phi^2 \rfloor )  =  n~. & \mbox{all} ~~ n \geq 1 ~.
\end{array}
\right.
$$
This function is well-defined because the sets
$A = \{ \lfloor n \phi \rfloor : n \geq 1 \}$ and
$B = \{ \lfloor n \phi^2 \rfloor : n \geq 1 \}$
form a partition of the positive
integers.
(This is a special case of Beatty's theorem.)
The function $f$ is analogous to the $3x+1$ function in that it is
increasing on the domain $A$ and decreasing on the
domain $B$.
The paper shows that almost all trajectories diverge to $+ \infty$,
the only exceptions being a certain set $\{3,7,18,47, \ldots \}$
of density zero which converges to the two-cycle $\{3,7\}$.
Stolarsky determines the
complete symbolic dynamics of $f$.
The possible symbol sequences of orbits are
$B^\ell (AB)^k A^\infty$ for some $\ell , k \geq 0$,
for divergent trajectories, and
$B^\ell (AB)^\infty$ for some $\ell \geq 0$, for
trajectories that reach the two-cycle $\{3,7\}$.

\item
Gy\"{o}rgy Targ\'{o}nski (1991), 
{\em Open questions about KW-orbits and iterative roots, }
Aequationes Math. {\bf 41} (1991), 277--278. \\
\newline
\hspace*{.25in}
The author suggests the possibility of applying a result of H. Engl (An
analytic representation for self-maps of a countably infinite set and its
cycles,
Aequationes Math.
{\bf 25}
(1982), 90-96, MR85d.04001) to bound the
number of cycles of the $3x+1$ problem.
Engl's result expresses the number of cycles 
as the geometric multiplicity of 1
as an eigenvalue of a map on the sequence space $l^1$.

\item
Riho Terras (1976),
{\em A stopping time problem on the positive integers,}
Acta Arithmetica {\bf 30} (1976), 241--252. 
(MR 58 \#27879). \\
\newline
\hspace*{.25in}
This is the first significant research paper 
to appear that deals directly
with the $3x+1$ function. The $3x+1$
function was however the motivation for the paper
Conway (1972).
The main result of this paper was obtained independently 
and contemporaneously  by Everett (1977).

A positive integer $n$ is said to have  
{\em stopping time} $k$ if the $k$-th iterate $T^{(k)}(n) < n,$
and $T^{(j)}(n) \ge n$ for $1 \le j < k$. 
The author shows that the set of integers having
stopping time $k$ forms a set of congruence classes
$(\bmod~ 2^k)$, minus a finite number of elements.
He shows that the set of integers having a  finite
stopping time has natural density one.
Some further 
details of this proof were supplied later in Terras (1979).

This paper introduces the notion of the {\em coefficient
stopping time} $\kappa(n)$ of an integer $n > 1.$
Write  $T^{(k)}(n) = \alpha(n) n + \beta(n)$
with $\alpha(n) = \frac{3^{a(n)}}{2^n}$,
where $a(n)$ is the number of iterates $T^{(j)}(n) \equiv 1 (\bmod~2)$
with $0 \le j < k$.
Then $\kappa(n)$ is 
defined to be the least
$k \ge 1$ such that
$T^{(k)}(n) = \alpha(n) n + \beta(n)$ has $\alpha(n) < 1$,
and  $\kappa(n) =\infty$ if no such value exists.
It is clear that $\kappa(n) \le \sigma(n)$,
where $\sigma(n)$ is the stopping time of $n$.
Terras formulates the {\em Coefficient Stopping Time
Conjecture}, which asserts that $\kappa(n) = \sigma(n)$
for all $n \ge 2$. He proves this conjecture
for all values  $\kappa(n) \le 2593$.
This can be done for  $\kappa(n)$ below a fixed
bound  by upper bounding $\beta(n)$ and showing 
bounding $\kappa(n) < 1 - \delta$ for suitable $\delta$
and determining  the maximal  value $\frac{3^l}{2^j} < 1$
possible with $j$ below the given bound. The 
convergents of the continued
fraction expansion of $\log_2 3$ play a role in 
determining the values of $j$ that must be checked.

\item
Riho Terras (1979),
{\em On the existence of a density}, 
Acta Arithmetica {\bf 35} (1979), 101--102.
(MR 80h:10066). \\
\newline
\hspace*{.25in}
This paper supplies additional details concerning
 the proof in Terras (1976) that
the set of integers having an infinite stopping time
has asymptotic density zero.
The proof in Terras (1976)
had been criticized by M\"{o}ller (1978).

\item
Bryan Thwaites (1985), 
{\em My conjecture,}
Bull. Inst. Math. Appl. {\bf 21} (1985), 35--41. \\
(MR86j:11022). \\
\newline
\hspace*{.25in}
The author states that he invented the
$3x  +  1$ problem in 1952.
He derives basic results about iterates ,and
makes conjectures on the ``average'' behavior of trajectories.

{\em Note.} The mathematical community generally credits 
L. Collatz as being  the first to propose the $(3x  +  1)$-problem, see 
Collatz (1986), and the comment of Shanks (1965). 
This would be an independent discovery of the problem.

\item
Bryan Thwaites (1996),
{\em Two Conjectures, or how to win $\pounds$ 1000,} 
Math. Gazette {\bf 80} (1996), 35--36. \\
\newline
\hspace*{.25in}
One of the two conjectures is the $3x+1$-problem, for which
the author offers the stated reward.

\item
Robert Tijdeman (1972),
{\em Note on Mahler's $3/2$-problem},
Det Kongelige Norske  Videnskabers Selskab Skrifter No. 16, 1972, 4 pages.
(Zbl. 227: 10025.) \\
\newline
\hspace*{.25in}
This paper concerns the Z-number problem of Mahler (1968), which asks
whether  there
exists any nonzero real number $\eta$ such that 
the fractional parts $0 \le \{\{\eta (\frac{3}{2})^n\} \} \le \frac{1}{2}$
for all $n \ge 0.$ By an elementary argument Tijdeman shows
that analogues of $Z$-numbers  exist in a related problem.
Namely, for every $k \ge 2$ amd $m \ge1$ there exists a real number
$\eta \in [m, m+1)$ such that
$0\le \{\{ \eta (\frac{2k+1}{2})^n\}\} \le \frac{1}{2k+1}$ holds
for all $ n \ge 0$.

\item
Charles W. Trigg, Clayton W. Dodge and Leroy F. Meyers (1976),
{\em Comments on Problem 133},
Eureka (now Crux Mathematicorum) 
{\bf 2}, No. 7 (August-Sept.) (1976), 144--150. \\
\newline
\hspace*{.25in}
Problem 133 is the $3x+1$ problem. 
It was proposed  by K. S. Williams 
(Concordia Univ.), who said that he was shown it by one
of his students. 
C. W. Trigg gives some earlier history of the
problem. He remarks that Richard K. Guy wrote to him stating 
that Lothar Collatz had given a lecture on the problem at Harvard in
1950 (informally at  the International Math. Congress).
He  reported that in 1970 H. S. M. Coxeter offered
a prize of $\$50$ for proving the $3x+1$ Conjecture and 
$\$100$ for finding a counterexample, in his talk:
``Cyclic Sequences and Frieze Patterns''
(The Fourth Felix Behrend Memorial Lecture in Mathematics),
The University of Melbourne,  1970, see Coxeter (1971).
He also referenced  a discussion of the problem in several issues
of {\em Popular Computing} No. 1 (April 1973) 1--2;
No. 4 (July 1973) 6--7; No. 13 (April 1974) 12--13;
No. 25 (April 1975), 4--5.
Dodge references the work of Isard and Zwicky (1970).

{\em Note.} Lothar Collatz was present at the 1950 ICM
as part of the DMV delegation, while
Coxeter, Kakutani and Ulam each delivered a lecture
at the 1950 Congress that appears in its proceedings.

%\item
%Evangelos Tzanakis (2003^+),
%{\em Collatz conjecture: properties and algorithms},
%preprint. \\
%\newline
%\hspace*{.25in}
%This paper is concerned with algorithms for testing
%the $3x+1$ Conjecture. The natural numbers are divided into
%classes according to their iteration behavior, with
%each class labelled by a binary tree. A search algorithm
%for testing the $3x+1$ conjecture up to some bound is
%presented.

\item
Toshio Urata  and Katsufumi Suzuki (1996),
{\em Collatz Problem}, (Japanese),
Epsilon (The Bulletin of the Society of Mathematics Education of Aichi Univ. of Education) 
[Aichi Kyoiku Daigaku Sugaku Kyoiku Gakkai shi]
{\bf 38} (1996), 123--131.\\
\newline
\hspace*{.25in}
The authors consider iterating the Collatz function $C(x)$ (denoted $f(x)$) and the
speeded-up Collatz function $\phi(x)$ mapping  odd integers to odd integers by 
$$
\phi(x) = \frac{3x+1}{2^e},~~\mbox{with} ~2^e || 3x+1.
$$ 
They let $i_n$ (resp. $j_n$) denote the number of iterations to get from $n$ to $1$ of
the Collatz function (resp. the function $\phi(x)$). They study the Cesaro means
$$
h_n := \sum_{k=1}^n i_n, ~~~ m_n := \sum_{k=1}^n j_k.
$$
and observe empiricially that
$$
h_n \approx 10.4137\log n  - 12.56
$$
$$
m_n \approx 2.406 \frac{\log n}{\log 2}  -2
$$
fit the data up to $4 \times 10^8$. Then they introduce an entire function which interpolate
the Collatz function at positive integer values
$$
f(z) = \frac{1}{2} z (\cos \frac{\pi z}{2})^2 + (3z+1) (\sin \frac{\pi z}{2})^2.
$$
They also introduce more complicated entire functions that interpolate the function $\phi(x)$
at odd integer values. 
$$
F(z) = \frac{4}{\pi^2} (\cos \frac{\pi z}{2})^2 \left( \sum_{n=0}^{\infty}[ 
 \alpha_n \frac{1}{(z-(2n+1))^2} - \frac{1}{(2n+1)^2}] \right)
$$
in which $\alpha(n):= \phi(2n+1)$.
They raise the problem of determining the Fatou and Julia sets of these functions.
These functions are constructed so that the positive  integer values (resp. positive odd integer values)
fall in the Fatou set.

\item
Toshio Urata, Katsufumi Suzuki and Hisao Kajita (1997),
{\em Collatz Problem II}, (Japanese),
Epsilon (The Bulletin of the Society of Mathematics Education of Aichi Univ. of Education) 
[Aichi Kyoiku Daigaku Sugaku Kyoiku Gakkai shi]
{\bf 39} (1997), 121--129.\\
\newline
\hspace*{.25in}
The authors study the speeded-up Collatz function $\phi(x)= \frac{3x+1}{2^e}$ which takes
odd integers to odd integers. They observe that every orbit of $\phi(x)$ contains
some integer $n \equiv 1~(\bmod~4).$ They then introduce entire functions
that interpolate $\phi(x)$ at positive odd integers. They start with
$$
F(z) = \frac{4}{\pi^2} (\cos \frac{\pi z}{2})^2\left( \sum_{n=0}^{\infty}[ 
 \alpha_n \frac{1}{(z-(2n+1))^2} - \frac{1}{(2n+1)^2} ] \right)
$$
in which $\alpha(n):= \phi(2n+1)$, so that
$F(2n+1)= \phi(2n+1)$ for $n \ge 0$. They observe  $F(z)$ has an attracting 
fixed point at $z_0=-0.0327$, a repelling fixed point at $z=0$, and a
superattracting fixed point at $z=1$. 
They show that the positive odd integers are in the Fatou set of $F(z)$.
They show that the immediate basins of the 
Fatou set around the positive odd integers of $F(z)$ are disjoint.
Computer drawn pictures are included, which include small copies of
sets resembling the Mandlebrot set.

The authors  also  introduce, for each integer $p \ge 2$, the entire functions 
$$
K_p(z) = \left( \frac{4}{\pi^2} (\cos \frac{\pi z}{2})^2\right)^p 
\left( \sum_{n=0}^{\infty} \frac{\alpha_n}{(z-(2n+1))^{2p}}\right),
$$
and, for $p \ge 1$, the entire functions. 
$$
L_p(z) = 
\left( \frac{4}{\pi^2} \cos^2 \frac{\pi z}{2} \right)^p \frac{\sin \pi z}{\pi}
\left( \sum_{n=0}^{\infty} \frac{\alpha_n}{ (z-(2n+1))^{2p+1} } \right)
$$
These functions also interpolate $\phi(x)$
at odd integers, i.e. $\phi(2n+1)= K_p(2n+1)= L_p(2n+1)$ for $n \ge 0$.

%[add more]
\item
Toshio Urata  and Hisao Kajita (1998),
% Suzuki second author cited  in Urata (2002), actual Japanese paper lists  Kajita as author
{\em Collatz Problem III}, (Japanese),
Epsilon (The Bulletin of the Society of Mathematics Education of  Aichi Univ. of Education)
[Aichi Kyoiku Daigaku Sugaku Kyoiku Gakkai shi]
{\bf 40} (1998), 57--65.\\
\newline
\hspace*{.25in}
The authors study the speeded-up Collatz function $\phi(x)= \frac{3x+1}{2^e}$ which takes
odd integers to odd integers. They describe the infinite number of preimages of a
given $x$ as $\{\frac{4^n x-1}{3}:~n \ge 0\}$
 if $x \equiv 1~(\bmod~3)$ and
$\{\frac{4^n (2x)-1}{3}: ~n \ge 0\}$
 if $x \equiv 2~(\bmod~3)$.
 They encode the trajectory of this function with a vector $(p_1, p_2, p_3, ...)$ which keeps
 track of the powers of $2$ divided out at each iteration.
 More generally they consider the speeded-up $(ax+d)$-function, with $a, d$ odd.
 For   $y_0$ an initial value (with $y_0 \ne -\frac{d}{a}$), 
 the iterates are $y_i= \frac{a y_{i-1} + d}{2^{p_{i+1}}}$.
They encode this  iteration for $n$ steps by an $n \times n$ matrix
 $$
\Theta:= \left[ \begin{array}{cccccc}
 2^{p_1} & 0 &0 & \cdots & 0 & -a\\
 -a & 2^{p_2} & 0 & \cdots & 0 & 0 \\
 0 & -a & 2^{p_3} &     \cdots    &  0  & 0 \\
  & \cdots &&\cdots &&\\
  0&0&0& \cdots &-a& 2^{p_n} \\
 \end{array}
 \right],
$$ 
which acts by 
$$
\Theta \left[ \begin{array}{c} 
y_1\\ y_2 \\ y_3\\ \cdots\\ y_n\\ 
\end{array} \right] =
 \left[ \begin{array}{c} 
d - a(y_n-y_0) \\ d\\ d\\ \cdots \\d \\ 
\end{array} \right] .
$$
The authors observe that $\det(\Theta) = 2^{p_1+\cdots + p_n} - a^n.$
A  periodic orbit, one with $y_n=y_0$,
corresponds to a vector of iterates mapping to the vector with
constant entries $d$. They compute examples for  periodic orbits of $5x+1$ problem, and the $3x+1$ problem
on negative integers.
Finally they study orbits with $y_0$ a rational number with denominator
relatively prime to $a$, and give some examples of periodic orbits. 

\item
Toshio Urata (1999),
{\em Collatz Problem IV}, (Japanese),
Epsilon (The Bulletin of the Society of Mathematics Education of Aichi Univ. of Education)
[Aichi Kyoiku Daigaku Sugaku Kyoiku Gakkai shi]
 {\bf 41} (1999), 111-116.\\
\newline
\hspace*{.25in}
The author studies a $2$-adic interpolation of the speeded-up Collatz function $\phi(n)$
defined on odd integers $n$ by 
dividing out all powers of $2$, i.e. for an odd integer $n$, 
$\phi(n) = \frac{3n+1}{2^p(3n+1)}$,
where $p(m)= ord_2(m).$  Let $\ZZ_2^{\ast}=\{ x \in \ZZ_2: ~x \equiv 1 ~(\bmod~2)\}$ denote the
$2$-adic units. The author sets $OQ:= \QQ \cap \ZZ_2^{\ast}, $ and one has
$\ZZ_2^{\ast}$ is the closure $\overline{OQ}$ of $OQ \subset \ZZ_2$. 
The author shows that the map $\phi$ uniquely
extends to a continuous function
 $\phi: \ZZ_2^{\ast} \smallsetminus \{-\frac{1}{3} \}  \to \ZZ_2^{\ast}$.
 % Furthermore he extends it
%to a map $\phi: \ZZ_2 \backslash \{0, -\frac{1}{3} \} \to  \ZZ_2^{\ast}.$
He shows that if $f(x) = 2x+ \frac{1}{3}$ then $f(x)$ leaves $\phi$ invariant,
in the sense that $\phi(f(x))= \phi(x)$ for all $x \in \ZZ_2^{\ast}  \smallsetminus \{-\frac{1}{3} \} .$
It follows that $f(f(x))= 4x+1$ also leaves $\phi$ invariant.

To each $x \in \ZZ_2^{\ast}   \smallsetminus \{-\frac{1}{3} \} $ he associates
 the sequence of $2$-exponents  $(p_1, p_2, ...)$ produced by iterating $\phi$.
 He proves that
an element $x \in \ZZ_2^{\ast}   \smallsetminus \{\frac{-1}{3} \} $ uniquely determine $x$;
and  that every possible sequence corresponds to some value
$x \in \ZZ_2^{\ast}   \smallsetminus \{-\frac{1}{3} \} $ 
He shows that all periodic points of $\phi$ on $\ZZ_2^{\ast}$ are rational numbers
 $x=\frac{p}{q} \in OQ$,and that there
 is a unique such periodic point for any finite sequence $(p_1, p_2, \cdots, p_m)$
 of positive integers, representing $2$-exponents, having  period $m$.
 If $ C(p_1, p_2, ..., p_{m}) = \sum_{j=0}^{m-1} 2^{p_1 + \cdots + p_j} 3^{m-1-j} $
 then this periodic pont is 
 $$
x= R(p_1, p_2..., p_m):= \frac{C(p_1, ..., p_m)}{2^{p_1 + \cdots + p_m} - 3^m}
$$
He shows that an orbit is periodic if and only if its sequence of $2$-exponents is
 periodic.

%Math Reviews summary:
% ``The aim of this article is to make a link between 
%the congruential systems investigated by Conway and
%the theory of infinite graphs. We first compare the graphs
%of congruential systems with a well-known family of infinite
%graphs: the regular graphs of finite degree considered
%by Muller and Shupp, and by Courcelle. We first consider 
%congruential systems as word rewriting systems to extract
%some subfamilies of congruential systems, the
%$q$-$p$ congruential systems, representing the regular 
%graphs of finite degree. We then prove the nonregularity of
%the Collatz graph.''

\item
Giovanni Venturini (1982), 
{\em Sul Comportamento delle Iterazioni di Alcune Funzioni Numeriche,}
Rend. Sci. Math. Institute Lombardo {\bf A 116} (1982), 115--130.
(MR 87i:11015; Zbl. 583.10009). \\
\newline
\hspace*{.25in}
The author studies functions $g(n) =  a_r n  +  b_r$
for $n  \equiv  r~ (\bmod  ~p)$ where $a_r$
$(0  \leq  r \leq  p)$ are positive rationals with denominator $p$.
He mainly treats the case that the $a_r$ take two distinct values.
If $\tau =  (a_0 a_1 \cdots a_{p-1} )^{1/p}$ has
$\tau  <  1$ then for almost all $n$ there is some $k$
with $g^{(k)} (n)  <  n$, while if $\tau  >  1$ then the iterates 
tend to increase.
[The Zentralblatt reviewer says that proofs are incomplete but 
contain an interesting idea.
Rigorous versions of these results have since been established, see
Lagarias\ (1985), Sect.\ 3.2.]

\item
Giovanni  Venturini (1989), 
{\em On the $3x  +  1$ Problem}, 
Adv. Appl. Math {\bf 10} (1989), 344--347.
(MR 90i:11020). \\
\newline
\hspace*{.25in}
This paper shows that
 for any fixed $\rho$ with $0  <  \rho  < 1$  the set of
$m \in \ZZ^+ $ which either have some $T^{(i)} (m) =  1$ or some
$T^{(i)} (m)  <  \rho m$ has density one.
This result improves on Dolan, Gilman and Manickam (1987).

\item
Giovanni  Venturini (1992), 
{\em Iterates of Number Theoretic Functions with Periodic Rational
 Coefficients (Generalization of the $3x  +  1$ problem), }
Studies in Applied Math. {\bf 86} (1992), 185--218.
(MR 93b:11102). \\
\newline
\hspace*{.25in}
This paper studies iteration of maps $g  : \ZZ  \rightarrow  \ZZ$
of the form $g(m) =  \frac{1}{d} (a_r m  +  b_r )$ if
$m  \equiv  r~ (\bmod  ~d$) for $0 \leq r \leq  d  -  1$,
where $d  \geq  2$ is an
arbitrary integer, and 
all $a_r, ~b_r  \in  \ZZ$.
These maps generalize the $(3x  +  1)$-function, and include a
wider class of such functions than in Venturini\ (1989).
The author's methods, starting in section 3,  are similar in spirit to
the Markov chain methods  introduced by Leigh (1985), which in turn
were motivated by work of 
Matthews and Watts (1984, 1985).

The author is concerned with classifying $g$-ergodic sets $S$
of such $g$ which are finite unions of
congruence classes.
For example, the mapping $g(3m) =  2m$,
$g(3m  +  1) =  4m  +  3$ and
$g(3m  +  2) =  4m  +  1$ is a
permutation and has $S =  \{ m : m   \equiv  0$ or $5 (\bmod  ~10)\}$ as a
$g$-ergodic set.
He then associates a (generally finite) Markov chain to such a
$g$-ergodic set,
whose stationary distribution is used to derive a
conjecture for the distribution of iterates
$\{ g^{(k)} (n_0 : k  \geq  1 \}$ in these
residue classes for a randomly chosen initial value $n_0$ in $S$.
For the example above the stationary distribution is
$p_0 =  \frac{1}{3}$ and $p_5 =  \frac{2}{3}$.
One can obtain also conjectured growth rates $a( g |_S )$ for
iterates of a randomly chosen initial value in a
$g$-ergodic set $S$.
For the example above one obtains
$a(g| s) =  ( \frac{4}{3} )^{2/3} { 2 \choose 3} \cong 1.0583$.

The author classifies  maps $g$ into classes $G_v (d)$, for
$v =  0,1,2 , \ldots ,$ with an additional  class $G_\infty (d)$.
The parameter $v$ measures the extent to which the 
numerators $a_r$ of iterates of $g(v)$ have common factors with $d$.
The class $G_0(d)$  consists of those maps  $g$ having 
$gcd(a_0a_1 ...a_{d-1},d)=1$, which is exactly
the {\em relatively prime case}  treated in Matthews and Watts (1984)
and the class  $G_1(d)$  are exactly those
 maps having $gcd(a_r, d^2)= gcd(a_r, d)$ for all $r$, which were  the class
 of maps treated in
Matthews and Watts (1985).
His Theorem 7  shows that for each finite $v$ 
both of Leigh's Markov chains for auxilary modulus $m=d$ are finite for all maps in $G_v(d)$ 
(strengthening Leigh's Theorem 7).   The maps  in the exceptional class
$G_\infty(d)$  sometimes, but not always,  lead to infinite Markov chains.
The class $G_\infty (d)$ presumably contains functions 
constructed by J. H. Conway
[Proc. 1972 Number Theory Conf., U. of Colorado,
Boulder 1972, 39-42] which have all $b_r =  0$ and which can encode
computationally undecidable problems.

The proof of the author's Corollary to Theorem~6 is incomplete and the
result is not established:
It remains an open problem whether a $\ZZ$-permutation having an
ergodic set $S$ with $a(g |_S )  >  1$ contains any orbit that is infinite.
Sections~6 and 7 of the paper contain many interesting examples of
Markov chains associated to such functions $g$; these examples are worth
looking at for motivation before reading the rest of the paper.

\item
Giovanni  Venturini (1997),
{\em On a Generalization of the $3x  +  1$ problem,}
Adv. Appl. Math. {\bf 19} (1997), 295--305. 
(MR 98j:11013). \\
\newline
\hspace*{.25in}
This paper considers mappings $T(x) = \frac {t_r  x  -  u_r }{p}$
when $x  \equiv  r~ (\bmod  ~p$), having
$t_r  \in \ZZ^+ $, and $u_r   \equiv  rt_r~ (\bmod  ~p$).
It shows that if gcd $(t_0 t_1 \ldots t_{p-1} , p) =  1$ and
$t_0 t_1 \ldots  t_{p-1}  <  p^p$ then for any fixed
$\rho$ with $0  <  \rho  <  1$ almost all $m  \in  \ZZ$ have an
iterate $k$ with
$| T^{(k)} (m) |  <  \rho | m |  $.
The paper also considers the question:
when are such mappings $T$ permutations of $\ZZ$?
It proves they are if and only if
$\sum_{r=0}^{p-1}\frac{1}{t_r} =  1$ and
$T(r)  \not\equiv  T(s)~ (\bmod  ~(t_r ,  t_s )$)
for
$0  \leq  r < s \leq  p  -  1$.
The geometric-harmonic mean inequality implies
that
$t_0 \ldots  t_{p-1}  >  p^p$ for such
permutations, except in the trivial case that all $t_i =  p$.

\item
Carlo Viola (1983),
{\em Un Problema di Aritmetica
(A problem of arithmetic)} (Italian),
Archimede {\bf 35} (1983), 37--39. \\
(MR 85j:11024). \\
\newline
\hspace*{.25in}
 The author states the $3x+1$ problem,
 and gives a very brief survey of known results on
 the problem, with  pointers to the literature. 
 
\item
Stanley Wagon (1985), 
{\em The Collatz problem,} 
Math. Intelligencer {\bf 7}, No. 1, (1985), 72-76. \\
(MR 86d:11103, Zbl. 566.10008). \\
\newline
\hspace*{.25in}
This article studies  a random walk imitation of the ``average'' behavior 
of the $3x  +  1$ function, computes its expected value and 
compares it to data on $3x  +  1$ iterates.

\item
Wang, Shi Tie (1988),
{\em Some researches for transformation over recursive programs}, (Chinese, 
English summary)
J. Xiamen University, Natural Science Ed. [Xiamen da xue xue bao. Zi ran ke xue ban]
{\bf 27}, No. 1 (1988), 8--12. [MR 89g:68057, Zbl. 0689.68010]\\
\newline
\hspace*{.25in}
English Abstract: "In this paper the transformations over recursive programs
with the fixpoint theory is reported. The termination condition of the duple-recursive
programs to compute the "$91$" function and to prove the $3x+1$ problem is discussed." \\

{\em Note.} The author gives a general iteration scheme for computing certain
recursively defined functions. The $91$ function $F(x)$ is defined recursively on
positive integers by the condition, if $x>100$ then $F(x) := x-10$, otherwise
$F(x) := F( F(x+11))$.

\item
Blanton C. Wiggin (1988), 
{\em Wondrous Numbers -- Conjecture about the $3n  +  1$ family,} 
J. Recreational Math. {\bf 20}, No. 2 (1988), 52--56. \\
\newline
\hspace*{.25in}
This paper calls Collatz function iterates ``Wondrous Numbers'' and
attributes this name to D. Hofstadter,
{\em G\"odel, Escher, Bach.}
He proposes studying iterates of the class of ``MU'' functions
$$
F_D(x) =  \left\{
\begin{array}{ll}
\df{x}{D} & \mbox{if}~~x \equiv 0~ (\bmod  ~D) ~, \\
~~~ \\
(D+1)x-j & \mbox{if}~~x \equiv j ~ (\bmod  ~D) ,~~ 1 \leq j \leq D-2~, \\
~~~ \\
(D+1)x+1 & \mbox{if}~~x \equiv -1~ (\bmod  ~D) ~,
\end{array}
\right.
$$
$F_D $ is the Collatz function for $D =  2$.
Wiggin's analogue of the $3x  +  1$ conjecture for a given 
$D  \geq  2$ is that all iterates of $F_D(n)$
for $n  \geq  1$ reach some number
smaller than $D$.
Somewhat surprisingly, no $D$ is known for which this is false.
It could be shown false for a given $D$ by exhibiting a cycle with all
members $ > D$; no such cycles exist for $x < 3 \times 10^4$ for 
$2 \leq  D  \leq  12$.

\item
G\"{u}nther  J. Wirsching (1993), 
{\em An Improved Estimate Concerning $3N  +  1$ Predecessor Sets,}
Acta Arithmetica, {\bf 63} (1993), 205--210.
(MR 94e:11018). \\
\newline
\hspace*{.25in}
This paper shows that, for all $a \not\equiv 0~~(\bmod~3)$, the set
$\theta_a (x) : =  \{ n \leq x:$ some $T^{(k)} (x) =  a \}$ 
has cardinality at
least $x^{.48}$, for sufficiently large $x$.
This is achieved by exploiting the inequalities of Krasikov (1989).

\item
G\"{u}nther J. Wirsching (1994),
{\em A Markov chain underlying the backward Syracuse algorithm, }
Rev. Roumaine Math. Pures Appl. {\bf 39} (1994), no. 9, 915--926.
(MR 96d:11027). \\
\newline
\hspace*{.25in}
The author constructs from the inverse iterates of the $3x+1$ function
(`backward Syracuse algorithm') a Markov chain defined on the state space
$[0,1] \times \ZZ_3^\times$,
in which $\ZZ_3^\times$ is the set of invertible 3-adic integers.
He lets $g_n (k,a)$ count the number of ``small sequence''
preimages of an element $a \in \ZZ_3^\times$
at depth $n+k$, which has $n$ odd iterates among its preimages, and
with symbol sequence 
$0^{\alpha_0} 1 0^{\alpha_1} \cdots 1 0^{\alpha_{n-1}} 1$,
with $\sum_{i=0}^n \alpha_i = k$ satisfying the ``small sequence'' condition
$0 \leq \alpha_j < 2 \cdot 3^{j-1}$ for $0 \leq j \leq n-1$.
These quantities satisfy a functional equation
$$
g_n (k,a) = \df{1}{2 \cdot 3^{n-1}} \sum_{j=0}^{2 \cdot 3^{n-1}} g_{n-1}
\left( k-j ~,~ \df{2^{j+1} a-1}{3} \right) ~.
$$
He considers the ``renormalized'' quantities
$$
\hat{g} \left(  \frac{k}{n} , a \right) : = \frac{1}{\Gamma_n} g_n (k,a)~,
$$
with
$$
\Gamma_n = 2^{1-n} 3^{- \frac{1}{2} (n-1)(n-2)} (3^n -n )~,
$$
He obtains, in
a weak limiting sense, a Markov chain whose
probability density of being at $(x, a )$ is a limiting average of
$\hat{g}(x, a )$ in a neighborhood of $(x,a)$
as the neighborhood shrinks to $(x, \alpha )$.
One step of the backward Syracuse algorithm induces (in some limiting
average sense) a limiting Markovian transition measure,
which has a
density taking the form of a product measure 
$\frac{3}{2} \chi_{\left[ \frac{x}{3} , \frac{x+2}{3} \right]} \otimes \phi$,
in which $\chi_{ \left[ \frac{x}{3} ,~ \frac{x+2}{3} \right]}$ 
is the characteristic function of the
interval $\left[ \frac{x}{3} ,~ \frac{x+2}{3} \right]$
and
$\phi$ is a nonnegative integrable function on $\ZZ_3^\times$ 
(which may take the value $+ \infty$).
\newline
\hspace*{.25in}
The results of this paper are included in Chapter~IV of Wirsching (1996b).

\item
G\"{u}nther J. Wirsching (1996), 
{\em On the Combinatorial Structure of $3N+1$ Predecessor Sets,}
Discrete Math. {\bf 148} (1996), 265--286. 
(MR 97b:11029). \\
\newline
\hspace*{.25in}
This paper studies the set $P(a)$ of inverse images of an integer $a$
under the $3x+1$ function.
Encode iterates $T^{(k)} (n) =  a$ by a set of
nonnegative integers $( \alpha_0 , \alpha_1, \ldots , \alpha_\mu )$
such that the 0-1 vector
$\bv$ encoding
$\{ T^{(j)} (n)~ (\bmod~2)$:
$0 \leq j \leq k -1 \}$
has $\bv = 0^{\alpha_0} 1~0^{\alpha_1} \cdots 1~0^{\alpha_\mu}$.
Wirsching studies iterates corresponding to ``small sequences,''
which are ones with $0 \leq \alpha_i < 2 \cdot 3^{i-1}$.
He lets $G_\mu (a)$ denote the set of ``small sequence'' 
preimages $n$ of $a$
having a fixed number $\mu$ of iterates
$T^{(j)} (n)  \equiv  1~~(\bmod~ 2)$, and shows that
$| G_\mu (a)| = 2^{\mu -1} 3^{1/2 ( \mu -1 \ldots \mu -2)}$.
He lets $g_a (k, \mu )$ denote the number of such sequences with
$k =  \alpha_0 + \alpha_1 +  \cdots  +  \alpha_\mu$, and introduces
related combinatorial quantities $\psi ( k, \mu )$ satisfying
$\sum_{0 \leq l \leq k} \psi ( l, \mu -1 )  \geq  g_a (k, \mu )$.
The quantities $\psi ( k, \mu )$ can be asymptotically estimated, and have a
(normalized) limiting distribution
$$
\psi (x)  =  \lim_{\mu  \rightarrow  \infty }
\frac{3^\mu - \mu}{2^\mu 3^{1/2\mu (\mu -1)}} 
\psi ([3^\mu - \mu ) x] , \mu ) .
$$
$\psi (x)$
is supported on $[0,1]$ and is a $C^\infty$-function.
He suggests a heuristic argument to estimate the number of ``small sequence''
preimages of $a$ smaller than $2^n a$, in terms of a double integral 
involving the
function $\psi (x)$.
\newline
\hspace*{.25in}
The results of this paper are included in Chapter~IV of Wirsching (1998).

\item
G\"{u}nther J. Wirsching (1997), 
{\em $3n+1$ Predecessor Densities and
Uniform Distribution in $\ZZ_3^\ast$,}
{\em Proc. Conference on Elementary and Analytic Number Theory}
(in honor of E.~Hlawka), Vienna, July 18-20, 1996, 
(G. Nowak and H. Schoissengeier, Eds.), 1996,
pp. 230-240. (Zbl. 883.11010). \\
\newline
\hspace*{.25in}
This paper formulates a kind of equidistribution hypothesis on
3-adic integers under backwards iteration by the $3x+1$
mapping, which, if true,
would imply that the set of integers less than $x$ which
iterate under $T$ to a fixed integer
$a \not\equiv 0$ (mod~3) has size
at least
$x^{1- \epsilon}$ as
$x \ra \In$, for any fixed $\epsilon > 0$.

\item
G\"{u}nther J. Wirsching (1998a), 
{\em The Dynamical System Generated by the $3n + 1$ Function}, 
Lecture Notes in Math. No. 1681, Springer-Verlag: Berlin 1998.
(MR 99g:11027). \\
\newline
\hspace*{.25in}
This volume is a revised version of  the author's  Habilitationsscrift
(Katholische Universit\"at Eichst\"att 1996). 
It studies the problem of showing that
for each positive integer $a \not\equiv 0 ~ ( \bmod ~3)$ a
positive proportion of integers less than $x$ iterate to $a$, as 
$x \rightarrow \infty$.
It develops an interesting 3-adic approach to this problem.
\newline
\hspace*{.25in}Chapter~I contains a history of work
on the $3x+1$ problem,
and a summary of known results.
\newline
\hspace*{.25in}Chapter~II studies the graph of iterates of $T$ on the positive
integers, (``Collatz graph'')
and, in particular studies the graph $\Pi^a ( \Gamma_T )$ connecting
all the inverse iterates $P(a)$ of a given positive integer
$a$.
This graph is a tree for any noncyclic value $a$.
Wirsching uses a special encoding of the symbolic dynamics of
paths in such trees, which
enumerates symbol sequences by keeping track of the successive blocks of 0's.
He then characterizes which graphs $\Pi^a ( \Gamma_T )$
contain a given symbol sequence reaching the root node $a$.
He derives ``counting functions'' for such sequences, and uses them to
obtain a formula giving a lower bound for the function
$$
P_T^n (a) : = \{m : 2^n a \leq m < 2^{n+1} a ~~
\mbox{and~some~iterate ~ $T^{(j)} (m) = a \}$} ~,
$$
which states that
\beql{eq1}
|P_T^n (a) | : = \{ m : 2^n (a) : = \sum_{\ell = 0}^\infty e_\ell
(n + \lfloor \log_2 \left( \df{3}{2} \right) \ell \rfloor ~,  ~~ a)~,
\eeq
where $e_\ell (k,a)$ is a ``counting function.''
(Theorem~II.4.9).
The author's hope is to use (\ref{eq1}) to prove that
\beql{eq2}
|P_T^n (a) | > c(a)2^n ~,~~ n = 1,2 , \ldots ~.
\eeq
for some constant $c(a) > 0$.
\newline
\hspace*{.25in}Chapter~III studies the counting functions appearing in the
lower bound above, in the hope of proving (2).
Wirsching observes that the ``counting functions''
$e_\ell (k,a)$, which are ostensibly defined for positive integer variables
$k , \ell , a$,
actually are well-defined when the variable $a$ is a 3-adic integer.
He makes use of the fact that one can
define a `Collatz graph' $\Pi^a ( \Gamma_T )$ in which $a$ is a
3-adic integer, by taking a suitable limiting process.
Now the right side of (\ref{eq1}) makes sense for all 3-adic integers, and
he proves that actually $s_n (a) = + \infty$
on a dense set of 3-adic integers $a$.
(This is impossible for any integer $a$ because
$|P_T^n (a) | \leq 2^n a$.)
He then proves that $s_n (a)$ is a nonnegative integrable function 
of a 3-adic
variable, and he proves that its expected value $\bar{s}_n$ is explicitly
expressible using binomial coefficients.
Standard methods of asymptotic analysis are used to estimate $\bar{s}_n$,
and to show that
\beql{eq3}
\liminf_{n \rightarrow \infty} \df{\bar{s}_n}{2^n} > 0 ~,
\eeq
which is Theorem~III.5.2.
In view of (\ref{eq1}), this says that (\ref{eq2}) ought to hold in some
``average'' sense.
He then proposes that (\ref{eq2}) holds due to the: \\

{\bf Heuristic Principle.}
{\em As $n \rightarrow \infty$, $s_n (a)$ becomes relatively close to
$\bar{s}_n$, in the weak sense that there is an absolute constant
$c_1 > 0$ such that}
\beql{eq4}
s_n (a) > c_1 \bar{s}_n ~~ for ~ all ~ n > n_0 (a) ~.
\eeq
One expects this to be true for all positive integers $a$.
A very interesting idea here is that it might conceivably be true 
for all 3-adic
integers $a$.
If so, there is some chance to rigorously prove it.
\newline
\hspace*{.25in}Chapter~IV studies this heuristic principle further 
by expressing
the counting function $e_\ell (k,a)$ in terms of simpler counting
functions $g_\ell (k,a )$ via a recursion
$$
e_\ell (k,a) = \sum_{j=0}^k p_\ell (k-j) g_\ell (j,a ) ~,
$$
in which $p_\ell (m)$ counts partitions of $m$ into parts of a 
special form.
Wirsching proves that properly scaled versions of the functions 
$g_\ell$ ``converge''
(in a rather weak sense) to a limit function which is independent of $a$,
namely
\beql{eq5}
g_\ell (k,a) \approx 2^{\ell-1} 3^{ \frac{1}{2} ( \ell^2 - 5 \ell +2)}
\psi \left( \df{k}{3^\ell} \right) ~,
\eeq
for ``most'' values of $k$ and $a$.
(The notion of convergence involves integration against test functions.)
The limit function $\psi : [0,1] \rightarrow \RR$
satisfies the functional-differential equation
\beql{eq6}
\psi'(t) = \df{9}{2} ( \psi (3t) - \psi (3t -2 )) ~.
\eeq
He observes that $\psi$ has the property of being $C^\infty$
and yet is piecewise polynomial, with infinitely many pieces, off a set of
measure zero.
There is a fairly well developed theory for related functional-differential
equations, cf. V. A. .~Rvachev,
Russian Math. Surveys
{\bf 45}
No.~1 (1990), 87--120.
Finally, Wirsching observes that if the sense of convergence in (\ref{eq5})
can be strengthened, then the heuristic principle can be deduced, and the
desired bound (\ref{eq2}) would follow.

\item
G\"{u}nther J. Wirsching (1998b),
{\em Balls in Constrained Urns and Cantor-Like Sets,}
Z. f\"ur Analysis u. Andwendungen {\bf 17} (1998), 979-996.
(MR 2000b:05007). \\
\newline
\hspace*{.25in}
This paper studies solutions to an integral equation that arises
in the author's analysis of the $3x+1$ problem (Wirsching (1998a)). Berg and
Kr\"uppel ( J. Anal. Appl. {\bf 17} (1998), 159--181) showed that the
integral equation
$$ \phi (x) = \frac {q}{q -1} \int_{qx -q +1}^{qx} \phi (y) dy $$
subject to $\phi$ being supported in the interval $[0,~ 1]$ and having
$\int_0^1 \phi(y) dy = 1$ has a unique solution
whenever $q > 1.$ The function $\phi(y)$ is a $C^{\infty}$-function.
In this paper Wirsching shows that a certain iterative procedure
converges to this solution when $q > \frac {3}{2}.$ He also shows that
for $q > 2$ the function $\phi(y)$ is piecewise polynomial off a 
Cantor-like  set of measure zero.
The case of the
$3x+1$ problem corresponds to the choice $q = 3.$

\item
Wu, Jia Bang (1992),
{\em A tendency of the proportion  of consecutive numbers of 
the same height in the Collatz problem} (Chinese), 
J. Huazhong (Central China) Univ. Sci. Technol.  [Hua zhong gong xue yuan] 
{\bf 20}, No. 5, (1992), 171--174. (MR 94b:11024, Zbl. 766.11013). \\
\newline
\hspace*{.25in}
English Abstract: "The density distribution and length of consecutive numbers of
the same height in the Collatz problem are studied. The number of integers, $K$,
which belong to $n$-tuples ($n \ge 2$) in interval $[1, 2^N)$ $(N=1, 2, ..., 24)$
is accurately calculated. It is found that the density $d(2^N)$
$(=K/2^N-1)$ of $K$ in $[1, 2^N)$ is increased with $N$. This is a correction
of Garner's inferences and prejudgments. The longest tuple in $[1, 2^{30})$, 
which is the $176$-tuple with initial number $722,067,240$ has been found.
In addition, two conjectures are proposed."

{\em Note:} The author refers to  Garner (1985). 
%From Math Reviews:
%This paper studies the density of consecutive numbers having the same height.
%Set $d(2^N ) : = \frac{1}{2^N -1} \{ n < 2^N : h (n) = h (n+1) \}$.
%It computes $d(2^N )$ for $N \leq 24$ and obsrves that $d(2^N)$ is increasing
%with $N$ for $N \leq 24$.
%The longest $k$-tuple of consecutive numbers of the same height 
%for $n \leq 2^{30}$
%is a 176-tuple with initial value $n = 722067240$.

\item
Wu, Jia Bang (1993), 
{\em On the consecutive positive integers of the same height in the 
Collatz problem} (Chinese),
Mathematica  Applicata, suppl. {\bf 6}
 [Ying yung shu hs\"{u}eh]  (1993), 150--153. (MR 1 277 568)  \\
\newline
\hspace*{.25in}
English Abstract: "In this paper, we prove that, if a $n$-tuple
($n$ consecutive positive integers of the same height) is found, an
infinite number of $n$-tuples can be found such that each of $n$ numbers
in the same tuple has the same height. With the help of a computer, the
author has checked that all positive integers up to $1.5 \times 10^8$, verified
the $3x+1$ conjecture and found a $120$-tuple,
which is the longest tuple among all checked numbers. Thus there exists an
infinite number of $120$-tuples."

{\em Note.} A table of record holders  over intervals of length $10^{7}$
is included.

\item
Wu, Jia Bang  (1995),
{\em The monotonicity of pairs of coalescence numbers in the Collatz
problem.} (Chinese),
 J. Huazong Univ. Sci. Tech. 
 [Hua zhong gong xue yuan] {\bf 23} (1995), suppl. II,
170--172. 
(MR 1 403 509) \\
\newline
\hspace*{.25in}
English Abstract: "A $310$-tuple with an initial number 6,622,073,000 has been
found. The concepts such as the pairs of coalescence numbers, conditional
pairs of coalescence numbers and unconditional pairs of coalescence
numbers are suggested. It is proved that, in interval $[1, 2^N]$, the 
density of conditional pairs of coalescence numbers 
\[
\bar{d}(2^N)= \frac{1}{2^N}\#\{n<2^N: n~ \mbox{and}~n+1~\mbox{are}~ k(\le N)~
\mbox{times~the~pair~of~coalescence~numbers}\}
\]
is increased with $N$. A conjecture that, in interval $[1,2^N]$, the density of
\[
\bar{d}_{0}(2^N)= \frac{1}{2^N}\#\{n<2^N: n~ \mbox{and}~n+1~
\mbox{pair~of~coalescence~number}\}
\]
is also increased with $N$, is proposed."\\

%[I have not seen this paper.]

\item
Masaji Yamada (1980),
{\em A convergence proof about an integral sequence},
Fibonacci Quarterly {\bf 18} (1980), 231--242.
(MR 82d:10026) \\
\newline
\hspace*{.25in}
This paper claims a proof of the $3x+1$ Conjecture.
However the  proof is faulty, with  specific mistakes pointed
out in  Math. Reviews. In particular, Lemma 7 (iii)
and Lemma 8 are false as stated.

\item
Yang, Zhao Hua  (1998),
{\em An equivalent set for the $3x+1$ conjecture} (Chinese),
J. South China Normal Univ. Natur. Sci. Ed. 
[ Hua nan shi fan da xue xue bao. Zi ran ke xue ban] 1998, no. 2, 66--68.
(MR 2001f:11040). \\
\newline
\hspace*{.25in}
%MR
 The paper shows that the $3x+1$ Conjecture is true if,
for any fixed $k \ge 1$, it is true for
all positive integers in the arithmetic progression 
$n \equiv 3 + \frac{10}{3}(4^k - 1) ~(\bmod~2^{2k+2}).$
Note that Korec and Znam (1987) gave analogous conditions 
for the $3x+1$ Conjecture being true, using
a residue class $(\bmod~p^k)$ where $p$ is an odd prime for
which $2$ is a primitive root. 

\item
Yang, Zhi and Zhang, Zhongfu (1988),
{\em Kakutani conjecture and graph theory representation of the black hole problem}, 
Nature Magazine [Zi ran za zhi (Shanghai)]  {\bf 11} (1988), No. 6, 453-456. \\
\newline
\hspace*{.25in}
This paper considers the Collatz function (denoted $J(n)$)  viewing iteration of the
function as a directed graph, with edges from $n \to J(n)$. The set $N_J$ of Kakutani
numbers is the set of all numbers that iterate to $1$; these form a connected graph.  
It phrases the $3x+1$ conjecture
as asserting that all  positive integers are ``Kakutani numbers." 

 It also discusses the ``black hole" problem, which concerns iteration on  $n$-digit 
 integers (in base $10$), taking $K(n) = n^{+} - n^{-}$ where $n^{+}$ is the $n$-digit
 integer obtained from $n$, arranging digits in decreasing order,
 and $n^{-}$ arranging digits in increasing order. This iteration can
 be arranged in a  directed graph with edges from $n$ to $T(n)$. It
 is known that  for $3$-digit numbers, $T(\cdot)$ iterates to the fixed point $n=495$,
 and for $4$-digit numbers it iterates to the fixed point $n=6174$
 (the "Kaprekar constant").  For
 more than $4$-digits, iterations may enter cycles of period exceeding
 one, called  here ``black holes".
 For $n=6$ there is a cycle of length $7$ with
 $n=840852$, and for $n=9$ there is a cycle of length $14$ starting with $n=864197532$. 
The paper is descriptive, with examples but no proofs.

{\em Note.} The ``black hole" problem has
a long history. The starting point was: D. R. Kaprekar, {\em An interesting property of the
number 6174}, Scripa Math. {\bf 21} (1955) 244-245.
Relevant papers include H. Hasse and G. D. Prichett,
{\em The determination of all four-digit Kaprekar constants}, J. reine Angew. Math. {\bf 299/300}
(1978), 113--124; G. Prichett, A. Ludington and J. F. Lapenta, {\em The determination of
all decadic Kaprekar constants}, Fibonacci Quarterly {\bf 19} (1981), 45--52. 
The latter paper proves that  ``black holes"  exist for every $n \ge 5$.

\item
Zhang, Cheng Yu  (1990),
 {\em A generalization of Kakutani's conjecture.} (Chinese)
Nature Magazine [Zi ran za zhi (Shanghai)] , {\bf 13} No. 5, (1990), 267--269. \\
\newline
\hspace*{.25in}
Let $p_1=2, p_2=3, ...$ list the primes in increasing order. 
This paper suggests the analysis of the mappings $T_{d,k}(x)$, indexed by $d \ge 1$
and $k \ge 0$, with
$$
T_{d,k}(x) = \left\{
\begin{array}{cl}
\df{x}{2}  & \mbox{if} ~~ 2 | x ~~~ \\
~~~ & ~~~\\
\df{x}{3} & \mbox{if}~~ 3 | x~\mbox{and}~2 \nmid x \\
~~~& ~~~\\
... & \\
~~~& ~~~\\
\df{x}{p_{d-1}} & \mbox{if}  ~ p_{d}|x~\mbox{and}~ p_j \nmid x, ~1 \le j \le d-1, \\
~~~& ~~~\\
p_{d+1}x + (p_{d-1})^k & \mbox{otherwise}.
\end{array}
\right.
$$
This function for  $(d,k)=(1,0)$ is the Collatz  function.  
The author actually does not give a tie-breaking rule for defining the
function when several $p_j$ divide $x$. However, this is not
important in studying  the conjecture below. 

The author formulates the {\em Kakutani $(d, k)$-Conjecture}, which states that the
function $T_{d,k}(x)$ iterated on the positive integers always reaches the
integer $(p_{d+1})^k$. Note  that $(p_{d+1})^k$ always belongs to 
a cycle of the function $T_{d,k}(x)$.
This conjecture for $(1,k)$ was presented
earlier by the author,
%in ZHANG, Cheng Yu (1987), 
where it was shown equivalent to the $3x+1$ Conjecture. 

The author shows this conjecture is false for $d=4,5, 6, 8$ because
there is a nontrivial cycle. The cycles he finds, for given $p_{d+1}$ and  starting
value $n$ are:  $p_5=11 ~(n=17), ~p_6=13 ~(n=19),~ p_7=17~ (n=243),
~p_9=29 ~(n=179).$

The author then modifies the Kakutani $(d,k)$-Conjecture in these cases to say that every
orbit on the positive integers enters one of the known cycles. 
For example he modifies Conjecture $(4,k)$ to say that all orbits
reach the value $(11)^k$ or else enter the cycle containing $ 17 \cdot (11)^k$. 
He proves that the  Conjecture $(d, k)$ is equivalent to the
modified Conjecture $(d, 0)$ for $d=1, 2, 3, 4.$

%\item
%Zhang, Zhongfu and Yang, Zhi (1988),
%[unknown title]
%Nature Magazine [Zi ran za zhi (Shanghai)] , {\bf 11} (1988), 453. \\
%\newline
%\hspace*{.25in}
%[I have not seen this paper.]

\item
Zhang, Zhongfu  and Yang, Shiming (1998),
{\em Problems on mapping sequences} (Chinese),
Mathmedia [Shu xue chuan bo ji kan] {\bf 22} (1998), No. 2, 76--88.\\
\newline
\hspace*{.25in}
This paper presents  basic results on a number of iteration problems
on integers, represented as directed graphs, and lists twelve open
problems.  This includes a discussion of 
the $3x+1$ problem. It also discusses iteration of the map $V(n)$
which equals  the difference of numbers
given by the decimal digits of $n$ rearranged in decreasing, resp.
increasing order. On four digit numbers $V(n)$ iterates to a fixed
point $6174$, the Kaprekar constant. \\

\item
Zhou, Chuan Zhong  (1995),
{\em Some discussion on the $3x+1$ problem} (Chinese),
J. South China Normal Univ. Natur. Sci. Ed. 
[ Hua nan shi fan da xue xue bao. Zi ran ke xue ban] {\bf 1995}, No. 3, 103--105.
(MR 97h:11021). \\
\newline
\hspace*{.25in}
  Let $H$ be
the set of positive integers that eventually reach $1$ under iteration
by the Collatz function $C(n)$; the $3x+1$ conjecture states that  $H$ consists
of all positive integers.  This paper  extends a result of B. Y. Hong
[J. of Hubei Normal University. Natural Science
Edition (Hubei shi fan xue yuan xue bao. Zi ran ke xue) 1986, no. 1, 1--5.]
by showing that if the $3x+1$ Conjecture is not true
the minimal positive $n \not\in H$ that does not iterate to $1$ 
%$\min \{ n:~n \not\in H\}$ 
satisfies 
$n \equiv 7, 15, 27, 31, 39, 63, 79$ or $91~(\bmod ~96).$ It also proves
that if $2^{2m+1} - 1 \in H$ then $2^{2m+2} - 1 \in H.$ Finally the paper observes 
 that the numerical computation of $C(n)$ could be simplified
if it were performed in base $3$.

\item
Zhou, Chuan Zhong (1997),
{\em Some recurrence relations connected with the $3x+1$ problem} (Chinese),
J. South China Normal Univ.  Natur. Sci. Ed.
[Hua nan shi fan da xue xue bao. Zi ran ke xue ban]  {\bf 1997}, no. 4, 7--8. \\
\newline
\hspace*{.25in}
English Abstract:  ``Let $\NN$ denote the set of natural numbers, $J$ denote
the $3x+1$ operator, and $H=\{ n \in \NN: \mbox{there~is~k} \in \NN
~\mbox{so~that}~ J^{k}(n)=1\}.$ The conjecture $H=\NN$ is the so-called
$3x+1$ problem. In this paper, some recurrence relations on this
problems are given."

%[I have not seen this paper.]
\end{enumerate}

\noindent
{\bf Acknowledgements.}
I am indebted to B. Bailey, A. M. Blas, M. Chamberland, G.-G. Gao, 
J. Goodwin, C. Hewish,
S. Kohl, Wang Liang, D. Merlini, P. Michel, C. Reiter, T. Urata, 
G. Venturini, S. Volkov and 
G. J. Wirsching for references. I am indebted to the staff at Vinculum
for supplying the Coxeter reference. For help with Chinese references I thank
W.-C. Winnie Li,  Xinyun Sun, Yang Wang and Chris Xiu. For help with Japanese references
I thank T. Ohsawa. \\

%****************************************************
%
% Section 4. References
%
%****************************************************
%
%\section{References}
%\hsp

%\noindent
%{\bf Remark.}
%I would appreciate knowing of any other papers on the $3x +  1$ problem.
%Please write or send email to: {\tt lagarias@umich.edu}

\end{document}